\documentclass{gtpart}

\usepackage{hyperref}
\hypersetup{breaklinks=true}

\numberwithin{equation}{section}
\theoremstyle{plain}
\newtheorem{theorem}[equation]{Theorem}

\newtheorem{lemma}[equation]{Lemma}
\newtheorem{corollary}[equation]{Corollary}
\newtheorem{conjecture}[equation]{Conjecture}

\theoremstyle{remark}

\theoremstyle{definition}
\newtheorem{definition}[equation]{Definition}

\newtheorem*{question*}{Question}

\usepackage{datetime} 

\usepackage{graphicx}

\usepackage{url}

\title[Quantitative Bi-Lipschitz embeddings]{Quantitative Bi-Lipschitz embeddings of bounded curvature manifolds and orbifolds}

\author{Sylvester Eriksson-Bique}
\givenname{Sylvester}
\surname{Smith}
\address{Courant Institute \\ New York University \\ \newline 251 Mercer St \\ \newline New York, NY-10012 \\USA}
\email{ebs@cims.nyu.edu}
\urladdr{http://www.cims.nyu.edu/~ebs}

\arxivreference{1507.08211 }
\arxivpassword{8qupe}



\newcommand{\QED}{
\begin{flushright}
$\square$
~\\ ~\\
\end{flushright}}

\newcounter{prob}
\newcommand{\N}{\ensuremath{\mathbb{N}}}

\newcommand{\mc}[1]{\ensuremath{\mathfrak{#1}}}

\newcommand{\diam}{\ensuremath{\mathrm{diam}}}

\def\XXint#1#2#3{{\setbox0=\hbox{$#1{#2#3}{\int}$ }
\vcenter{\hbox{$#2#3$ }}\kern-.58\wd0}}

\def\co{\colon\thinspace}

\begin{document}

\begin{abstract}
We construct bi-Lipschitz embeddings into Euclidean space for bounded diameter subsets of manifolds and orbifolds of bounded  curvature. The distortion and dimension of such embeddings is bounded by diameter, curvature and dimension alone. We also construct global bi-Lipschitz embeddings for spaces of the form $\R^n / \Gamma$, where $\Gamma$ is a discrete group acting properly discontinuously and by isometries on $\R^n$. This generalizes results of Naor and Khot. Our approach is based on analyzing the structure of a bounded curvature manifold at various scales by specializing methods from collapsing theory to a certain class of model spaces. In the process we develop tools to prove collapsing theory results using algebraic techniques.
\end{abstract}

\maketitle

\section{Introduction}

\subsection{Statement of results}

It is an old problem to describe, in some insightful manner, the metric spaces which admit a bi-Lipschitz embedding into Euclidean space \cite{Heinonen2000}. Recall the definition of a Bi-Lipschitz map.

\begin{definition} Let $L>0$. A map  $f\co X\to Y$ between metric spaces is called a \textit{(L-)bi-Lipschitz} function (also embedding) if

\begin{equation}
\frac{1}{L} d(a,b) \leq d(f(a),f(b)) \leq L d(a,b) \label{eq:bilip}
\end{equation}
 for all $a,b \in X$\footnote{Where ambiguity does not arise we will use the letter $d$ for the distance on any metric space that arises in the course of this paper.}.  The smallest constant $L$ satisfying the previous equation  is referred to as \emph{the distortion} of $f$. By a bi-Lipschitz function of distortion, or with distortion, $L$, we mean a function whose distortion can be bounded by $L$.
\end{definition}

We are generally not concerned with the optimal distortions of a bi-Lipschitz map, but simply some upper bounds. When $Y=\R^N$, a necessary condition for embeddability is that $(X,d)$ is metric doubling (see Definition \ref{doublingdef}). As a partial converse, Assouad proved that any doubling metric space admits an $\alpha-$bi-H\"older embedding into Euclidean space for any $0<\alpha<1$ \cite{assouad} (also see \cite{naorassouad}). But due to counter-examples by  Pansu \cite{pansu89} (as observed by Semmes \cite{semmes1996nonexistence})  and Laakso \cite{Laakso2000} additional assumptions are necessary for bi-Lipschitz embeddability. 

Most of the spaces considered in this paper admit some bi-Lipschitz embedding but our interest is in proving uniform, or quantitative, embedding theorems where the distortion is bounded in terms of the parameters defining the space. While a number of authors have addressed the question of quantitative bi-Lipschitz embeddability, natural classes of spaces remain for which such embeddings haven't been constructed.

Inspired by a question originally posed by Pete Storm we focus on embedding certain classes of manifolds and orbifolds. In our discussion, we will always view Riemannian manifolds (and orbifolds) as metric spaces by equipping them with the Riemannian distance function.  We emphasize that the problem of finding \textit{quantitative} bi-Lipschitz embeddings of Riemannian manifolds into Euclidean spaces has little to do with classical results such as the Nash isometric embedding theorem for Riemannian manifolds. The notion of isometric embedding in that context only requires the map to  preserve distances infinitesimally, whereas a bi-Lipschitz map controls the distances at any scale up to a factor. A Nash embedding for a \textit{compact} manifold/subset will, however, result in a bi-Lipschitz embedding with a possibly very large distortion. To give a bound for the distortion of our embedding we will need to resort to very different techniques. We also note, that Nash embeddings don't directly give bi-Lipschitz embeddings for non-compact spaces, while some of our results apply for such examples as well.

\paragraph{Example:} For intuition, consider the family of surfaces
$$S_N=\{x^2 + y^2 + N^2z^2 = 1\}.$$
Denote by $g_{N}$ the restricted metric on this two-dimensional surface. This defines boundaries of ellipsoids in $\R^3$ which are isometrically embedded in $\R^3$ but the natural embedding is not bi-Lipschitz. The distance along the surface between the points $x_\pm=\left(0,0, \pm\frac{1}{N} \right)$ is $d(x_+, x_-) \sim 2$, but in $\R^3$ the straight-line distance is $|x_+-x_-| = \frac{2}{N}$. Thus, the isometric embedding of $S_N$ does not have a well-controlled bi-Lipschitz constant. However, if we define the map $L \co S_N \to \R^3$ by sending $(x,y,z) \to \left(x,y, z + \left(1-\sqrt{x^2+ y^2}\right)H(z)\right)$. Here $H(z) = 1$ for $z \geq 0$ and otherwise $H(z)=0$. Such a map can be shown to be $C$-bi-Lipschitz, but it is clearly not an isometry. This constant $C$ can be chosen independent of $N$.  

Manifolds comprise a large class of metric spaces and we need to place some assumptions in order to ensure uniform embeddability. In order to ensure doubling it is natural to assume a diameter bound as well as a lower Ricci-curvature bound. Our results require the somewhat stronger assumption of bounding the sectional curvature in absolute value. Thus, we are led to the following theorem.

\begin{theorem} \label{thm:bilip} Let $D>0$. Every bounded subset $A$ with $\diam(A)\leq D$ in a $n$-dimensional complete Riemannian manifold $M^n$ with sectional curvature $|K| \leq 1$ admits a bi-Lipschitz embedding $f\co A \to \R^N$ with distortion less than $C(D,n)$ and dimension of the image $N \leq N(D,n)$. 
\end{theorem}

We emphasize that there are no assumptions on the injectivity radius or lower volume bound. If such an assumption were placed, the result would follow directly from Cheeger-Gromov compactness or a straightforward doubling argument (see Lemma \ref{doublingarg} and \cite{greene1988lipschitz}). While somewhat unexpected, we are also able to prove a version of the previous theorem for orbifolds.

\begin{theorem} \label{thm:orbifoldmain} Let $D>0$. Every bounded subset $A$ with $\diam(A)\leq D$ in a $n$-dimensional complete Riemannian orbifold $O^n$ with sectional curvature $|K| \leq 1$ admits a bi-Lipschitz embedding $f\co A \to \R^N$ with distortion less than $C(D,n)$ and dimension of the image $N \leq N(D,n)$. 
\end{theorem}

In particular, we have the following non-trivial result.

\begin{theorem}\label{thm:flatellipt} Every connected complete flat and elliptic orbifold $O$ of dimension $n$ admits a bi-Lipschitz mapping $f\co O \to \R^N$ with distortion less than $D(n)$ and dimension of the image $N \leq N(n)$. The constants $D(n)$ and $N(n)$ depend only on the dimension.
\end{theorem}

\paragraph{Remark:} We can give explicit bounds for the case of flat and elliptic orbifolds, but not for arbitrary bounded curvature orbifolds. This is because the proof of Fukaya's fibration theorem in \cite{chgrfu} uses compactness in a few steps and thus does not give an explicit bound for the parameters $\epsilon(n),\rho(n,\eta,\epsilon)$. Some related work in \cite{ghanaat}, when combined with techniques from \cite{chgrfu}, is likely to give an explicit bound, but to our knowledge this has not been published. In any case, the bounds extracted by known means would grow extremely rapidly (see the bounds in \cite{ruh1982, buser1981gromov}). For flat orbifolds and manifolds we attain a bound on the distortion $D(n)$ which is of the order $O(e^{Cn^4 \ln(n)})$.

For mostly technical reasons we prove an embedding result for certain classes of ``model'' spaces which occur in our proofs. These involve the notion of a quasiflat space which is a slight generalization of a complete flat manifold to the context of quotients of nilpotent Lie groups. They correspond to model spaces for collapsing phenomena. The precise definition is in Definition \ref{quasiflatdef}.

\begin{theorem}\label{thm:vector}For sufficiently small $\epsilon(n)$, every $\epsilon(n)$-quasiflat $n$-dimensional orbifold $M$ admits a bi-Lipschitz embedding into $\R^N$ with distortion and dimension $N$ depending only on $n$. Further, every locally flat Riemannian orbivector bundle over such a base with its natural metric admits such an embedding.  \end{theorem}

It may be of interest, that some of these spaces are not compact and do not have non-negative sectional curvature. Yet, they admit global bi-Lipschitz embeddings. This theorem could also be stated for orbifolds, and the proof below will directly apply to them. For terminology related to orbifolds we recommend consulting for example \cite{thurston2002geometry, ratcliffehyper,ding2011restriction,Kleiner2014}. Some of this terminology is covered in the appendix.

We remark that the results above can trivially be generalized to Gromov-Hausdorff limits of bounded curvature orbifolds, which themselves may not be orbifolds. The Gromov-Hausdorff limits of bounded curvature Riemannian manifolds have been intrinsically described using weak Alexandrov-type curvature bounds by Nikolaev in \cite{nikolaev1991}, and as such our results also lead to embeddings of such spaces.

\subsection{Outline of method}

Bi-Lipschitz embedding problems can be divided to two subproblems: embedding locally at a fixed scale and embedding globally. Similar schemes for constructing embeddings have been employed elsewhere, such as in \cite{langplaut,naormarkov,seo2011char}. Constructing global embeddings may be difficult, but by assuming a diameter bound and lower curvature bound we can reduce it to embedding a certain fine $\epsilon$-net. 

We are thus reduced to a local embedding problem at a definite scale. In particular, since we don't assume any lower volume bound, we need to apply collapsing theory in order to establish a description of manifolds at small but fixed scales. For general manifolds we first use a result from Fukaya \cite{fukaya90vect} (see also \cite{chgrfu}). Recall, that if $g$ is a metric tensor of a Riemannian manifold, and $T$ is any tensor, then its norm with respect to the metric is denoted by $||T||_g$. There are different tensor norms, but all of them are quantitatively equivalent. To fix a convention, $||T||_g$ will refer to the operator norm.

\begin{theorem}(Fukaya, \cite{chgrfu,fukayaICM,rong1996fundamental})\label{thm:fukayamanifold} Let $(M,g)$ be a complete Riemannian manifold of dimension $n$ and sectional curvature $|K|\leq 1$, and $\epsilon>0$ be an arbitrary constant. For every $\eta > 0$, there exists a universal $\rho(n,\eta,\epsilon)>0$ such that for  any point $p \in M$ there exists a metric $g'$ on the ball $B_p(10\rho(n,\eta,\epsilon))$ and a complete manifold $M'$ with the following properties.

\begin{itemize} 
\item $$||g-g'||_{g}<\eta$$ 
\item  $(B_p(5\rho(n,\eta,\epsilon)),g')$ is isometric\footnote{By isometry we mean Riemannian isometry in the sense that the metric tensors are transformed to each other via the push-forward of the differential. What we really use is that the ball $(B_p(\rho(n,\eta,\epsilon)),g')$ will have a distance preserving mapping to a subset of $M'$.} to a subset of a complete Riemannian manifold $M'$. 
\item $M'$ is either a $\epsilon$-quasi flat manifold or a locally flat vector bundle over a $\epsilon$-quasiflat manifold.
\end{itemize}
\end{theorem}

The constant $\epsilon$ will always be very small. We will set $\epsilon=\epsilon(n)$ for some dimension dependent $\epsilon(n)$, which is fixed below in section \ref{prelimin}. This choice guarantees that the model spaces are sufficiently nice. In other words their collapsing can be studied algebraically. Throughout the paper the constant $\epsilon(n)$ will have the same value, which depends on $n$. The constant $\eta$ will be fixed as $1/4$, and it is used to guarantee that $(B_p(\rho(n,\eta,\epsilon)),g')$  is bi-Lipschitz to $(B_p(\rho(n,\eta,\epsilon),g)$. Since $\epsilon$ and $\eta$ are chosen fixed, dependent on dimension, the scale $\delta=\rho(n,\eta,\epsilon)$ at which we have a nice model space, depends only on the dimension.

This theorem isn't stated as such in the references. Our terminology is different and is introduced in Definitions \ref{quasiflatdef} and \ref{quasiflatvect}. In \cite{fukayaICM} a topological version is stated and in \cite{chgrfu} the main results concern a fibration structure. 

The claims of our theorem are contained in a similar form in the appendix of \cite[Appendix 1]{chgrfu}, and the proof of our version can be derived from it. The main difference is that in \cite{chgrfu} the role of $M'$ is played by the normal bundle of a singular orbit (of the $N$-structure), $\nu(\mathcal{O}_q)$, with a natural metric. Here we have merely explicated that such a normal bundle can be described as a vector bundle over a $\epsilon(n)$-quasiflat manifold. In the appendix we outline a proof of the analogous statement for orbifolds. We also use the results in \cite{rong1996fundamental} that guarantee us a dimension independent bound on the sectional curvatures of the model spaces. 

Fukaya's fibration theorem permits us to reduce the embedding problem to embedding certain vector bundles (``model spaces''). These spaces have the ``twisted product'' structure of a vector bundle over a quasiflat base. Further, these spaces also have the advantage of possessing an algebraic description. Thus, their metric geometry can be studied using the additional structure. 

We next give a rough description of our approach to embedding the vector bundles $M$ that arise. We first use an approximation argument from Lemma \ref{approxvect} and Lemma \ref{gromhaus} to reduce the problem to embedding locally flat Riemannian vector bundles over flat manifolds. Next, the flat manifold $M$ is decomposed into pieces, and embeddings for each piece are patched together using Lipschitz extension theorems. These decompositions are similar in spirit but were developed independently of arguments by Seo and Romney in \cite{romneymatt,seo2011char}.

To be more explicit, assume that the diameter of the base of the vector bundle is unity. Let us assume $L=V(S)$ is a locally flat Riemannian vector bundle over a flat base $S$. The space $S$ is also the zero-section of the vector bundle $L$ equipped with the restricted metric. Define for $i \geq 1$ the sets $T_i$ to be the points at distance roughly $2^i$ from the zero section. Additionally, define $T_0$ to be the points  $x \in L$ with $d(x,S)\leq 4$. By using the radial function $r(x)=d(x, S)$, we are able to subdivide the embedding problem for the entire space to that of embedding each $T_i$ individually. 

To embed each $T_i$, we re-scale the space by a definite amount to reveal a simpler vector bundle structure at a small, but definite scale. This is done using an algebraic collapsing theory argument.  The difficult aspect of the proof is showing that these local descriptions exist at a scale comparable to the diameter of $T_i$, which is necessary for applying a doubling Lemma \ref{doublingarg} \footnote{If $T_i$ is scaled down to unit size, the curvature bound is scaled up by $2^{2i}$, and thus applying collapsing theory is difficult. If one only wants to prove Theorem \ref{thm:bilip}, then this scaling up could be avoided. In fact, the algebraic arguments below can be avoided by such an approach, but this would lead to other technicalities, as well as weaker results.}. Further, the nature in which the scaling produces a simpler space is somewhat delicate. For $i =0$, the resulting spaces are of the same dimension, but resemble vector bundles over lower dimensional spaces, i.e. their collapsing occurs more uniformly. For $i \geq 1$, the local descriptions split into a product of an interval and a lower dimensional space. In either case, by induction we can assume that the simpler spaces embed, and apply the doubling Lemma \ref{doublingarg} to combine said embeddings to an embedding of $T_i$.

The algebraic collapsing theory argument involves quantitative versions of the local group arguments applied in collapsing theory \cite{fukaya90vect}. Many of these arguments were inspired by proofs in \cite{ratcliffehyper, thurston1997three}. The core insight is that, up to a finite index, the fundamental groups are lattices. Lattices have ``nice'' bases, and the lengths of these basis elements give essentially canonically defined collapsing scales. These collapsing scales can be used to identify invariant submanifolds in the space. 

We also construct embeddings for bounded curvature orbifolds. For such spaces we need a version of Fukaya's theorem for orbifolds. The proof of this result is essentially contained in \cite[Appendix 1]{chgrfu}, and \cite{ding2011restriction}, but we provide a rough outline in the appendix.

\begin{theorem}\label{thm:fukayaorbifold} Let $(O,g)$ be a complete Riemannian orbifold of dimension $n$ and sectional curvature $|K|\leq 1$, and $\epsilon>0$ be an arbitrary constant. For every $\eta > 0$ there exists a universal $\rho(n,\eta,\epsilon)>0$ such that for  any point $p \in O$ there exists a metric $g'$ on the ball $B_p(10\rho(n,\eta,\epsilon))$ and a complete Riemannian orbifold $O'$ with the following properties.

\begin{itemize} 
\item $$||g-g'||_{g}<\eta$$ 
\item  $(B_p(5\rho(n,\eta,\epsilon)),g')$ is isometric to a subset of $O'$.
\item $O'$ is either a $\epsilon$-quasiflat orbifold or a locally flat Riemannian orbivector bundle over a $\epsilon$-quasiflat orbifold $S$.
\end{itemize}
\end{theorem}
The same comments on the choice of these parameters apply as before. When we apply this theorem, we fix $\epsilon=\epsilon(n)$ and $\eta=1/4$, which allows us to use that the scale of the model spaces $\rho(n,1/4,\epsilon(n))$ only depends on the dimension.

The methods to embed bounded curvature orbifolds are essentially the same as for the manifold case, except for modifying terminology and adding some additional cases.

\subsection{Previous work}

Bi-Lipschitz embedding problems have been extensively studied for finite metric spaces, see \cite{bourgain85,indyk2004,matsek90}. The problem has usually been to study the asymptotics of distortion for embedding certain finite metric spaces, either in general or restricting to certain classes of metrics (e.g., the earth mover distance). The image is either a finite dimensional $l_p^n$ or an infinite dimensional Banach space such as $L^1$. Motivation for such embeddings stems, for example, from approximation algorithms relating to approximate nearest neighbor searches or querying distances. There has also been some study on bi-Lipschitz invariants such as Markov type, and Lipschitz extendability, which can be applied once a bi-Lipschitz map is constructed to a simple space \cite{naorintro}.  

For non-finite metric spaces the research is somewhat more limited. We mention only some of the contexts within which embedding problems have been studied: subsets of $\R^n$ with intrinsic metrics by Tatiana Toro \cite{toro95,toro94}, ultrametric spaces \cite{Luosto1996, Luukkainen1994}, certain weighted Euclidean spaces by David and Semmes \cite{semmes1993bilipschitz, semmes1996nonexistence}; and an embedding result for certain geodesic doubling spaces with bicombings \cite{langplaut} by Lang and Plaut. Some of these works concern bi-Lipschitz parametrizations, which is a stronger problem. Concurrently with this paper Matthew Romney \cite{romneymatt} has considered related embedding problems on Grushin type spaces and applied methods from \cite{seo2011char} as well as the current paper. The different abstract question of existence of any bi-Lipschitz parametrization for non-smooth manifolds has been discussed in \cite{heinonen2011flat}. In contrast to these previous results, our goal is to prove embeddability theorems for natural classes of smooth manifolds, with uniformly bounded distortion and target dimension.

The closest results to this work are from Naor and Khot, who construct embeddings for flat tori and estimate the worst possible distortion \cite{naorfourier}. Improvements to their results were obtained in \cite{haviv2013euclidean}. In comparison, our work provides embeddings for any compact or non-compact complete flat orbifold. Our bounds are arguably weaker due to the higher generality.  Also, noteworthy is the paper of Bonk and Lang \cite{bonk2003bi}, where they prove a bi-Lipschitz result for Alexandrov-surfaces with bounded integral curvature. Related questions on bi-Lipschitz embeddability have been discussed in \cite{andoni2015snowflake}, where one considers an Alexandrov space target. 

A number of spaces can be shown not to admit any bi-Lipschitz embedding into Euclidean space. In addition to classical examples such as expanders, nontrivial examples were found by Pansu \cite{pansu89}  and Laakso \cite{Laakso2000}. These examples are related to a large class of examples that are covered by Lipschitz differentiation theory, which was initially developed by Cheeger \cite{ChDiff99}. 

\subsection{Open problems}

A natural question is to further study the dependence of our results on the curvature bounds assumed. At first one might ask whether the upper curvature bound is necessary. We conjecture that it may be dropped. In fact, in light of our results on orbifolds, it seems reasonable to presume that Alexandrov spaces admit such an embedding. For the definition of an Alexandrov space we refer to \cite{BuragoGeom, BGP}. The conjecture is also interesting within the context of Riemannian manifolds.

\begin{conjecture} \label{conj:alex} Let $D>0$. Every bounded subset $A$ with $\diam(A)\leq D$ in a $n$-dimensional complete Alexandrov space $X^n$ with curvature $K \geq -1$ admits a bi-Lipschitz embedding $f\co A \to \R^N$ with distortion less than $C(D,n)$ and dimension of the image $N \leq N(D,n)$. 
\end{conjecture}

If we assume in addition that the space $X^n$ is volume non-collapsed, the theorem follows from an argument in \cite{alexbook}. One of the main obstacles in proving embedding theorems for these spaces is the lack of theorems in the Lipschitz category for Alexandrov spaces. Most notably, proving a version of Perelman's stability theorem would help in constructing embeddings \cite{kapovitch2007perelman}. Perelman claimed a proof without publishing it and thus it is appropriately referred to as a conjecture. See Definition \ref{def:grhaus} for the definition of the Gromov-Hausdorff-distance $d_{GH}$.

\begin{conjecture} \label{conj:perelman} Let $X^n$ be a fixed $n$-dimensional compact Alexandrov space with curvature $K
\geq -1$. Then there exists an $0<\epsilon_0$ (depending possibly on $X^n$) such that for any other Alexandrov space $Y^n$ with 

$$d_{GH}(X^n,Y^n) < \epsilon_0,$$
we have a bi-Lipschitz map $f\co X^n \to Y^n$ with distortion at most $L$ (which may depend on $X^n$).
\end{conjecture}

Both of these conjectures seem hard. But mostly for lack of counter examples, we also suggest that the lower sectional curvature bound could be weakened to a lower Ricci-curvature bound. Weakening the curvature assumption to a simple Ricci-curvature bound results in great difficulty in controlling collapsing phenomena. However, even without collapsing the conjecture is interesting and remains open. Thus, we state the following conjecture.

\begin{conjecture} \label{conj:ricci} Let $A$ be a subset with $\diam(A)\leq D$ in a $n$-dimensional complete Riemannian manifold $(M^n,g)$ with curvature $\mathrm{Ric}(g) \geq -(n-1)g$. Assume $\text{Vol}(B_1(p))>v$ for some $p \in A$. Then there exists a bi-Lipschitz embedding $f\co A \to \R^N$ with distortion less than $C(D,n,v)$ and dimension of the image $N \leq N(D,n,v)$. 
\end{conjecture}

Finally, we remark on the problem of optimal bounds for our embeddings in the Theorem for flat and elliptic orbifolds \ref{thm:flatellipt}. As remarked, we obtain a bound of $D(n)$ of the order $O(e^{Cn^4\ln(n)})$. The main source of distortion is the repeated and inefficient use of doubling arguments at various scales, which result in multiplicative increases in distortion. In a related paper, Regev and Haviv  improve the super-exponential upper bound from \cite{naorfourier} and obtain $O(n\sqrt{\log(n)})$ distortion for $n$-dimensional flat tori \cite{haviv2013euclidean}. Thus, it seems reasonable to suspect that the true growth rate of $D(n)$ is polynomial. As pointed out to us by Assaf Naor, this problem is also related to \cite{andoni2015snowflake} because certain finite approximations to Wasserstein spaces arise as quotients of permutation groups.

\subsection{Outline}

In the next section we give explicit embeddings for three types of bounded curvature spaces. These examples illustrate the methods used to prove the embedding results. We will not explain all the details, as some of them are presented well in other references and will become more apparent in the course of the proof of the main theorem. Following this, we collect some general tools and lemmas that will be used frequently in the proofs of the main results. Some of these are very similar to \cite{langplaut,naormarkov}. Finally, in the fourth section we give full proofs of the main embedding theorems. This section proceeds by increasing generalities. First, flat manifolds are embedded, followed by flat orbifolds, and quasiflat orbifolds. Ultimately the results are applied to Riemannian manifolds and orbifolds. The appendix collects a few of the most technical results on quotients of nilpotent Lie groups and collapsed orbifolds. \\

\noindent \textbf{Acknowledgments:} The author thanks his adviser Bruce Kleiner for suggesting the problem and for numerous discussions on the topic. Discussions and comments from Jeff Cheeger, Zahra Sinaei, Or Hershkovits, Matthew Romney and Tatiana Toro have also been tremendously useful. We are grateful for Anthony Cardillo for drawing Figure \ref{fig:normals1}. We also thank the referee for many detailed comments that helped us improve the paper tremendously, as well as the reviewing committee for the patience in going through different versions of the paper. This research was supported by a NSF grant DGE 1342536.

\section{Embedding some key examples} \label{examples}

In the next section we give nearly explicit embeddings for three types of bounded curvature spaces. These examples illustrate the methods used to prove the embedding results. We will not explain all the details, as some of them are presented well in other references and will become more apparent in the course of the proof of the main theorem.

\subsection{Lens spaces}

Take two distinct co-prime numbers $p,q \in \N$. Consider the lens space $L(p,q)$, which is defined as the quotient of $S^3 = \{(z_1,z_2) \in \C^2\co |z_1|^2 + |z_2|^2 =1\}$ by an action of $\Z_p$. This action is defined by $t(z_1,z_2) \to (e^{2\pi t/p i}z_1, e^{2\pi tq/p i} z_2)$ for $t \in \Z_p$. Denote the equivalence class of an element $(z_1,z_2)$ by $[(z_1,z_2)]$. The action is free and isometric, wherefore $L(p,q)$ inherits a constant curvature metric. For an insightful discussion of these spaces see \cite{thurston1997three}. We ask for a bi-Lipschitz embedding $f\co L(p,q) \to \R^N$ for some $N$, and distortion, independent of $p$ and $q$. 

Collapsing theory states that the space can be covered by charts in which it resembles a normal bundle over some simpler space. For the space $L(p,q)$, such charts can be realized explicitly by the sets $U_1=\{[(z_1,z_2)] \in L(p,q) \co |z_2|<\sqrt{3}|z_1|\}$ and $U_2 = \{[(z_1,z_2)] \in L(p,q) \co |z_1| < \sqrt{3}|z_2|\}$. These two sets cover $L(p,q)$. Next, we will describe their geometry.

The sets $U_j$, for $j=1,2$, contain the ``polar'' circles $S_j$ defined by $|z_j|=1$. Further, $U_j$ can be expressed as normal bundles over $S_j$ with a curved metric. On the lift of $U_1$ to $S^3$ we can locally introduce co-ordinates $(\alpha,\theta,\phi)$ by $(\alpha,\theta,\phi) \to (\cos(\alpha)e^{i\theta}, \sin(\alpha)e^{i\phi})$ for $\alpha \in [0, \pi/3)$ and $(\theta, \phi) \in S^1 \times S^1$. The induced metric, becomes 

\begin{eqnarray*}
g &=& |\sin(\alpha)|^2|d\alpha|^2+ |\cos(\alpha)|^2 |d\theta|^2 + |\cos(\alpha)|^2|d\alpha|^2 + |\sin(\alpha)|^2 |d\phi|^2  \\
&=&  |\cos(\alpha)|^2|d\theta|^2 + |d\alpha|^2 + |\sin(\alpha)|^2 |d\phi|^2.
\end{eqnarray*}

We can also define another metric on the lift of $U_1$ by

$$g_f = |d\theta|^2 + |d\alpha|^2 + \alpha^2 |d\phi|^2.$$

Note that the norms induced by $g$ and $g_f$ differ by a factor at most $4$ (since $\sin(\alpha)^2 \geq \frac{1}{4}\alpha^2$ and $\cos(\alpha)^2 \geq \frac{1}{4}$ for $\alpha \in [0,\pi/3)$). The group $\Z_p$ acts on the lift of $U_1$ by isometries with respect to either of these metrics, and thus the latter descends to a metric $g_f'$ on $U_1$, which is up to a factor $4$ equal to the push-forward of the original metric $g$. In fact, $(U_1,g_f')$ is isometric to a subset of the holonomy bundle, which is given by the $\R^2$-vector bundle over $\frac{2\pi}{p}S^1$ (the circle of radius $2\pi/p$), with holonomy generated by a $\frac{2\pi q}{p}$-rotation. Since the metric tensors differ only up to a factor and both $(U_1,g)$ and $(U_1, g_f')$ are geodesically convex, we can see that the identity map $\iota_1\co (U_1, g) \to (U_1, g_f')$ is bi-Lipschitz with factor $4$. Similar analysis can be performed for $U_2$, where the bundle has holonomy generated by a rotation of angle $2\pi s/p$, where $sq\equiv 1 (\text{mod}(p))$. One is thus led to consider the problem of embedding flat vector bundles with non-trivial holonomy. These spaces will be discussed in the next subsection, where we indicate a construction for the desired bi-Lipschitz maps $F_j\co (U_j,g_f') \to \R^{n_j}$. 

Finally, given the embeddings $F_j$, we define the map $f=(F_1,F_2, d(\cdot, U_1))$, where $F_1$ and $F_2$ are extended using the McShane lemma. To observe that this map is bi-Lipschitz, we refer to similar arguments in the proof of Lemma \ref{doublingarg}, or the proof of \cite[Theorem 3.2]{langplaut}, which we present below. We remark, that there are many ways by which to patch up bi-Lipschitz embeddings from subsets into a bi-Lipschitz embedding of a whole. Our choice here is somewhat arbitrary. In fact, later in Lemma \ref{doublingarg} a slightly different choice is employed.

For convenience of the reader, we restate and prove Theorem 3.2 from \cite{langplaut}. This result is not necessary for the proofs below as we will use slightly different arguments. However, understanding this argument is helpful.

\begin{lemma}\label{langplaut} \cite[Theorem 3.2]{langplaut} Let $X$ be a metric space, $A, B \subset X$ be two sets such that $X \subset A \cup B$ and $f \co A \to \R^n$ and $g \co B \to \R^m$ be two maps. If $f$ and $g$ are $L$-bi-Lipschitz and are extended arbitrarily to be $L$-Lipschitz on $X$, then $F(x)=(f(x),g(x), d(x,A))$ defines a $CL$-bi-Lipschitz map $F\co X \to \R^{n+m+1}$ for some $C$, which depends only on $L$. 
\end{lemma}

\noindent \textbf{Proof:} Clearly $F$ is $3L$-Lipschitz since its three components are $L$-Lipschitz. Thus, we only need to prove the lower bound in Equation \ref{eq:bilip}. Let $x,y \in X$ be arbitrary. If we have $x,y \in A$ or $x,y \in B$, then 

$$|F(x)-F(y)| \geq |f(x)-f(y)| \geq \frac{1}{L}d(x,y),$$
or

$$|F(x)-F(y)| \geq |g(x)-g(y)| \geq \frac{1}{L}d(x,y),$$
respectively. Thus, assume that both $x$ and $y$ don't lie in the same subset $A$ or $B$. By symmetry, assume $x \in A$ and $y \not\in A$. There are two cases, $10L^2d(y,A) \leq d(x,y)$, or $d(x,y) < 10L^2d(y,A)$. In the latter case 

$$|F(x)-F(y)| \geq d(y,A) \geq \frac{1}{10L^2}d(x,y).$$

Consider next the first case. Define $y' \in A$ to be such that $d(y,y') \leq 2d(y,A)$. Then $|f(y)-f(y')| \leq Ld(y,y') \leq 2Ld(y,A) \leq \frac{1}{5L}d(x,y)$.

\begin{eqnarray*}
|F(x)-F(y)| &\geq& |f(x) - f(y)| \\
             &\geq& |f(x) - f(y')| - |f(y')-f(y)| \\
          &\geq& \frac{1}{L}d(x,y) - \frac{1}{5L}d(x,y) \geq \frac{1}{2L}d(x,y).
\end{eqnarray*}

Combining the two lower bounds, and the upper bound, we get that $F$ is $10L^2$-Bi-Lipschitz. 

\QED

\subsection{Flat vector bundles}

In the previous example, we reduced the problem of embedding a lens space to that of embedding flat vector bundles. We are thus motivated to consider the embedding problem for them. Consider the space $E_\theta^3$, which is a flat $\R^2$-bundle over $S^1$ with holonomy $\theta$. One has $E_\theta^3 = \R \times \C / \Z$, where the $\Z$-action is defined as $t(x,z)=(x+2\pi t,e^{it\theta} z) $ for $t\in \Z$ and $z \in \C$. Denote the orbit of $(x,z)$ in the quotient by $[x,z]$. Let $r$ be the distance function to the zero section, which is given by $r([x,z])=|z|$. To embed this space, we use a decomposition argument. Let $T_0=\{[x,z] \co |z| \leq 4\}$ and $T_j = \{[x,z] \co 2^{j-1} < |z|  < 2^{j+1}\}$, for $j \geq 1$ an integer. Assume first that each $T_j$ can be embedded with uniform bounds on distortion and dimension by mappings $f_j\co T_j \to R^n$. Then, we can collect the maps $f_j$ with disjoint domains by defining four functions $F_s$ for $s=0, \dots 3$ as follows. Consider the sets $D_s = \bigcup_{k=0}^\infty T_{s+4k}$, and define $F_s(x)=f_{s+4k}(x)$ for $x \in T_{s+4k}$. The resulting functions will be Lipschitz on their respective domains $D_s$. By an application of McShane extension Theorem \ref{thm:extension} we can extend them to $E_{\theta}^3$. The embedding we consider is $F=(r,F_0,F_1,F_2,F_3)$. 

To see that $F(x)$ is bi-Lipschitz, we need to establish the lower bound in Equation \ref{eq:bilip}. Take $x,y \in E_\theta^3$ arbitrary. First, if $d(x,y) \leq 4|r(x)-r(y)|$, then the lower bound is trivial. Thus, assume $d(x,y) \geq 4|r(x)-r(y)|$. We will identify a $j$ such that $x,y \in T_j$. This is obvious if $r(x)=r(y)$. We assume by symmetry that $r(x)>r(y)$ and $r(x)>0$. Then, the lower bound is obtained by using the lower-Lipschitz bounds for $f_j$ on $T_j$. We have the estimate $4 |r(x)-r(y)| \leq d(x,y) \leq r(x) + r(y) + 1$. Thus, $r(y) \geq r(x) - |r(x)-r(y)| \geq 3/4 r(x) - 1/4r(y) - 1/4$, and $r(y) \geq 3/5 r(x) - 1/5$. In particular, if $x \in T_0$, then $r(y) \leq r(x) \leq 4$, and thus $y \in T_0$. But, if $2^{j}\leq r(x) \leq 2^{j+1}$ for some $j \geq 2$, then $r(y) \geq 1/2 r(x) \geq 2^{j-1}$, and thus $x,y \in T_j$.

Next, we discuss the construction of the embeddings $f_j$. On $T_0$ the injectivity radius is bounded from below. Thus, we can find $f_0$ by patching up a finite number of bi-Lipschitz charts. More specifically, we use doubling and the Lemma \ref{doublingarg}. For $T_j$ when $j \geq 1$, we will first modify the metric. Use co-ordinates $(x, r,\theta) \to (x,re^{i\theta}) \in \R^3$. The metric on $T_j$ (or its lift in $\R^3$) can be expressed as

$$g=|dx|^2 + |dr|^2 + r^2|d\theta|^2.$$

We change the metric to $g_f=|dx|^2 + |dr|^2 + 2^{2j}|d\theta|^2=g_1$. As metric spaces, the space $(T_j, g)$ is 10-bi-Lipschitz to $(T_j, g_f)$. The new metric space $(T_j,g_f)$ is isometric to $(2^{j-1},2^{j+1}) \times T^j_\theta$, where $T^j_\theta$ is a torus defined by $\R^2 / \Z^2$, where the $\Z^2$ action is given by $(n,m)(x,y) = (x+2\pi n, y + n\theta + 2\pi m)$. The product manifold $(2^{j-1},2^{j+1}) \times T^j_\theta$ can be embedded by embedding each factor separately. Thus, the problem is reduced to embedding a flat torus. This can be done by choosing a short basis and is explained in detail in \cite{naorfourier}.

We comment briefly on the problem of embedding general flat vector bundles. Consider a flat $\R^n$-bundle over $S^1$. The argument remains unchanged for $T_0$, but for $T_j$ when $j \geq 1$ one can no longer modify the metric to be a flat product metric. Instead, one needs a further decomposition argument of the $\R^d$ factor similar to that used for lens spaces. If the base is more complicated than $S^1$, we are led to an induction argument where these decompositions are used to reduce the embedding problem for spaces that are simpler in some well-defined sense. 

Finally, we highlight the main ideas used in the previous embedding constructions for lens spaces and holonomy-bundles. First, the space is decomposed into sets (``charts'') on which we have in some sense a simpler geometry. The embeddings for these charts can be patched together to give an  embedding of the space. In the case for lens spaces, the simpler geometry was that of a vector bundle. For vector bundles, the simpler geometry was that of a product manifold. The simpler geometries may need to be decomposed several times, but we can bound the number of iterated decompositions required by the dimension of the space. Ultimately, the problem is reduced to embedding something very simple such as a non-collapsed space with a lower bound on the injectivity radius. The main technical issues arise from the fact that high-dimensional spaces may require several decompositions and that for abstract spaces it is complicated to identify good charts for which a simpler structure exists.

\subsection{Quotients of the Heisenberg Group and Nilpotent groups}
\label{examp:nilpotentlie}
The Heisenberg group may be described as the simply connected Lie group of upper triangular $3\times 3$-matrices with diagonal entries equal to one.

\begin{equation}
\mathbb{H} = \left \{ \left. \left (\begin{array}{ccc} 1 & a &c \\ 0 & 1 & b \\ 0 & 0 & 1 \end{array} \right) \right| a,b,c \in \R \right \}
\end{equation}

Define a lattice 

$$\Gamma  = \left \{ \left. \left(\begin{array}{ccc} 1 & a & c \\ 0 & 1 & b \\ 0 & 0 & 1 \end{array} \right) \right| a,b,c \in \Z.\right \} \subset \mathbb{H}$$
which acts on $\mathbb{H}$  by left multiplication. We can define the quotient space $\mathbb{H} / \Gamma$. Next metrize this base as follows. Take a basis for the Lie algebra of left-invariant vector fields given by 

$$X = \left(\begin{array}{ccc} 0 & 1 & 0 \\ 0 & 0 & 0 \\ 0 & 0 & 0 \end{array} \right) , Y = \left(\begin{array}{ccc} 0 & 0 & 0 \\ 0 & 0 & 1 \\ 0 & 0 & 0 \end{array} \right), Z = \left(\begin{array}{ccc} 0 & 0 & 1 \\ 0 & 0 & 0 \\ 0 & 0 & 0 \end{array} \right).$$

Consider the space $\mathbb{H}_\epsilon = (\mathbb{H}/ \Gamma,g_\epsilon)$ with the metric defined by setting $X \perp Y \perp Z$, $|X|=|Y|=1$ and $|Z|=\epsilon$. Initially, this metric is defined on $\mathbb{H}$, but since $\Gamma$ acts by isometries the metric descends onto the quotient. A direct computation shows that the sectional curvatures of this space lie in $$\left[-\frac{3\epsilon^2}{4}, \frac{\epsilon^2}{4}\right].$$ As $\epsilon \to 0$, the spaces $\mathbb{H}_\epsilon$ Gromov-Hausdorff converge to a torus $T^2$. The Gromov-Hausdorff approximation is given by the $S^1$-fibration map $\pi \co \mathbb{H}_\epsilon \to T^2$ as follows,

$$\pi\co \left[ \left(\begin{array}{ccc} 1 & a & c \\ 0 & 1 & b \\ 0 & 0 & 1 \end{array} \right) \right] \to ([a],[b]) \in \frac{1}{2\pi}S^1 \times \frac{1}{2\pi}S^1 = T^2,$$

\noindent where $[a]$ is the fractional part of $a \in \R$. By $\frac{1}{2\pi} S^1$ we mean the circle of unit length. This map is easily seen to be well-defined on $\mathbb{H}_\epsilon$. We wish to find a bi-Lipschitz embedding for $\mathbb{H}_\epsilon$ with distortion independent of $\epsilon>0$, as long as its sufficiently small.

An embedding for $\mathbb{H}_\epsilon$ can be visualized using Lemma \ref{doublingarg}. Consider a small ball $B_p(\delta) \subset T^2$, with some $\delta>0$ independent of $\epsilon$. Since $T^2$ is non-collapsed, there exists a $\delta>0$ such that $B_p(\delta)$ is isometric to a ball in the plane. Then, since this ball is contractible, we can conclude that $\pi^{-1}(B_p(\delta))$ is diffeomorphic to $B_p(\delta) \times S^1$. Because the curvature is almost flat, we can show that this diffeomorphism can be chosen to be a bi-Lipschitz map with small distortion. This claim can be proven by a similar argument to Lemma \ref{metricnilp} in the appendix. On the other hand, the space $B_p(\delta) \times S^1$ is easy to embed since it splits as a product of two simple spaces. The result can then be deduced by covering $T^2$, and thus $\mathbb{H_\epsilon}$, by a controlled number of these sets, and patching up the embeddings using Lemma \ref{doublingarg} or a repeated application of Lemma \ref{langplaut}.

\section{Frequently used results}

\label{prelimin}

Below, there are a number of constants determined from specific theorems. Since we wish to prove embedding theorems with quantitative bounds on the distortions, we need to be careful about the dependence on these parameters. In order to keep track of these different constants and functions throughout the chapter, we fix their meaning. The special constants are: $C'(n)$ from Lemma \ref{almostflat} and Definition \ref{quasiflatdef}, $c(n)$ is the maximum of the constants necessary for Lemmas \ref{jordan}, \ref{almostflat} and \ref{flatmanif}, $\epsilon(n)$ is fixed via Lemmas \ref{strat}, \ref{approx} and \ref{approxvect}, $\delta(n)$ from Lemma \ref{metricnilp}, and $C(n)$ from Lemma \ref{nilpdiam}. Also, we will choose $c(n)$ increasing and $\epsilon(n)$ decreasing in $n$. Their values will be considered fixed throughout the paper, although the precise value is left implicit. In general, we will use the convention that a constant $M(a_1, a_2, \dots, a_n)$ is a quantity depending on $a_1, \dots, a_n$ only.

\subsection{Embedding lemmas}

Consider any metric space $X$. Throughout this paper we will denote open balls by $B_x(r) =B(x,r)= \{y \in X | d(y,x) < r\}$. For a metric space $X$ we call $N \subset X$ a \emph{$\epsilon$-net} for $X$ if for every $x \in X$ there is a $n \in N$ such that $d(x,n)\leq \epsilon$, and if $x,y \in N$ then $d(x,y) > \epsilon$. Such nets can always be constructed using Zorn's lemma. For a subset $A \subset X$ and $r>0$, we denote by $$N_r(A)=\{x \in X | d(x,A) < r\}$$ the $r$-tubular neighborhood of a subset $A$. Recall, $d(x,A) = \inf_{y \in A} d(x,y)$. The metric spaces in this paper will be complete connected manifolds $M$ and orbifolds $O$ equipped with a Riemannian distance function $d$. 

Our Lipschitz mappings are initially defined locally and for a global embedding we need some type of extension theorem.

\begin{theorem}\label{thm:extension}(McShane-Whitney, \cite{heinonen2005lectures}) Let $X$ be a metric space and $A \subset X$. Then any L-Lipschitz function $f\co A \to \R^n$ has a $\sqrt{n}L$-Lipschitz extension $\tilde{f}\co X \to \R^n$ s.t $\tilde{f}|_A = f$.
\end{theorem}

\textbf{Remark:} McShane is a good extension result to use because of its generality. However, for bounded curvature manifolds better constants could be attained by the use of a generalized Kirzbraun's theorem \cite{kirsz}. For most of the paper Kirzbraun could be used instead and slightly better constants would ensue. 

In addition to extension, we will use decomposition arguments in two ways. On the one hand we have spaces that admit certain splittings, such as cones and products.

\begin{lemma}\label{product} Let $(X,d_X)$, $(Y,d_Y)$ be metric spaces and take $X \times Y$ with the product metric $d((x_1,y_1), (x_2,y_2))=\sqrt{d_X(x_1,x_2)^2 + d_Y(y_1,y_2)^2}$. If $X$ and $Y$ admit bi-Lipschitz embeddings $f\co X \to \R^N$ and $g\co Y \to \R^M$ with distortion $L$, then $X \times Y$ admits a bi-Lipschitz embedding $(f,g)\co X \times Y \to \R^{N + M}$ with distortion at most $\sqrt{2}L$.
\end{lemma}

The proof of this result is trivial. The product space could also be equipped with different bi-Lipschitz equivalent metrics. In other cases, the space admits a conical splitting. We define a cone over a metric space.

\begin{definition} Let $X$ be a compact metric space with $\diam(X) \leq \pi$. Then we define the cone over $X$ as the space $C(X) = \{(t,x) \in \R \times X\} / \sim$, where $(0,x) \sim (0,y)$ for any $x, y\in X$, equipped with  the metric

$$d([t,x], [s,y]) = \sqrt{t^2+s^2-2ts\cos(d(x,y))}.$$ 
\end{definition}

The metric cones which arise in our work will always be quotients of $\R^n$ by a group with a fixed point. Thus, we state the following lemma whose proof is immediate. Note that all the groups of the following lemma are forced to be finite.

\begin{lemma}\label{euclideanquot} Let $\Gamma$ be a group acting properly discontinuously and by isometries on $\R^n$ which fixes $0 \in \R^n$. Then the group $\Gamma$ acts properly discontinuously and by isometries on $S^{n-1} \subset \R^n$, and $\R^n / \Gamma$ is isometric to $C(S^{n-1} / \Gamma)$.
\end{lemma}

\begin{lemma}\label{cone} Let $X$ be a compact metric space with $\diam(X)\leq \pi$, and $Y = C(X)$. If $X$ admits a bi-Lipschitz embedding to $\R^n$ with distortion $L$, then $Y$ admits a bi-Lipschitz embedding to $\R^{n+1}$ with distortion $20L$.
\end{lemma} 

\textbf{Proof:} Let $r$ be the radial function on $Y$, where points on $Y$ are given by equivalence classes $[(r,x)]$ for $r \geq 0, x \in X$. We will simplify notation by suppressing the brackets. Assume that $f\co X \to \R^n$ is bi-Lipschitz with distortion $L$ and $0 \in \text{Im}(f)$. Define $g(r,x) = (Lr, rf(x))$ on $Y$ and show that it is the desired embedding. If $(s,x)$ and $(t,y)$ are points in $Y$ and $s\leq t$, then 

$$|g(s,x)-g(t,y)| \leq L|s-t| + s|f(x)-f(y)| + |s-t||f(y)|.$$

Note that $|f(y)| \leq L\pi$ (because $\diam(X) \leq \pi$). Further, $$d((s,x),(t,y))=\sqrt{t^2+s^2-2ts\cos(d(x,y))}\geq 2sd(x,y)/\pi.$$ Using $s|f(x)-f(y)| \leq Lsd(x,y)$ we get 

$$|g(s,x)-g(t,y)| \leq (L\pi+L+2L/\pi) d((s,x),(t,y)) \leq 10L d((s,x),(t,y)).$$

Next, we derive the necessary lower bound from Equation \ref{eq:bilip}. Take any two points $(s,x)$ and $(t,y)$ in $Y$ with $s\leq t$. Either $|s-t| \geq \frac{1}{20L^2}d((s,x),(t,y))$ or not. In the first case

$$|g(s,x)-g(t,y)| \geq L|s-t| \geq  \frac{1}{20L}d((s,x),(t,y)).$$

In the latter case, $|s-t| \leq \frac{1}{20L^2}d((s,x),(t,y))$, and simple estimates give

\begin{eqnarray*}
|g(s,x)-g(t,y)| &\geq& |g(t,x)-g(t,y)| - |s-t||f(x)| \\
                      &\geq& \frac{1}{2L} d((t,x),(t,y)) - \frac{\pi}{20L}d((s,x),(t,y)) \geq  \frac{1}{4L}d((s,x),(t,y)).
\end{eqnarray*}

These estimates complete the proof.
\QED

\noindent \textbf{Remark:} Similar constructions would also apply for a spherical suspension of a metric space  (see \cite{BGP}) but we do not need that result. 

When a simple splitting structure doesn't exist at a given scale, it may exist at a smaller scale. To enable us to accommodate for this we need the following ``doubling argument''. To state it, we need to define doubling metric spaces. We emphasize that a metric space must be doubling in order to admit a bi-Lipschitz embedding into $\R^N$. 

\begin{definition} \label{doublingdef} A doubling metric space with doubling constant $D$ is a metric space $X$ such that for any ball $B_x(r) \subset X$ there exist points $p_1, \dots, p_{D}$ such that

$$B_x(r) \subset \bigcup_{i=1}^D B_{p_i}\left( \frac{r}{2} \right).$$
\end{definition}

The following Lemma could be proven by repeatedly applying Lemma \ref{langplaut}, but we give a different and slightly more direct argument.

\begin{lemma} (Doubling Lemma) \label{doublingarg} Assume $0<r < R, 0<N,l$ are given and that $X$ is a doubling metric space with doubling constant $D$. Fix a point $p \in X$. If for every point $q \in B_p(R)$ there is a $l$-bi-Lipschitz embedding $f_q\co B_q(r) \to \R^{N}$, then there is a $L$-bi-Lipschitz embedding $G \co B_p(R) \to \R^M$ for some $M>0$. Further, we can bound the distortion of the bi-Lipschitz embedding $G$ by $L\leq L(N,l,D, R/r)$ and the target dimension by $M \leq M(N,l,D, R/r)$.
\end{lemma}

\noindent \textbf{Proof:} Construct a $r/8$-net $Q=\{q_i\}$ in $B_p(R)$ and let $f_i\co B_{q_i}(r) \to \R^N$ be the bi-Lipschitz mappings assumed. The size of the net is $|Q| = H \leq D^{\log_2(R/r)+4}$. Translate the mappings so that $f_i(q_i)=0$, and scale so that the Lipschitz-constant is $1$. This will ensure that for any $a,b \in B_{q_i}(r)$ we have $|f_i(a)-f_i(b)|\geq 1/l^2 d(a,b)$. Further, let $F\co x\to (d(x,q_i))_i \in \R^H$ be a distance embedding from the net. Group the $q_i$ into $K$ groups $J_k$ such that if $q_i, q_j \in J_k$ then $d(q_i,q_j) > 4r$. The number of points in the net $Q$ within distance at most $4r$ from any $q_i$ is at most $D^4$. Thus, a standard graph coloring argument such as in \cite{Heinonen2000} will furnish the partition $J_k$ with $K$ sets with $K \leq D^4$. 

For each $J_k$, construct a map $g_k: \bigcup_{q_i \in J_k} B_{q_i}(r) \to \R^N$ by setting it equal to $f_i$ on $B_{q_i}(r)$. Because the distance between the balls is at least $2r$, the Lipschitz constant will not increase. To see this take arbitrary $a \in B_{q_i}(r)$ and $b \in B_{q_j}(r)$, where $i \neq j$ and $q_i,q_j \in J_k$. Then $d(a,b) \geq 2r$ and $|f(a)-f(b)|\leq |f(a) - f(q_i)| + |f(b) - f(q_j)| \leq 2r \leq d(a,b)$. Extend $g_k$ using Lemma \ref{thm:extension} to give a $\sqrt{N}$-Lipschitz map from the entire ball $B_p(R)$, and denote it by the same name. For the embedding combine these all into one vector

$$G = (g_1, \dots, g_K, F).$$

Clearly this is $2\sqrt{N}D^4 + 2D^{\log_2(R/r)+4}$-Lipschitz. Also, it is a map $G \co B(p,R) \to \R^M$ with $M=KN + H$. For the lower Lipschitz bound in \eqref{eq:bilip} take arbitrary $a,b \in B_p(R)$. First assume $d(a,b) < r/2$. Then, there is a $q_i$ such that $d(a,q_i) < r/8$ and $d(b,q_i) < r$. Both belong to the same ball $B_{q_i}(r)$ and for the $k$ such that $q_i \in J_k$ we have $|G(a)-G(b)| \geq |g_k(a)-g_k(b)| = |f_i(a)-f_i(b)| \geq 1/l^2 d(a,b)$. Next assume $d(a,b) \geq r/2$. Then, there is a $q_i \in Q$ such that $d(a,q_i)< r/8 \leq d(a,b)/4$ and thus $d(b,q_i) \geq d(a,b) - d(a,b)/4 \geq 3/4 d(a,b)$. Thus, $|d(a,q_i)-d(b,q_i)| > d(a,b)/2$ and we get $|G(a)-G(b)| \geq |F(a)-F(b)| \geq d(a,b)/2$.
\QED

\begin{definition} \label{def:quasi} Let $X,Y$ be metric spaces and $s>0$. A mapping $f \co X \to Y$ is called a ($s$-)quasi-isometry if for all $x,y \in X$

$$d_X(x,y) - s \leq d_Y(f(x),f(y)) \leq d_X(x,y) + s,$$
and $Y \subset N_s(f(X))$. 
\end{definition}

\paragraph{Remark:} It is easy to see that a $s$-quasi-isometry will be a $2$-bi-Lipschitz map when restricted onto a $2s$-net.

Using quasi-isometries, one can define a distance between compact metric spaces and an associated notion of convergence. For more detailed discussion on convergence of metric and metric measure spaces, we refer the reader to \cite{Sormani2012}. See also the excellent references \cite{BuragoGeom,Gbook}.

\begin{definition} \label{def:grhaus} Let $X,Y$ be compact metric spaces. We denote by $d_{GH}(X,Y)$ the infimum of numbers $s>0$ such that there are maps $f \co X \to Y$ and $g \co Y \to X$ which are $s$-quasi-isometries. Further, we say that a sequence of compact metric spaces $X_i$ Gromov-Hausdorff converges to another metric space $X$ if $$\lim_{i \to \infty} d_{GH}(X_i, X)=0.$$ 
\end{definition}

\begin{lemma} (Gromov-Hausdorff-lemma) \label{gromhaus} Let $X$ be a doubling metric space with doubling constant $D$ and let $\epsilon, \epsilon',l',l,m,N>0$ be constants. Assume that $Y$ is a metric space admitting a bi-Lipschitz map $h\co Y \to \R^m$ with distortion $l$ and that $d_{GH}(X,Y)\leq \epsilon'$. If for every  $p \in X$ there is a $l'$-bi-Lipschitz embedding $f_p \co B_p(\epsilon) \to \R^N$, then there is a $L$-bi-Lipschitz map $G\co X \to \R^K$, with distortion $L \leq L(D, l,l',m, \epsilon'/ \epsilon, N)$ and target dimension $K \leq K(D,l,l',m, \epsilon'/\epsilon, N)$.
\end{lemma}

\textbf{Proof:} Fix $M=\epsilon'/\epsilon$. Let $g\co X \to Y$ be a Gromov-Hausdorff approximation. By Lemma \ref{doublingarg} we can first construct for every $p \in X$ an embedding $f'_p\co B_p(100 l \sqrt{m}\epsilon') \to \R^{N'}$. By a scaling and dilation, we can assume that $f'_p$ are 1-Lipschitz and $f'_p(0)=0$. The distortion and the dimension $N'$ will depend on $l,l',D,N,m$ and $M$. Grouping $f_p$ similarly to the proof in Lemma \ref{doublingarg}, we can define corresponding maps $g_k \co X \to \R^N$ for $k=1\dots H'$ for some $H'$ depending on the parameters of the problem. Then, we can construct a map $F\co X \to \R^{K'}$ with $F=(g_1, \dots, g_H)$ with $K'=H'N$. Moreover, this construction guarantees that for every $x,y \in X$ with $d(x,y)\leq 10l\epsilon' \sqrt{m}$ we have a $p \in X$ and an index $i$ such that $x,y \in B_p(100l\sqrt{m}\epsilon')$ for some $p$ and $g_i|_{B_p(100l\sqrt{m}\epsilon')}=f_p$. This can be used to give the lower bi-Lipschitz bound for Equation \ref{eq:bilip} for any pair $x,y \in X$ such that $d(x,y) \leq 10l\epsilon'\sqrt{m}$.

Next, take a $\delta$-net $N_\delta$ for $\delta= 4\epsilon'$. The map $h \circ g|_{N_\delta}\co N_{\delta} \to \R^m$ is a $2l$-bi-Lipschitz map because $g$ is $2$-bi-Lipschitz on a $4\epsilon'$-net. Extend this to a map $H\co X \to \R^m$ which is $2l$-bi-Lipschitz on $N_\delta$ and $2l\sqrt{m}$-Lipschitz on $X$. The map $H$ is used to give a lower Lipschitz bound for pairs $x,y \in X$ with $d(x,y) > 10l\sqrt{m}\epsilon'$. The proof can now be completed by similar estimates as in Lemma \ref{doublingarg} and defining $G(x)=(H(x),F(x))$. This gives a map $G \co X \to \R^{K'+m}$. 

\QED

The main source of Gromov-Hausdorff approximants will be via quotients so we state the following Lemma. 

\begin{lemma} \label{closeness} Assume $(X,d_X),(Y,d_Y)$ are metric spaces and $\delta>0$. If $\pi \co X \to Y$ is a quotient map and if we define $d_Y(a,b) = d_X(\pi^{-1}(a), \pi^{-1}(b))$ for all $a,b \in Y$ and $\diam(\pi^{-1}(a)) \leq \delta$ for all $a \in X$, then $d_{GH}(X,Y) \leq 2\delta$.
\end{lemma}

\subsection{Group actions and quotients} \label{subsec:group}

Assume that $X$ is a proper metric space, i.e. that $B_p(r)$ is precompact for every $p \in X$ and $r>0$. Further, assume that $\Gamma$ is a discrete Lie group acting on $X$ by isometries. For any group, its identity element is denoted by either $e$ or $0$, depending on if it is abelian. For any $p \in X$, its \textit{isotropy} group is denoted by $\Gamma_p =\{\gamma \in \Gamma | \gamma p = p\}$. We say that the action is \textit{properly discontinuous} if for every $p \in X$ there exists a $r>0$, such that $\{\gamma \in \Gamma | \gamma(B_p(r)) \cap B_p(r) \neq \emptyset \} = \Gamma_p$ and $\Gamma_p$ is a finite group. We will study quotient spaces, which will be denoted by $X/ \Gamma$, where the action of $\Gamma$ is properly discontinuous. These spaces are also in some cases referred to as orbit spaces.

For $a \in X$ we will denote its orbit or equivalence class in $X / \Gamma$ by $[a]$. Since the action is by isometries, we can define a quotient metric by 
$$d([a], [b])=\inf_{\gamma, \gamma' \in \Gamma} d(\gamma a, \gamma' b).$$ 
This distance makes $X / \Gamma$ a metric space. For more terminology see \cite[Chapter 5]{ratcliffehyper}. The global group action may be complicated, but the local action can often be greatly simplified. The following lemma is used to make this precise. First recall, that for any set of elements $S \subset G $ the smallest subgroup containing them is denoted by $\langle S \rangle$. We say that $S$ generates $\langle S \rangle$. For a group $G$, we say that $\gamma_1, \dots, \gamma_n$ are generators of the group, if $G = \langle \gamma_1, \dots, \gamma_n \rangle$. 

\begin{lemma}\label{small} Let $X$ be a metric space and $\Gamma$ a discrete Lie group of isometries acting properly discontinuously on $X$. Take a point $p \in X$. Let $\Gamma_p(r)=\langle g  \in \Gamma | d(gp,p) \leq 8r \rangle$ be the subgroup generated by elements such that $d(p,gp) \leq 8r$. Then $B_{[p]}(r) \subset Y = X / \Gamma$ is isometric to $B_{[p]'}(r) \subset X / \Gamma_p(r)$ with the quotient metric.
\end{lemma}

\textbf{Proof:} We will denote the cosets, functions and points in  $Y' = X / \Gamma_p(r)$ using primes and corresponding objects in $Y$ without primes. For example, elements of $X / \Gamma_p(r)$ can be represented by a $\Gamma_p(r)$-cosets $[x]'$, for $x \in X$, and elements of $Y = X / \Gamma$ will be represented as $[x]$. The action is still properly discontinuous and the distance function on the orbit space is given by $d'([x]',[y]')=\inf_{\gamma, \gamma' \in \Gamma_p(r)} d(\gamma x, \gamma' y)$. Since $\Gamma_p(r) \subset \Gamma$, we have a continuous map $F\co (Y',d') \to (Y,d)$ which sends a  $\Gamma_p(r)$-coset to the $\Gamma$-coset that contains it, i.e. $[x]' \to [x]$. We will next show that $F$ maps $B_{[p]'}(r) = \{y' \in Y' | d'([p]', y') < r\}$  isometrically onto $B_{[p]}(r)$,  from which the conclusion follows.

Let $[a]', [b]' \in B_{[p]'}(r) \subset X/\Gamma_p(r)$ and as before $F([a]')=[a], F([b]')=[b]$. Clearly $d([a]',[b]') \geq d([a], [b])$, as the $\Gamma$-cosets are super-sets of $\Gamma_p(r)$-cosets. We can choose the representatives of the cosets $a,b$ such that $a,b \in B_p(r) \subset X$. By triangle inequality $d([a],[b])<2r$, and for any small $\epsilon >0$ there is a $\gamma_\epsilon \in \Gamma$ such that $ d(a,\gamma_\epsilon b) < d([a],[b]) + \epsilon < 2r $. But then $\gamma_\epsilon b \in B_p(3r) \subset X$ and $d(\gamma_\epsilon b, b) \leq 4r.$ Therefore, $d(\gamma_\epsilon p, p) \leq d(\gamma_\epsilon p, \gamma_\epsilon b) + d(\gamma_\epsilon b, b) + d(b,p) \leq r + r + 4r \leq 6r$. Thus, $\gamma_\epsilon \in \Gamma_p(r)$ and therefore $d([a]',[b]') \leq d(a,\gamma_\epsilon b) \leq d([a], [b]) + \epsilon$. The result follows since $\epsilon$ is arbitrary. 

Finally, the map is onto, because if $d([p],[a]) < r$, then we can choose the representatives such that $d(p,a)<r$. In particular $[a]' \in B_{[p]'}(r)$, and $F([a]')=[a]$ by definition.

\QED

Motivated by the previous Lemma, we define a notion of a local group. 

\begin{definition} Let $\Gamma$ be a discrete Lie group acting properly discontinuously on $X$. We denote by $\Gamma_p(r)=\langle g  \in \Gamma | d(gp,p) \leq 8r\rangle$, and call it the local group at scale $r$ (and location $p$).
\end{definition}

In collapsing theory it is often useful to find invariant points and submanifolds. However, often it is easier to first construct almost invariant points/submanifolds, and then to apply averaging to such a construction to give an invariant one. For this reason, we will recall Grove's and Karcher's center of mass technique. See \cite{grove1973conjugatec,center} for a more detailed discussion and proofs. 

\begin{lemma}(Center of mass Lemma)\label{centerofmass} Let $M$ be a complete Riemannian manifold with a two-sided sectional curvature bound $|K| \leq \kappa$ and $p \in M$ a point. Further, let $\mu$ be a probability measure supported on a ball $B_p(r)$ of radius
$r \leq \min\left\{\frac{\pi}{2\sqrt{\kappa}}, \mathrm{inj}(M,p)/2\right\}$. There exists a unique point $C_\mu \in B_p(r)$ which minimizes the functional $$F(q) = \int d(q,x)^2 d\mu_x.$$ 

The minimizer $C_\mu$ is called the center of mass. The center of mass is invariant under isometric transformations of the measure. If $f\co M \to M$ is an isometry, and $f_*(\mu)$ is the push-forward measure of $\mu$, then

$$f(C_{\mu})=C_{f_*(\mu)}.$$
\end{lemma}

As a corollary, one obtains the existence of fixed points for groups with small orbits. If $G$ is any group acting on a manifold $M$, a point $p \in M$ is called a fixed point of a group action if $gp=g$ for all $g \in G$. A more traditional form of this lemma appears in \cite[Lemma 5.9]{ratcliffehyper}.

\begin{corollary}(Fixed point Lemma)\label{fixedpoint}
 Let $M$ satisfy the same assumptions as in Lemma \ref{centerofmass}. Let a compact Lie group $G$ act by isometries on $M$. If $$\mathrm{diam}(G_p) < \min\left\{\frac{\pi}{2\sqrt{\kappa}}, \mathrm{inj}(M,p)/2\right\},$$ then there exists a fixed point $q$ of the $G$-action such that $d(q,p) \leq \mathrm{diam}(G_p)$.
\end{corollary}

In fact, one can use the Haar measure $\mu$ on $G$, and consider the push-forward measure $f_*(\mu)$ by the map $f\co G \to M$ which is given by $f(g)=gp$. The center of mass $C_{f_*(\mu)}$ will be invariant under the group action by $G$.

We will repeatedly use an argument that allows us to derive metric conclusions from a stratification of a finite index subgroup. The idea of the statement is that a semidirect product with the stratified normal subgroup will itself  admit a stratification. Such stratifications arise from Lemma \ref{strat} below. Assume $\Gamma$ is a group. Then we call $|\cdot| \co \Gamma \to \R$ a subadditive norm if the following hold.

\begin{itemize}
\item For all $g \in \Gamma$, $|g| \geq 0$ and $|e|=0$.
\item $|ab| \leq |a| + |b|$.
\item $|a| = |a^{-1}|$
\end{itemize}

For groups acting on spaces, one can choose a natural class of subadditive norms.

\begin{definition} Let $\Gamma$ act on a metric space $X$, and let $p \in X$ be fixed. We define the subadditive norm at $p$ for $g \in \Gamma$ as

$$|g|_p = d(gp, p).$$
\end{definition}

If $g \in \text{Isom}(X)$ is an isometry, we denote $|g|_p=d(gp,p)$. Further, if $\Gamma$ acts on a Lie group $N$, then we will often choose $p=e$ and denote $|\gamma|=|\gamma|_e$.

\begin{lemma}(Local Group Argument)\label{loc} Fix an arbitrary scale parameter $l > 0$. Let $\Gamma$ be a group with a subadditive norm $|\cdot|$, and assume it admits a short exact sequence $$0 \to \Lambda \to \Gamma \to H \to 0,$$
where $H$ is a finite group of size $|H| < \infty$. Assume $\Lambda_0 <\Lambda$ is a normal subgroup  in $\Gamma$ generated by all $g \in \Lambda$ of length $|g| \leq l$, and further assume that every element $g \in \Lambda \setminus \Lambda_0$ has length $|g| \geq 10|H| l$. Then, the subgroup $\Gamma_0<\Gamma$ generated by all $g \in \Gamma$ with $|g| \leq l$ admits a short exact sequence

$$0 \to \Lambda_0 \to \Gamma_0 \to H_0 \to 0,$$
where $H_0 < H$. Further, any element $g \in \Gamma_0$ is can be expressed as
$$g = \prod_{i=1}^{|H|+1} g_i \lambda ,$$
where $|g_i| \leq l$ and $\lambda \in \Lambda_0$.
\end{lemma} 

\textbf{Proof:} Define the group by generation: $\Gamma_0=\langle g \in \Gamma | |g| \leq l \rangle$. Define $S'$ to be the collection of left-cosets of $\Lambda$ in $\Gamma$ that contain a representative $s \in \Gamma$ with $|s| \leq l$. Define $S$ to consist of a representative $s$ for each left coset in $S'$ with $|s| \leq l$. Clearly $|S|=|S'| \leq |H|$. Also, the cosets in $S'$ coincide with $\Lambda_0$-cosets, because if $s=\lambda s'$ for some $\lambda \in \Lambda$ and $|s|, |s'| \leq l$, then $|\lambda| \leq |s|+|s'| \leq 2l$, and thus $\lambda \in \Lambda_0$. 

Construct $W$ using all the ``words'' $w$ of length at most $|H| + 1$ formed by multiplying elements $s$ and $s^{-1}$ for $s \in S$, including the empty word which represents the identity
element of $\Gamma$. None of these words represent elements of $\Lambda \setminus \Lambda_0$ since for any $w \in W$ we have $|w|\leq (|H|+1)l$. Denote by $W'=\{[w]_{\Lambda_0} | w \in W\}$ the left-cosets of $\Lambda_0$ represented by elements in $W$. Also, let $\overline{\Gamma}=\bigcup_{A \in W'} A \subset \Gamma$. We will show that $\overline{\Gamma}=\Gamma_0$. In other words, we show that $W$ exhausts the left-cosets of $\Lambda_0$ in $\Gamma_0$, and that these cosets have a product given by the product of their representatives. To do so, it is sufficient to show that any product of two elements in $W$ is equivalent modulo $\Lambda_0$ to an element in $W$. This shows that $\overline{\Gamma}$ is closed under products. Since also $\overline{\Gamma}<\Gamma_0$ as it is generated by $S \cup \Lambda_0$, we have $\Gamma_0=\overline{\Gamma}$.

Take any word $w'$ with length longer than $|H|+1$. Consider the subwords $w_1', \dots, w_{|H|}'$ formed by the first $1, \dots, |H|+1$ letters of $w'$, respectively. Because there are only $|H|$ cosets, two of them will be in the same $\Lambda$-coset in $\Gamma$. Thus, assume $w_i$ and $w_j$ with $i<j$  belong to the same coset. Let $w'=w_j v_j$, where $v_j$ is the product of the remaining letters. We have $w_iw_j^{-1} \in \Lambda$. But also $|w_iw_j^{-1}|\leq (2|H|+1)l$ so $w_iw_j^{-1} \in \Lambda_0$. Thus, we have for the $\Lambda_0$-cosets $[w_i]_{\Lambda_0}=[w_j]_{\Lambda_0}$, and thus $[w']_{\Lambda_0}=[w_j]_{\Lambda_0}[v_j]_{\Lambda_0}=[w_i]_{\Lambda_0}[v_j]_{\Lambda_0}=[w_i v_j]_{\Lambda_0}$. The word $w_iv_j$ is shorter than $w'$. This process can be continued as long as $w'$ is longer than $|H|+1$. Thus, any word longer than $|H|+1$ is equivalent modulo elements in $\Lambda_0$ to an element in $W$. In particular, the cosets formed by products of cosets corresponding to elements in $W$ can be represented by elements contained in $W$. This allows one to define a group structure for the cosets represented by $W$, and therefore a subgroup
$H_0 < H$. It is immediate, that $H_0$ contains all the cosets generated by elements in $S$.

 By the previous argument, $\Gamma_0$ is a finite extension of $\Lambda_0$ of index at most $|H_0|$. Also, $\Lambda_0$ is normal in $\Gamma_0$ since it is normal in $\Gamma$. The final statement about representing elements of $g$ as products follows from the above, since $g$ belongs to a coset represented by an element in $W$.

\QED

\subsection{Nilpotent Lie groups, lattices and collapsing theory}

Nilpotent geometries occur as model spaces for collapsing phenomena. Thus, we start by an estimate relating to the curvature of a general simply connected nilpotent Lie group $N$ equipped with a left-invariant metric $g$. For definitions of nilpotent Lie groups as well as the computations involved here, see \cite{corwin2004representations}. Associate to the Lie group $N$ its nilpotent Lie algebra $\mc{n}$. We construct a triangular basis for $\mc{n}$ as follows. Let $\mc{n}_0 = \mc{n}, \mc{n}_1 = [\mc{n},\mc{n}], \dots, \mc{n}_{k+1} = [\mc{n}, \mc{n}_k]$. For $k$ large enough $\mc{n}_k=\{0\}$. Let $X_i$ be an orthonormal basis, such that there are integers $k_j$ and 
$X_i \in \mc{n}_j$ for $i \geq k_j$ and $X_i \perp \mc{n}_j$ for $i < k_j$. Let $c_{ij}^k$ be the Maurier-Cartan structural constants with respect to this basis, i.e.

$$[X_i, X_j] = \sum_{k} c_{ij}^k X_k.$$ 
By construction of the basis, $c_{ij}^k \neq 0$ only if $k>j$. Computing the sectional curvatures $K(X_i, K_j)$ using \cite{leftinv} we get for $j \geq i$:

$$K(X_i,X_j) = \frac{1}{4} \sum_{k< j} \left(c_{ik}^j \right)^2 - \frac{3}{4} \sum_{k>j} \left( c_{ij}^k \right)^2.$$

As such, bounding the curvature is, up to a factor, the same as bounding the structural constants. Most of the time we will scale spaces so that sectional curvature satisfies $|K| \leq 2$. Therefore, we will assume throughout this chapter that $||[\cdot, \cdot]||_g \leq C'(n)$, where $C'(n)$ depends only on the dimension $n$. For comparison and later use we will mention the following theorem.

\begin{theorem} (Bieberbach, \cite{bibergromov, ratcliffehyper})\label{flatmanif} There exists a dimension dependent constant $c(n)$ such that the following holds. If $M$ is a complete (not necessarily compact) connected flat Riemannian orbifold of dimension $n$, there exists a discrete group of isometries $\Gamma$ acting on $\R^n$ such that $M = \R^n / \Gamma$. Also, there exists an affine $k$-dimensional subspace $O \subset \R^n$ which is invariant under $\Gamma$, and on which the action of $\Gamma$ is properly discontinuous and co-compact. Further, there exists a subgroup $\Lambda < \Gamma$ with index $[\Gamma : \Lambda]\leq c(n)$ which acts on $O$ by translations. Moreover, $M$ is isometric to a locally flat Riemannian $n-k$-dimensional vector bundle over $O$ (see Definition \ref{quasiflatvect} below).
\end{theorem}
Collapsing theory involves working with certain classes of spaces, which we define here. All notions in this chapter are defined for orbifolds. The versions for manifolds are obtained by a slight variation of terminology. 

\begin{definition} A compact Riemannian orbifold $S$ is called $\epsilon-$almost flat if $|K| \leq 1$ and $\diam(S) \leq \epsilon$.
\end{definition}

Throughout this chapter $\epsilon$ will be small enough, i.e we will only consider $\epsilon(n)'$-almost flat $n$-dimensional manifolds for small enough $\epsilon(n)'$. The choice is dictated by Theorem \ref{almostflat}. For such manifolds a detailed structure theory was developed by Gromov and Ruh (\cite{gromov78almost,ruh1982, buser1981gromov, bibergromov}). Later, Ding observed \cite{ding2011restriction} that with minor modifications the proof as presented in \cite{ruh1982, buser1981gromov} generalizes to orbifolds\footnote{The metric closeness part is not explicit in \cite{ding2011restriction}, but a direct consequence of applying \cite{ruh1982}.}. We state this structure theory here. 

For definitions of orbifolds and some of the other terminology, we refer to the appendix and \cite{thurston2002geometry, ratcliffehyper,Kleiner2014}. For any Lie group $G$ we denote by $\text{Aff}(G)$ the Lie group of affine transformations of $G$, which preserve the flat left-invariant connection $\nabla_{\text{can}}$. The left-invariant connection $\nabla_{\text{can}}$ is defined as the flat connection making left-invariant fields parallel. By $\text{Aut}(G)$, we denote the Lie group of its automorphisms. Throughout $N$ will denote a simply connected nilpotent Lie group.

\begin{theorem}(\cite{ding2011restriction,ruh1982, buser1981gromov} Gromov-Ruh Almost Flat Theorem\footnote{We also recommend consulting the preprint version of the paper by Ding which contains a more complete proof \cite{ding2011restrictionarxiv}. })\label{almostflat} For every $n$ there are constants $\epsilon'(n),c(n)$ such that for any $\epsilon'(n)$-almost flat orbifold $S$ with metric $g$ there exists a simply connected nilpotent Lie group $N$, a left-invariant metric $g_N$ on $N$, and a subgroup $\Gamma < \text{Aff}(N)$ with the following properties.

\begin{itemize}
\item $\Gamma$ acts on $(N,g_N)$ isometrically, co-compactly and properly discontinuously. In particular $N / \Gamma$ is a Riemannian orbifold.
\item We have the bound $|K| \leq 2$ for the sectional curvature of $N$, and $||[\cdot, \cdot]||_{g_N} \leq C'(n)$ for the norm of the bracket.\footnote{For the precise sectional curvature bound, it is necessary to also use a remark in \cite{rong1996fundamental}.}  
\item There is an orbifold diffeomorphism $f\co S \to N / \Gamma$, such that $||f_* g - g_N||_g<\frac{1}{2}$.(The bound here could be replaced by any $\delta<1$, but we fix a choice to reduce the number of parameters in the theorem.)
\item There exists a normal subgroup $\Lambda < \Gamma$, such that $[\Gamma : \Lambda] \leq c(n)$ and $\Lambda$ acts on $N$ by left-translations. In other words, by slight abuse of notation, $\Lambda < N$.
\end{itemize}
\end{theorem}

We remark that we consistently use $N / \Gamma$ to indicate a quotient under a group action even though here $\Gamma$ is acting on $N$ on the left and the translational subgroup $\Lambda$ acts on $N$ by \textit{left} multiplication.  Again $||\cdot||_g$ denotes the operator norm of the various tensors. In the work of \cite{ding2011restriction} and \cite{ruh1982} no explicit bounds are derived for the Lie bracket and curvature, but these can be obtained by using the tools from \cite{rong1996fundamental} and curvature bounds at the beginning of this section. 

In the course of the proof of the main theorem, we will need certain special classes of orbifolds that we will call quasiflat. These will be certain quotients of a nilpotent Lie group $N$ by a group $\Gamma < \text{Aff}(N)$. Each element $\gamma \in \text{Aff}(N)$ can be represented by $\gamma=(a,A)$, where $a \in N$ and $A \in \text{Aut}(N)$. The action of $(a,A)$ on $N$ is given by

$$(a,A)m=aAm$$
where $m \in N$. The component $A$ is called \textit{the holonomy} of $\gamma$. If $A$ is reduced to identity, the operator can be expressed as $(a,I)$, and is called translational. Translational elements $\gamma \in \Gamma$ can be canonically identified by an element in $N$, and we will frequently abuse notation and say $\gamma \in N$. There is a holonomy homomorphism $$h\co \Gamma \to \text{Aff}(N)$$ given by $h(\gamma)=A$.  

The Lie group $N$ will always be equipped with a left-invariant metric $g$. As in the previous section, if $p \in N$ is given and $g \in \text{Isom}(N)$, then we denote

$$|g|_p=d(gp,p).$$

\begin{definition} \label{quasiflatdef} A Riemannian orbifold $S$ is called a \textit{($\epsilon$-)quasiflat} orbifold if $S$ is isometric to a quotient $N / \Gamma$ of a simply connected nilpotent Lie group $(N,g)$ equipped with a left-invariant metric $g$ by a discrete co-compact group of isometries $\Gamma$ satisfying the following. 

\begin{itemize}
\item For the left-invariant metric $g$ the sectional curvature is bounded by $|K| \leq 2$  and $||[\cdot, \cdot]||_g \leq C'(n)$.
\item The group $\Gamma$ acts on the space $N$ isometrically, properly discontinuously and co-compactly.
\item The quotient space $N/ \Gamma$  is a Riemannian orbifold with $\diam(S)=\diam(N/ \Gamma) <  \epsilon$. The space $N$ is referred to as the (orbifold) universal cover and $\Gamma$ is called the orbifold fundamental group.
\end{itemize}

\end{definition} 

We also define a flat orbifold as follows. Another definition, used in Theorem \ref{flatmanif}, requires that the space be complete and that the sectional curvatures vanish. This is, however, equivalent by the same theorem to the following.

\begin{definition} If $\Gamma$ is a group acting co-compactly, properly discontinuously and isometrically on $\R^n$, then we call $S = \R^n / \Gamma$ a flat orbifold.
\end{definition}

The manifold versions of the previous definitions only differ by assuming that the actions are free. Recall, that a group action of $\Gamma$ on a manifold $M$ is \textit{free} if the isotropy group is trivial for any $p \in M$, i.e. $\Gamma_p=\{e\}$, where $e$ is the identity element of $\Gamma$.

By Theorem \ref{almostflat} any $\epsilon(n)'$-almost flat manifold is bi-Lipschitz to an $\epsilon(n)''$-quasiflat manifold with a slightly different $\epsilon(n)''$. Also, by Theorem \ref{flatmanif} every compact flat manifold is also a flat orbifold in the previous sense. In the following, we will choose $\epsilon(n)$ so small, that any $\epsilon(n)$-quasiflat manifold is $\epsilon(n)'$-almost flat, and so that \eqref{eq:commest} holds. Further, we will assume, by possibly making $\epsilon(n)$ smaller, that $\epsilon(n)<\epsilon(m)$ for $n>m$. 

We first prove a result concerning the algebraic structure of $\Gamma$. This can be thought of as a generalization of Bieberbach theorem (see \cite{bibergromov}).

\begin{lemma}\label{finiteindex} Let $c(n)$ be the constant in Gromov-Ruh Almost Flat theorem \ref{almostflat}. Let $S=N / \Gamma$ be a $\epsilon(n)$-quasiflat $n$-dimensional orbifold, then $\Gamma$ has a finite index subgroup $\Lambda$ with $[\Gamma : \Lambda]\leq c(n)$, and $\Lambda = N \cap \Gamma$, i.e. $\Lambda$ acts on $N$ by translations.
\end{lemma}

\noindent \textbf{Proof:} Consider the left invariant metric $g$ of $N$. Since a $\epsilon(n)$-quasiflat space is also $\epsilon(n)'$-almost flat, by Theorem \ref{almostflat} we can construct a new group $\Gamma'$, and a nilpotent Lie group $N'$ such that $S$ is orbifold-diffeomorphic to $N' / \Gamma'$. The proof also gives that $\Gamma'$ acts by affine transformations on $N'$ and has a finite index normal subgroup $\Lambda' \triangleleft \Gamma'$, with $[\Gamma' : \Lambda'] \leq c(n)$, and $\Lambda'$ acts on $N'$ by translations.

Both $N'$ and $N$ constitute orbifold universal covers in the sense of \cite[Chapter 13]{thurston2002geometry}. Thus, there is a diffeomorphism $f\co N' \to N $ conjugating the action of $\Gamma'$ to that of $\Gamma$. Denote the induced homomorphism by $f_* \co \Gamma' \to \Gamma$.

By \cite[Theorem 2]{auslander1960bieberbach} there is a finite index $\Lambda < \Gamma$, such that $\Lambda$ acts by translations on $N'$. In particular, $\Lambda$ is a co-compact subgroup of $N$. We still need to show that it is possible to choose $\Lambda$ with a universal bound on the index $c(n)$, since \cite{auslander1960bieberbach} does not bound the index. We wish to translate the index bound for $\Lambda'$ to one for $\Lambda$. Applying \cite[Proposition 2]{auslander1960bieberbach} to the nilpotent group $f_*^{-1}(\Lambda')$ we get $f_*^{-1}(\Lambda')<\Lambda$. Thus, $[\Gamma : \Lambda]\leq c(n)$. 

\QED

In the previous section we defined generators for groups.  For lattices in nilpotent Lie groups there are special classes of generators that we introduce here. 

\begin{definition} Let $N$ be a simply connected nilpotent Lie group. A co-compact discrete subgroup $\Lambda < N$ is called a lattice.
\end{definition}

\begin{definition} \label{def:basis} Let $N$ be a $n$-dimensional simply connected nilpotent Lie group and $\Lambda < N$ a lattice. Then a set $\{\gamma_1, \dots, \gamma_n\} \subset \Lambda$ is called a triangular basis for $\Lambda$ if the following properties hold.
\begin{itemize}
\item $\gamma_1, \dots, \gamma_n$ generate $\Lambda$.
\item $\gamma_i = e^{X_i}$, and $X_i$ form a vector space basis for the  Lie algebra $\mc{n}\cong \R^n$ of $N$.
\item For $i<j$ we have $[\gamma_i, \gamma_j] \in \langle \gamma_{1}, \dots, \gamma_{i-1} \rangle$.
\item For $i<j$ $[X_i, X_j] \in \mathrm{span}(X_{1}, \dots, X_{i-1})$.
\end{itemize}

\end{definition}

Here, $e^X$ is the Lie group exponential map. In the case of $N=\R^n$ the previous definition reduces to the standard definition of a triangular basis for a Lattice. An immediate consequence of the previous definition is that the group $N$ can be given co-ordinates by $(t_1, \dots t_n) \to e^{t_1 X_1} \cdots e^{t_n X_n}$. This will be used in the appendix. 

\begin{lemma} \label{lem:basis} Let $S = N / \Gamma$ be an $\epsilon(n)$-quasiflat manifold, and $\Lambda = \Gamma \cap N$. For any $p \in N$ we can generate $\Gamma$ and $\Lambda$ by elements  $\gamma$ with $|\gamma|_p \leq 16c(n)\epsilon(n)$.
\end{lemma}

\noindent \textbf{Proof:} Consider the action of $\Lambda=\Gamma \cap N$. Then $N / \Lambda$ is a Riemannian manifold. Further by Lemma \ref{finiteindex}, the index of $\Lambda$ in $\Gamma$ is at most $c(n)$, so there is a $m$-fold orbifold covering map $\pi \co N / \Lambda \to N / \Gamma$ with $m \leq c(n)$. Since $\diam(N / \Gamma) < \epsilon(n)$ and since $N$ is connected, we have $\diam(N / \Lambda) < 2c(n) \epsilon(n)$.

 Let $\Gamma' = \langle g \in \Gamma | |g|_p \leq 16 \epsilon(n) c(n)\rangle$, and $\Lambda' =  \langle g \in \Lambda | |g|_p \leq 16 c(n)\epsilon(n) \rangle$. By Lemma \ref{small} the space $N / \Gamma'$ is isometric to $N / \Gamma$, and thus $\Gamma=\Gamma'$. Similarly, we can show $\Lambda=\Lambda'$. This completes the proof.

\QED

The main conclusion of the following statement is the existence of a basis $\gamma_1, \dots, \gamma_n$ with several technical properties.  The latter part of the statement concludes the existence of associated subgroups $\Lambda_k, L_k$ and scalars $l_k$ with certain geometric properties. These parameters $l_k$ are the ``collapsing scales'', and part of the conclusion is that they are independent of the chosen base point in a specific sense.

\begin{lemma}(Stratification lemma\label{strat}) Fix arbitrary integers $n\geq 1$ and $l>4$, and denote $L=2l^n$.  For every $\epsilon(n)$-quasiflat orbifold $S = N / \Gamma$ there exists following objects with certain desired properties described below.

\begin{enumerate}
\item A triangular basis $\gamma_1, \dots, \gamma_n$ for $\Lambda$.
\item A natural number $s$ (``the number of collapsing scales'').
\item For every $k=1, \dots, s$, subsets $I_k, J_k \subset \{1, \dots, n\}$.
\item Subgroups $\Lambda_k =  \langle \gamma_i | \gamma_i \in J_k\rangle \subset N$ corresponding to $J_k$, for $k=1, \dots, s$.
\item Connected normal subgroups $L_k \triangleleft N$ corresponding to $\Lambda_k$ for $k=1, \dots, s$.
\item Scalars $l_k=2|\gamma_{i_k}|_e>0$, where $i_k=\max{I_k}$ (``collapsing scales'') for each $k=1,\dots, s$.
\end{enumerate}

The basis $\gamma_i$ and the groups $I_k \subset \{1, \dots, n\}$ for $k=1, \dots,s$ satisfy the following.  \begin{itemize}
\item $|\gamma_i|_e<|\gamma_j|_e$ for $i<j$.
\item For all $k$ and $\gamma_i, \gamma_j \in I_k$:  $|\gamma_i|_e\leq L |\gamma_j|_e$.
\item For all $\gamma_i \in I_s$ and $\gamma_j \in I_t$  and $s<t$: $ l |\gamma_i|_e \leq  |\gamma_j|_e$.
\end{itemize}

We can set $J_k=\cup_{i=1}^k I_k.$ The subgroups $\Lambda_k$ are normal in $\Gamma$, and $\Lambda_k < L_k$ is a co-compact lattice in $L_k$. Further, the subgroups $L_k$ are invariant with respect to the action of the image of the holonomy $h \co \Gamma \to \text{Aut}(N)$. 

Finally, the scalars $l_k$ satisfy for every $p \in N$ and every $k=1, \dots, s$ that 
$$\Lambda_k = \langle \lambda \in \Lambda | |\lambda|_p< l_k \rangle$$ and that $$|\lambda|_p>l \cdot l_k/4$$ for every $\lambda \in \Lambda \setminus \Lambda_k$. Also, for any non-trivial $\lambda \in \Lambda$ we have $|\lambda|_p>l_1/(2L)$.
\end{lemma}

\paragraph{Remark:} The properties of $l_k$ are needed to be able to apply Lemma \ref{loc}.

\textbf{Proof:} The basis is constructed similarly to \cite{buser1981gromov}. Let $\gamma_1$ minimize $|\cdot|_e=|\cdot|$ in $\Lambda$. Set $G_0 = \{e\}$ to be the trivial group. Define the subgroup $G_1$ generated by $\gamma_1$. Proceed to choose $\gamma_2 \in \Lambda \setminus G_1$ as the shortest with respect to $|\cdot|_e$. Define $G_2 = \langle \gamma_1, \gamma_2\rangle$. We proceed inductively to define $\gamma_1, \dots, \gamma_n$. The element $\gamma_i$ is the shortest element with respect to $|\cdot|_e$ in $\Lambda \setminus G_{i-1}$ and $G_i = \langle \gamma_1, \dots, \gamma_i \rangle$. A priori, this process could last more than $n$ steps, but by an argument below the number of $\gamma_i$ agrees with the dimension $n$. Thus, we abuse notation slightly by using the same index.

 Take any $a,b \in N$ such that $|a|_e, |b|_e \leq 16c(n) \epsilon(n)$.  We have  \begin{equation}\label{eq:commest}|[a,b]|_e \leq 16 C'(n)c(n) \epsilon(n) \min\{ |a|_e, |b|_e \}.\end{equation} This estimate follows either from curvature estimates such as in \cite{buser1981gromov}, or directly from the bounds on the structural constants on $N$. Choose $\epsilon(n) < \frac{1}{10^2 c(n)C'(n)}$, so that $|[a,b]|_e < \min\{ |a|_e, |b|_e \}$. 

Because of Lemma \ref{lem:basis} we can choose generators $\alpha_i$ for the translational subgroup $\Lambda < \Gamma$ such that $|\alpha_i|_e<16c(n)\epsilon(n)$. Moreover, by construction we get then that $|\gamma_i| < 16c(n)\epsilon(n)$. By estimate \eqref{eq:commest}, we have $|[\alpha_i, \gamma_i]|_e< |\gamma_i|_e$, and thus $[\alpha_i, \gamma_i] \in G_{i-1}$ by construction. Since $\alpha_i$ generate $\Lambda$, we get for any $\lambda \in \Lambda$ that $[\lambda, \gamma_i] \in G_{i-1}$. In particular, for $i<j$ we have $[\gamma_i, \gamma_j] \in G_{i-1}$.  

We define a map $f\co (a_1, \dots, a_n) \to \gamma_1^{a_1} \dots \gamma_n^{a_n}$. By induction, and the product relation below, we can show that the map is injective. Using the commutator relations from above and reordering terms, we can show that $f$ is bijective onto $\Lambda$. Let $a=(a_1, \dots, a_n), b=(b_1, \dots, b_n)$. We can define polynomials $p_j$ by the defining relation

$$\gamma_1^{a_1} \dots \gamma_n^{a_n} \times \gamma_1^{b_1} \dots \gamma_n^{b_n} = \gamma_1^{a_1 + b_1 + p_1(a,b)} \dots \gamma_n^{a_n+b_n + p_n(a,b)}.$$

The polynomials are defined by appropriately applying commutator relations. Also, $p_i$ depends only on $a_j,b_j$ for $j < i$. This defines a product on $\Z^n$, which can be extended to $\R^n$ by the previous relation. This defines a nilpotent Lie group $N'$ with a co-compact lattice $\Z^n$. The lattice has a standard basis $e_i=f^{-1}(\gamma_i) \in \Z^n$, and $e_i=e^{X_i}$, where $X'_i=(0,\dots,1,\dots,0) \in T_0 \R^n$ with the $1$ in the i'th position. By a Theorem of Malcev (see \cite{malcevclass} and \cite{chgrfu}) we can extend the homomorphism $f\co \Z^n \to \Lambda$ to a bijective homomorphism $\overline{f}\co N' \to N$. Since the map is bijective and since $X'_i$ form a basis at $T_e N$, then $d\overline{f}(X'_i)=X_i$ will form a basis for $T_{e} N$. Further $\gamma_i = e^{X_i}$, and $[X_i,X_j]=d\overline{f}([X_i',X_j'])$ it is direct to verify that 
$[X_i, X_j] \in \mathrm{span}(X_{1}, \dots, X_{i-1})$ for $j>i$. In particular the basis is triangular.  Moreover, the number of basis elements has to agree the dimension of the space.

Next, we define the sets $I_j$. Define $i_0=1$. Let $i_1$ be the smallest index such that $|\gamma_{i_1+1}|>l|\gamma_{i_1}|$, or if there is none $i_1=n$. Continue inductively defining $i_k$ to be the smallest index bigger than $i_{k-1}$ such that $|\gamma_{i_k+1}|>l |\gamma_{i_k}|$, or if there is no such index $i_k=n$ and we stop. Define $s$ such that $i_s=n$, and $I_k=\{\gamma_i | i_{k-1}<i \leq i_k\}$. The properties in the statement are easy to verify. The $k$ most collapsed directions are collected into $J_k = \bigcup_{i=1}^k I_i$, and we define $\Lambda_k = \langle \gamma_i | \gamma_i \in J_k \rangle$. Define also $\Lambda_0=\{e\}$

By choosing $\epsilon(n)<\frac{1}{10^2 c(n)C'(n) L}$ and applying \eqref{eq:commest}, we can prove that, if $\gamma_i \in I_s, \gamma_j \in I_t$ and $s<t$, then $[\gamma_i, \gamma_j] \in \Lambda_{s-1}$. We remark that by virtue of the construction, if $\gamma_i \in \Lambda_k$, and $\gamma \in \Lambda$ such that $|\gamma|<l|\gamma_i|$, then also $\gamma \in \Lambda_k$.

Since $\gamma_i$ form a basis, we can show that $\Lambda_k$ is normal in $\Lambda$. Next, we want to show that it is also normal in $\Gamma$. Take an element $\sigma \in \Gamma$ with non-trivial holonomy, and a $\gamma_i \in \Lambda_k$ from the basis. Further, represent $\sigma=(t,A)$, where $A\in \text{Aut}(N)$ is in the holonomy and $t \in N $. By left-multiplying with an element in $\Lambda$, we can assume $|t|_e=|\sigma|_e \leq 16c(n)\epsilon(n)$. This follows from the diameter bound for $N / \Lambda$ in Lemma \ref{lem:basis}. Then by normality of $\Lambda$, $\sigma \gamma_i \sigma^{-1} = t (A\gamma_i) t^{-1} \in \Lambda$. Since $A$ is an isometry, $|A\gamma_i|=|\gamma_i|$. This combined with the commutator estimate \eqref{eq:commest} gives  $d(t (A\gamma_i) t^{-1}, A\gamma_i) \leq \frac{1}{10}|\gamma_i|$, and thus $|t (A\gamma_i) t^{-1}| \leq (1+\frac{1}{10})|\gamma_i|$. By the iterative definition of $\Lambda_k$, we must have $ t (A\gamma_i) t^{-1} \in \Lambda_k$. Note that we are assuming $l>2$.

The connected Lie subgroups $L_k$ can be generated by $X_1, \dots, X_{i_k}$, and normality in $\Lambda$ follows since the basis is triangular. It is also easy to see that $L_k$ is normal in $N$. Further, co-compactness of $\Lambda_k$ follows by considering the map $(t_1, \dots, t_{i+k}) \to e^{t_1 X_1} \cdots e^{t_{i_k} X_k}$. By the argument in the previous paragraph, for any element $\sigma=(t,A) \in \Gamma$ with $|\sigma|_e < 16 \epsilon(n)c(n)$ and $\gamma \in \Lambda_k$, we have $tA \gamma t^{-1} \in \Lambda_k \in L_k$. Since $L_k$ is normal in $N$ we have $A\gamma \in L_k$. From this, we can conclude that $A L_k = L_k$, i.e. that $L_k$ are invariant under the holonomy action of $\Gamma$. 

Finally, define $l_k = 2|\gamma_{i_k}|$. Fix $p \in N$. We will show that $\Lambda_k = \langle \lambda \in \Lambda | d(\lambda p, p)< l_k \rangle$ and $|\lambda|_p=d(\lambda p, p)>l/4 l_k$ for every $\lambda \in \Lambda \setminus \Lambda_k$. First of all, we can assume that $d(p, e)< 16 c(n) \epsilon(n)$. This follows again from the diameter bound for $N / \Lambda$ in Lemma \ref{lem:basis}.  By assumption for any $\gamma_i \in \Lambda_k$ we have $|\gamma_i|_e \leq |\gamma_{i_k}|_e < \frac{1}{2}l_k$. Thus by the commutator estimate \eqref{eq:commest},

$$d(\gamma_i p, p)=d(p^{-1}\gamma_i p, e)\leq d(p^{-1}\gamma_i p, \gamma_i) + d(\gamma_i,e) \leq \left(1+\frac{1}{10}\right)|\gamma_{i_k}|_e < l_k.$$

Thus, $\gamma_i \in \langle \lambda \in \Lambda | d(\lambda p, p)< l_k \rangle$, and further by varying $\gamma_i$ and generation we get $\Lambda_k < \langle \lambda \in \Lambda | d(\lambda p, p)< l_k \rangle$. Next if $\gamma \not\in \Lambda_k$, then $|\gamma|_e=d(\gamma, e)>l/2 l_k$ by the construction of $\Lambda_k$. Simple estimates give

$$d(\gamma p, p)=d(p^{-1}\gamma p, e)\geq  d(\gamma,e) - d(p^{-1}\gamma p, \gamma) \geq (1-\frac{1}{10})|\gamma|_e \geq l/4 l_k.$$

In particular, for every element such that $d(\lambda p, p)<l_k<l/4 l_k$ we must have $\lambda \in \Lambda_k$. Thus, $\Lambda_k = \langle \lambda \in \Lambda | d(\lambda p, p)< l_k \rangle$ and for any $\gamma \not\in \Lambda_k$ we have $|\lambda|_p>l/4 l_k$. The final conclusion, that $d(\lambda p, p)>l_1/(2L)$ for all $\lambda \in \Lambda$, follows from $|\lambda|_e>l_1/L$ which follows from the construction of $\gamma_i$.

\QED

It follows from the previous lemma that the scales $l_k$ and subgroups $\Lambda_k, L_k$ are essentially independent of the chosen base point. Further, they correspond to ``local groups'' at scales $l_k$.

\begin{definition} \label{shortbasis} Let $S=N / \Gamma$ be an $\epsilon(n)$-quasiflat orbifold. If $l=400c(n)$ and $L=2l^n$, then call a basis $\gamma_1, \dots, \gamma_n$ satisfying the previous lemma a short basis for $\Lambda$. The grouping $I_1, \dots, I_s$ is called canonical grouping and the number $s$ is called the number of collapsing scales. If such a basis exists, we say that $S$ has $s$ collapsing scales. For $k=1, \dots, s$, the sets $J_k$ are called the $k$'th most collapsed scales.  
\end{definition}

\begin{definition}\label{quasiflatvect} Let $S$ be a complete Riemannian orbifold and $L = V(S)$ an orbivector bundle over $S$. Then $L$ is called a ($k$-dimensional) \textit{locally flat (Riemannian) orbivector bundle} if the orbivector bundle $V(S)$ admits a bundle-metric and a metric connection such that it is locally flat (and such that the fibers have dimension $k$). In the case where $V$ and $S$ are complete Riemannian manifolds, we call $V$ a locally flat Riemannian vector bundle over $S$.
\end{definition}

A bundle metric is a metric on the fibers of $S$ which is parallel with respect to the connection. For the notion of an orbivector bundle one may consult \cite{Kleiner2014} and the appendix. The notions of a bundle metric and metric connection are naturally extended from vector bundles to orbivector bundles by considering local charts. Note that we permit the case when $S=\{\text{pt}\}$ is a single point (i.e. a zero-dimensional orbifold), and $L=\R^k / \Gamma$, where $\Gamma < O(k)$ is a finite group.  

Because the monodromy action\footnote{Monodromy may also be defined for orbifold. A similar discussion is used to prove that flat orbifolds are good in \cite{ratcliffehyper}.} is by isometries, and the bundle is locally trivial, there is a canonical Riemannian metric and distance function on $L$ (and not just on the fibers). To understand the geometry of $L$ we are reduced to understanding the geometry of $S$ and that of the monodromy action. These metrics are somewhat easier to visualize in the case where $S$ is a compact flat manifold, and $V$ is a locally flat vector bundle over such a base. In this case $L$ is given by a quotient of Euclidean space by an action of a free, discrete and isometric group action, which may fail to be co-compact.

In our case, all the orbifolds $L=V(S)$ will arise as group quotients. By assumption, $S$ will be $\epsilon(n)$-quasiflat. Thus, $S$ has an orbifold universal cover $N$, which is a simply connected nilpotent Lie group. Further, $S$ is given by the quotient $N / \Gamma$, where $\Gamma$ is the orbifold fundamental group, which acts discretely, co-compactly and by isometries on $N$. The orbivector bundle $L$ has a pull-back vector bundle $L'$ on $N$. Since $N$ is simply connected, and $L$ was assumed locally flat, we have a flat pull-back connection on $L'$ making it trivial. In particular, we have, by choosing an orthonormal frame on some fiber, a natural trivialization $L'=N \times \R^k$, where $k$ is the dimension of the fiber. The pull-back metric on $L'$ is simply a product metric, and the distance is given by the product structure. From this perspective, the distance on $L$ is given by a quotient metric on $L' / \Gamma$, where $\Gamma$ acts on the $\R^k$ factor via a monodromy action. Note, all these vector bundles are quotients and will be studied as such. However, they are not \emph{compact}, because the action of $\Gamma$ on the space $N \times \R^k$ is not co-compact.

The monodromy action is given by a homomorphism $\Gamma \to O(n)$, where $\Gamma$ corresponds to the orbifold fundamental group. The rest of this section collects some necessary facts about these group representations. Note that in our case $\Gamma$ is ``virtually'' nilpotent, and its action on $\R^n$ will turn out to be ``virtually'' abelian. By ``virtually'' we mean that there is a finite index subgroup with the desired property, and that the index can be controlled in terms of the dimension $n$.

\begin{definition} Let $H$ be an abelian subgroup of $O(n)$. A decomposition $\R^n = \bigoplus V_i$ is called the \textit{canonical decomposition} and the subspaces $V_0, \dots, V_k$ are called \textit{canonical invariant subspaces} if they satisfy the following properties.

\begin{itemize}
\item Each subspace has positive dimension, except possibly $V_0$.
\item For every $h \in H$, and every $j$ we have  $hV_j = V_j$.
\item There is one $V_j$ such that for all $h \in H$ we have $hv=v$ for all $v \in V_j$ (i.e. each $h \in H$ acts trivially on this subspace). We will always denote this subset $V_0$ possibly allowing $V_0=\{0\}$
\item On each other subspace $V_j$ either one of the following holds.
\begin{enumerate}
\item For every  $h \in H$  and every $v \in V_j$ either $hv=v$ or $hv=-v$
\item $V_j$ is even dimensional and we can associate isometrically complex co-ordinates $(z_1, \dots, z_n)$, such that for any $h \in H$ the action is given by  $h((z_1, \dots, z_n))=\left(e^{\theta_j(h) i}z_1, \dots, e^{\theta_j(h) i}z_n\right)$ for some $\theta_j(h)>0$. 
\end{enumerate}
\item The subspaces are maximal with respect to these properties and pairwise orthogonal.
\end{itemize} 

If $A \in O(n)$ is an isometry, the canonical decomposition associated to it is the one associated to the cyclic subgroup generated by $A$.
\end{definition}

Every abelian isometric action has a unique canonical decomposition. For an individual matrix, this should be thought of as a decomposition corresponding to the factorization of its minimal polynomial. Using this terminology we prove a minor modification of the Jordan lemma.

\begin{lemma}\label{jordan} (Modified Jordan Theorem) There is a constant $c''(n)$ depending on the dimension with the following property. Let $\Gamma$ be a group with a nilpotent normal subgroup $\Lambda$ with index $[\Gamma:\Lambda]<\infty$ and $h\co \Gamma \to O(n)$ a homomorphism. Then $h(\Gamma)$ has a normal abelian subgroup of index at most $c''(n)$.
\end{lemma}

\noindent \textbf{Proof:} Let $L = \overline{h(\Gamma)}$ be the closure of the image of $h$.  The subgroup $L$ is Lie subgroup of $O(n)$. It is either finite or has non-trivial identity component $T$. In the first case, we may proceed as in the standard Jordan Theorem \cite{tao2014hilbert}. In the latter case, since $\Gamma$ has a finite index nilpotent subgroup, $T \subset O(n)$ must be a compact nilpotent group. This is only possible if it is abelian, i.e. a torus. The torus $T$ corresponds to a decomposition of $\R^n$ into canonical invariant subspaces $V_j$. Since $T$ is a continuous group, except possibly for one $V_0$, every $V_j$ for $j=1,\dots,s$ will be even-dimensional, and the action will correspond to a Hopf-action $\theta(z_1, \dots , z_k) \to (e^{i\theta}z_1, \dots , e^{i \theta}z_k)$. On $V_0$ the action is trivial.

Since $T$ is normal in $L$, the decomposition to $V_j$'s corresponds to a homomorphism 
$$\phi = (\phi_1,\dots, \phi_k)\co L \to \prod_{j=1}^s O(V_j).$$ Also, the action of $L$ on each subspace conjugates the Hopf-action on each $V_i$ to itself (which agrees with the action of $T$). The infinitesimal generator in $\mc{o}(V_j)$ corresponding to the Hopf-action is $J$ which maps $J(z_1,\dots, z_k)=(iz_1, \dots, iz_k)$, if $\text{dim}(V_j)=2k$. In particular for any $g \in \phi(L)$ we have $gJg^{-1}=J$ or $gJg^{-1}=-J$ (because it preserves the action), which means that $g$ either preserves the complex structure, or it composed with the conjugate map does so. In other words, $g$ is either symplectic or a product of a symplectic matrix with a matrix $B$ giving the conjugation on $V_j$. Thus, $\phi_j(L) \subset U(V_j) \cup \overline{U(V_j)}$, where $U(V_j)$ is the symplectic group of isometries and is isomorphic to the unitary group of $k$ by $k$ matrices, and $\overline{U(V_j)}$ is obtained by composing such a matrix with $B$.

Consider $\phi(L) \cap \prod_j (SU(V_j) \cup \overline{SU(V_j)}) = K$, where $SU(V)$ and $\overline{SU(N)}$ are those matrices in $U(V_j)$ with unit determinant (using the $k$ by $k$ complex co-coordinates), or their product with $B$. This group must be finite because every left coset of $T$ contains a finite number of elements in $\prod_j SU(V_j)$. Here, we used that the group $T$ acts on $V_j$ via a Hopf-action $\theta(z_1, \dots , z_n) \to (e^{i\theta}z_1,\dots, e^{i \theta}z_n)$, which has determinant $e^{i\text{dim}(V_j)\theta}$. Also, $e^{i\text{dim}(V_j)\theta}=1$ only when $\theta=2\pi n/\text{dim}(V_j)$. Since $K$ is finite, conjugation by $T$ will leave it fixed, and so any element of $T$ commutes with any element in $K$.  Further, $\phi(L)=KT$ by the previous discussion.

By the standard Jordan theorem, the group $K$ has a finite index abelian normal subgroup $H \triangleleft K$ of index at most $d(n) \leq O((n+1)!)$. Define the abelian group $G=HT$, which is generated by $H$ and $T$. We have

$$[\phi(L):G]=[KT:HT]\leq [K:H].$$
Since we have an effective bound for the index $[K:H] \leq d(n)$, we get a bound for $[\phi(L):G] \leq d(n)$. It is also easy to see that $G$ is normal in $\phi(L)$.

\QED 

\textbf{Remark:}  The optimal bound for the index in the classical Jordan theorem is $(n+1)!$, which is valid for $n \geq 71$ \cite{collins2007jordan}. Thus, from the previous argument we get the optimal bound $[\Gamma : \Lambda] \leq c''(n) \leq (n+1)!$ for $n \geq 71$. Below, we will set $c(n)$ to be the larger of $c(n)''$ in Lemma \ref{jordan} and $c(n)$ in Theorem \ref{almostflat}.

\section{Construction of embeddings}

The proof will proceed by increasing generalities. First we discuss the case of flat manifolds, then the case of flat orbifolds, and finally the case of quasiflat orbifolds. At the end, we use Theorems \ref{thm:fukayamanifold} and \ref{thm:fukayaorbifold} to obtain the proof for Theorem \ref{thm:bilip}.

\subsection{Flat Manifolds}
In this section we construct embeddings for flat manifolds by proving the following.

\begin{theorem}\label{flatmanifemb} Every connected complete flat manifold $M$ of dimension $n$ admits a bi-Lipschitz embedding $f\co M \to \R^N$ with distortion less than $D(n)$ and dimension of the image $N \leq N(n)$.
\end{theorem}

We will prove this theorem by induction. The  somewhat complicated induction argument is easier to understand through examples for which we refer to Section \ref{examples}. While it is true that the proof could be directly done for orbifolds, the author believes it is easier to start with the manifold context and later highlight the small differences with the orbifold case.

The goal is to use embeddings at different scales and reveal ever simpler vector bundles at smaller scales. Simplicity is measured in the total dimensionality of the space or the dimensionality of the base space of a vector bundle. This is also a frequent idea of collapsing theory, and as such the proof can be seen as a refinement of collapsing theory to a family of explicit spaces. To avoid lengthy phrases, we use $V-(n,d)$ to abbreviate the induction statement, where $d \geq 0$ and $n > 0$. The case $V-(n,0)$ covers \emph{compact} cases of the theorem, and $V-(n,d)$ for $d>0$ covers the \emph{non-compact} cases, which are vector bundles. Bieberbach's theorem \ref{flatmanif} can be used to show that any connected complete flat manifold is covered by either of these two cases.  The notation $V-(n,0)$ is justified by imagining compact manifolds as 0-dimensional vector bundles over themselves.

Next, we define these statements more carefully. The case $d=0$ is treated in a special way.

\begin{itemize}
\item $V-(n,0):$ There exist constants $L(n,0)$ and $N(n,0)$ such that the following holds. Assume $M$ is an arbitrary $n$-dimensional connected compact flat manifold. Then, there exists a bi-Lipschitz embedding $f\co M \to \R^N$ with distortion $L \leq L(n,0)$ and target dimension $N \leq N(n,0)$ depending only on $n$.
\item $V-(n,d):$ There exist constants $L(n,d)$ and $N(n,d)$ such that the following holds. Assume $V$ is an arbitrary $d$-dimensional connected locally flat Riemannian vector bundle over a compact flat manifold $M$ of dimension $n$. Then, there exists a bi-Lipschitz embedding $f\co V \to \R^N$ with distortion $L \leq L(n,d)$ and target dimension $N \leq N(n,d)$ depending only on $n,d$.
\end{itemize}

In the proof below we will perform a double induction to prove $V-(n,d)$, which is based on reducing the embedding problem for $V-(n,d)$ to embedding instances of $V-(n-k,k+d)$ for some $k>0$ or of $V-(a,b)$ with $a+b<n+d$. For $V-(n,0)$ the reduction is similar to that of a lens space, and for $V-(n,d)$ the reduction is analogous to that for holonomy-bundles (see Section \ref{examples}).

The main steps of these reductions are choosing a length scale $\delta>0$, applying a stratification of a lattice from Lemma \ref{strat}, followed by an application of the local group Lemma \ref{loc} to expose the structure of a local group at the scale $\delta$, and then applying Lemma \ref{small} to describe the metric structure at those scales. The metric structure at a small but definite scale is always that of a vector bundle. The base of such a vector bundle is an embedded manifold of smaller dimension. This embedded manifold is found by finding submanifolds $\R^k$ in the universal cover - that is $\R^n$ - which remain invariant under the action of a local group. Here, the averaging from  Lemma \ref{fixedpoint} plays a crucial role. Once the embeddings are constructed at scale $\delta$, using the induction hypothesis, the embedding of the space is obtained by using Lemma \ref{doublingarg}. 

\paragraph{Proof or Theorem \ref{flatmanifemb}:} We first give a description of the induction argument, and then move onto defining the relevant notation and proving the induction claims.

\paragraph{Summary of induction argument:} We need to embed all complete connected flat manifolds. By Biberbach's theorem \ref{flatmanif} every complete connected flat manifold can be represented as a quotient by a discrete group of isometries. In particular, any such manifold is either isometric to $\R^n$, a quotient of $\R^n$ by a co-compact group action or a locally flat vector bundle over such a base. In the first case, the embedding is given by the identity mapping. Otherwise, either the space can be represented as $M=\R^n / \Gamma$, or $V= \R^n \times \R^d / \Gamma$, where $\Gamma$ is some discrete Lie group which acts co-compactly, properly discontinuously and by isometries on $\R^n$. Additionally, for $V$, the action of $\Gamma$ leaves $\R^n \times \{0\}$ invariant, and acts by orthogonal tranformations on the fibers $\R^d$. In the latter case, the action is \emph{not} co-compact on all of $\R^n \times \R^d$, and the quotient will be a non-compact space. The case $V-(n,0)$ covers the problem of embedding $M$, and $V-(n,d)$ covers the problem of embedding $V$. 

Thus, we need to show $V-(n,d)$ for $n \geq 1, d \geq 0$. First, we induct on $n+d$. If $n+d=1$, then we only have the case $V-(1,0)$, which is the base case for our induction. 

Next, assume the statement has been shown for $n+d \leq m$, and we prove it for $V-(n,d)$ with $n+d = m+1$, that is $V-(k, m+1-k)$ for $k=1, \dots, m+1$. We show this by induction on $k$. The base case is $V-(1,m)$ which is reduced to $V-(a,b)$ with $a+b\leq m$. Note $S=1$ in this case. 

Assume now that we have shown $V-(k,m+1-k)$ for all $k=1, \dots, k'-1$ with $k'<m+1$ and show it for $k=k'$. If $k'< m+1$, then the third reduction below reduces the embedding problem to instances of $V-(a,b)$ with either $a<k'$ or $a+b\leq m$. In both cases we have reduced the embedding problem to cases we have assumed to be true. For $k'=m+1$ we want to prove $V-(m+1,0)$, which is reduced to embedding instance of $V-(a,b)$ with $a<m+1$. By the induction principle, we can assume these are already embedded, and thus the statement follows.

\paragraph{Defining variables:}
 We will use $\Gamma$ to denote the group acting on $\R^n \times \R^d$ and $\R^n$ as described above. By Biberbach's theorem, since the action of $\Gamma$ on $\R^n$ is co-compact, there is a finite index lattice $\Lambda = \Gamma \cap \R^n $ acting by translations on $\R^n$ with index $[\Gamma : \Lambda]\leq c(n)$. Denote by $0$ the identity element of $\R^n$ and use the subadditive norm $|\gamma|_0=|\gamma|=d(\gamma 0, \gamma)$ on $\Gamma$. Denote by $e$ the identity element of $\Gamma$. Denote by $D = \diam(M)=\diam(\R^n / \Gamma)$. Apply Lemma \ref{strat} to choose the following objects. Choose a short basis $\gamma_1, \dots, \gamma_n$ for $\Lambda$, and the canonical grouping $I_1, \dots, I_S$ with parameters $l=400c(n)$ and $L=2l^n$ (see Definition \ref{shortbasis}). Let $S$ be the number of collapsed scaled. By $J_k$ for $k=1, \dots, S$ we denote the $k$'th most collapsed scales. By the proof for Lemma \ref{strat}, these correspond to normal subgroups $\Lambda_k < \Lambda$, which are normal in $\Gamma$. Also, they correspond to connected normal Lie-subgroups $L_k<\R^n$ which are invariant under the holonomy action of $\Gamma$, and such that $\Lambda_k < L_k$ are co-compact. Denote by $l_k$ the collapsing scales given by Lemma \ref{strat}.

Points in $\R^n$ are denoted by lower case letters such as $p$ and points in the quotient space by $[p]$. Cosets with respect to local groups will be denoted by $[p]_r$. We will now proceed by induction on $n$ and $n+d$ to show $V-(n,d)$ for all $n$.

\paragraph{Base case $V-(1,0)$:} In this case the flat manifold is one-dimensional and isometric to $S^1$. Such a manifold is trivial to embed.

\paragraph{Case $V-(n,0)$ when $S=1$:} From the definition of $l_1=2|\gamma_n|_0$ in Lemma \ref{strat}, and  Lemma \ref{nilpdiam} we obtain  $\diam(M)\leq C(n)l_1$. Let $p \in \R^n$ be arbitrary. Set $\delta=l_1/(320c(n)L)$. Define the group $\Gamma_\delta=\langle g \in \Gamma | d(gp,p)<8 \delta \rangle$. By Lemma \ref{small} for any $[p] \in M$ we have $B_{[p]}(\delta) \subset M$ is isometric to $B_{[p]_\delta}(\delta) \subset \R^n / \Gamma_\delta(p)$. We next describe the structure of the group $\Gamma_\delta(p)$. 

Apply Lemma \ref{loc} with $l=8\delta$. When we define $\Lambda_0=\langle \lambda \in \Lambda | d(\lambda p,p)<8\delta \rangle = \{e\}$, and observe that by Lemma \ref{strat} for any $\lambda \in \Lambda$ we have $|\lambda|_e > l_1/(2L) \geq 80c(n) \delta$, we obtain a short exact sequence

$$0 \to 0 \to \Gamma_\delta \to H_0 \to 0.$$

The group $H_0$ is finite and thus also $\Gamma_\delta$ is finite. Since $M$ is a manifold, $\Gamma_\delta$ acts freely on $\R^n$. Applying Lemma \ref{fixedpoint}, we get that, if $\Gamma_\delta$ were non-trivial, then it would have a fixed point. Since this would contradict that the action of $\Gamma_\delta$ is free, we must have $\Gamma_\delta$ is trivial. Thus,  $B_{[p]_\delta}(\delta)\subset \R^n / \Gamma_\delta(p)$ is isometric to a Euclidean ball of the same radius $B_0(\delta) \subset \R^n$. Thus, each ball $B_{[p]}(\delta)$ is isometric to a subset of $\R^n$, and obviously admits a bi-Lipschitz embedding to $\R^n$. Also, $\diam(M) / \delta \leq  320C(n)c(n) L$. Since $p$ was arbitrary, Lemma \ref{doublingarg} implies that $M$ also possesses an embedding.  

\paragraph{Reduction $V - (n, 0)$ to $V-(n - k, k)$ for some $n > k \geq 1$ when $S>1$}: The gain here is that the dimension of the new base space $n-k$ is smaller than $n$. Consider $\delta =\max\left(\frac{l_S}{200Lc(n)}, 2 l_{S-1}\right)$. Note that $\diam(M)<C(n) l_S \leq 200Lc(n) C(n) \delta$ by Lemma \ref{nilpdiam}. Thus, similarly to the previous case, it is sufficient by Lemma \ref{doublingarg} to embed each ball $B_{[p]}(\delta)$. Next fix $[p] \in M$, and its representative $p \in \R^n$. We will show that $B_{[p]}(\delta)$ is isometric to a subset of a vector bundle over a $(n-k)$-dimensional flat manifold for some $k$.

Define the local group $\Gamma_\delta=\langle g \in \Gamma | d(gp,p)<8 \delta \rangle$, and using Lemma \ref{small} we get $B_{[p]}(\delta)\subset M$ is isometric to $B_{[p]_\delta}(\delta) \subset \R^n  / \Gamma_\delta$. Define the subgroup $\Lambda_\delta = \langle \lambda \in \Lambda | d(\lambda p,p)<8\delta \rangle $. By Lemma \ref{strat} for each $\lambda \not\in \Lambda_{S-1}$ we have $d(\lambda p, p)>l \cdot l_{S-1}/4 > 80c(n)\delta$, and thus we can see that $\Lambda_\delta=\Lambda_{S-1}$. Further, by Lemma \ref{loc} applied with $l=8\delta$ there is a finite group $H_0$ and a short exact sequence

$$0 \to  \Lambda_{S-1} \to \Gamma_\delta \to H_0 \to 0.$$

Consider the connected Lie-subgroup $L_{S-1}$ corresponding to $\Lambda_{S-1}$ via Lemma \ref{strat}. The subgroup $L_{S-1}$ is invariant under $\Gamma_\delta$ and the holonomy. Thus, the action of $\Gamma_\delta$ descends to an action of $H_0$ on the orbit-space $E=\R^n / L_{S_1}$. Represent elements in $E$ as orbits $[x]_{L_{S-1}}$ for $x \in \R^n$. Thus by Lemma \ref{fixedpoint}, there is an orbit $[o]_{L_{S-1}} \in E$ fixed by the action of $H_0$. Let $O = [o]_{L_{S-1}} \subset \R^n$ be the orbit. By construction $O$ is invariant under the group action $\Gamma_\delta$.  Since $O$ is invariant under $\Gamma_\delta$ we have $\R^n / \Gamma_\delta = V'$, where $V'$ is a locally flat Riemannian vector bundle over $O / \Gamma_\delta$ of dimension $n-k$. Also, the action of $\Gamma_\delta$ is co-compact on $O$, because $\Lambda_{S-1}<\Gamma_\delta$. In particular $B_{[p]_\delta}(\delta)$ is isometric to a subset of a $k$-dimensional locally flat Riemannian vector bundle over a $n-k$-dimensional compact connected flat manifold. Thus, the ball $B_{[p]}(\delta)$ can be embedded by $V-(n-k,k)$. 

\paragraph{Reduction $V - (n, d)$ to $V-(n-k, d+k)$ (if $S>1$) and cases of the form $V-(a,b)$ with $a+b<n+d$:}  Denote by $V=\R^n \times \R^d / \Gamma$ the locally flat vector bundle over $M = \R^n \times \{0\}/ \Gamma \subset V$. This reduction guarantees that either the dimension of the base space drops, or the total dimension drops. The proof is similar to the holonomy bundle example discussed in Section \ref{examples}. 

The space $V$ is divided into sets defined by their distance to the zero-section $M$. Take the sets $T_0 = \{ p \in V | d(p, M) \leq 2D \}$ and $T_k = \{ p \in V | 2^{k+1}D \geq d(p, M) \geq 2^{k-1}D\}$ for $k \in \N$. Recall, $D=\diam(M)$. We first show that it is sufficient to construct an embedding $f_k$ for each $T_k$ of bounded dimension and distortion. Suppose such embeddings exist. First scale and translate $f_k$ so that $0 \in \text{Im}(f_k)$ and $f_k$ is $1$-Lipschitz. Consider the sets $D_s = \bigcup_{k=0}^\infty T_{s+4k}$, and define $F_s(x)=f_{s+4k}(x)$ on $T_{s+4k}$. The resulting functions will be Lipschitz on their respective domains $D_s$. By an application of McShane extension Theorem \ref{thm:extension}, we can extend them to all of $V$. A bi-Lipschitz embedding will result from

$$G(x)=(r(x),F_0(x),F_1(x),F_2(x),F_3(x)),$$
where $r(x)=d(x,M)$. To see that this is an embedding, repeat the argument from Section \ref{examples}. Next, construct embeddings $f_i$ for each $T_i$.

First, we consider $T_i$ for $i \geq 1$. Define $K = \{r(x) = 2^i D\}$. It is sufficient to find an embedding for the set $K$ with it's induced metric. The metric space $T_i$ is bi-Lipschitz equivalent to $[2^{i-1}, 2^{i+1}] \times K$, and we can use the product embedding from Lemma \ref{product} to embed it. Therefore, we can focus on embedding $K$. If $d = 1, 2$, $K$ will be compact and flat, and we can apply the case $V-(n + d - 1,0)$ or the case $V-(n + d - 2,0)$ to embed $K$. If $d > 2$, then $K$ is not flat. But it is lower dimensional and has the structure of a sphere bundle. For each $q \in K$ we find a ball with radius proportional to $2^iD$ and which is bi-Lipschitz to a lower dimensional vector bundle. This combined with the doubling Lemma \ref{doublingarg} gives an embedding for $T_i$.  

For a real number $t$ denote by $tS^{d-1}$ the $d-1$-dimensional sphere of radius $t$. The space $K$ can be represented as $K=\R^n \times D2^iS^{d-1} / \Gamma$, where $\Gamma$ acts on $S^{d-1}$ by isometries via the holonomy action $h \co \Gamma \to O(d)$. Any point in $K$ can be represented as an orbit $[(q,v)] \in K$ for an element $(q,v) \in \R^n \times D2^iS^{d-1}$. Fix $\delta=\frac{D2^i}{10^5 c(d)\sqrt{d}}$. We will show that $B_{[(q,v)]}(\delta) \subset K$ is bi-Lipschitz to a subset of locally flat vector bundle $V$ of total dimension $m=n-1$. This will be done by studying the local group $\Gamma_\delta = \langle g \in \Gamma | d(g(q,v), (q,v))\leq 8 \delta \rangle$, and observing by Lemma \ref{small} that $B_{[q]}(\delta)$ is isometric to $B_{[(q,v)]_\delta}(\delta) \subset \R^n \times D2^iS^{d-1} / \Gamma_\delta $ and describing in detail the action of $\Gamma_\delta$.

By the modified Jordan Theorem \ref{jordan}, the group $h(\Gamma)$ has a normal abelian subgroup $L$ of index at most $c(d)$, and $h^{-1}(L) = \Lambda_L$ is also normal and has finite index at most $c(d)$. Decompose $\R^d$ into canonical invariant subspaces $V_j$ corresponding to $L$. By uniqueness, the subspaces $V_j$ are invariant under the action of $\Gamma$. Decompose $v$ into its components $v_j \in V_j$. Consider the induced norm $||\cdot||$ on $\R^d$. There is an index $j$ such that $||v_j|| \geq \frac{2^i D}{3\sqrt{d}}$. To reduce clutter, re-index so that $j = 1$ and denote $\xi_1 = D2^i v_1 / ||v_1||$. Denote by $h_1 \co \Gamma \to O(V_1)$ the action of $\Gamma$ on $V_1$. Denote by $d_{S^{d-1}}$ the metric on the sphere $2^iD S^{d-1}$. Then 
$$d_{S^{d-1}}(v,\xi_1) \leq (2^i D) \cos^{-1}\left( \frac{1}{3\sqrt{d}} \right) \leq 2^iD \left(\pi/2 - \frac{1}{3\sqrt{d}} \right).$$ 
For each generator $g$ of $\Gamma_\delta$ we have $d(g(q,v),(q,v))\leq 8 \delta$. Then also 
\begin{equation} \label{sphereest}
d_{S^{d-1}}(g\xi_1,\xi_1) \leq 24\delta\sqrt{d}.
\end{equation} 
Denote by $\Lambda_0 = \Gamma_\delta \cap \Lambda_L$. Clearly $[\Gamma_\delta : \Lambda_0]\leq c(d)$. Also, define $\Lambda_\delta = \langle g \in \Lambda_L | d_{S^{d-1}}(g\xi_1, \xi_1)\leq 24\delta\sqrt{d} \rangle$.

Either the action of $L$, and of $\Lambda_L$, on the sub-space $V_1$ will consist of reflections and identity transformations, or it will be a Hopf-type action.  Consider the first case. We have $24c(d)\delta\sqrt{d} < D2^i \pi/4$, and thus the action of $\Lambda_\delta$ on $V_1$ is trivial. Consider the semi-norm $|g|'_{\xi_1}=d_{S^{d-1}}(g\xi_1, \xi_1)$, and the action of the groups $h_1(\Lambda_L)$ and $h_1(\Gamma)$ on $O(V_1)$ with the parameter $l=24\delta \sqrt{d}$. Applying Lemma \ref{loc}, and using \eqref{sphereest} we get a short exact sequence

$$0 \to  h_1(\Lambda_{\delta}) \to h_1(\Gamma_\delta) \to H_0 \to 0,$$
where $H_0$ is a finite group of order at most $c(d)$.

Since the action of $\Lambda_\delta$ on $V_1$ is trivial, we get that $h_1(\Gamma_\delta)$ is finite. By equation \eqref{sphereest} and the product bound of Lemma \ref{loc}, we get for any $g \in h_1(\Gamma_\delta)$ that $d_{S^{d-1}}(g\xi_1,\xi_1)<40c(d)\delta\sqrt{d} < D2^i\pi/(100\sqrt{d}).$ By the fixed point Lemma \ref{fixedpoint} applied to the sphere $D2^i S^{d-1}$ we find a $w \in V_1 \cap 2^i D S^{d-1}$ fixed by the action of $h_1(\Gamma_\delta)$. Moreover, this vector is fixed by the action of $h(\Gamma_\delta)$. Also, $d_{S^{d-1}}(\xi_1,w) \leq D2^i \pi/(100\sqrt{d})$. Thus,

\begin{eqnarray}
d_{S^{d-1}}(w,v)  &\leq & d_{S^{d-1}}(v,\xi_1) + d_{S^{d-1}}(\xi_1,w) \label{invest}\\
                                            &\leq & D2^i \left(\pi/(100\sqrt{d})+\pi/2 - \frac{1}{3\sqrt{d}}\right) \nonumber \\
                                            &\leq& D2^i\left(\pi/2 - \frac{1}{10\sqrt{d}}\right). \nonumber
\end{eqnarray}

 The submanifold $\R^n \times \{w\}$ is invariant under the action of $\Gamma_\delta$. Since the action of $\Gamma_\delta$ on $\R^n$ is also co-compact, the manifold $M_\delta = \R^n \times \{w\} / \Gamma_\delta$ is a compact flat manifold. Fix $R= D2^i(\pi/2 - \frac{1}{20\sqrt{d}})$. The normal exponential map gives a map of the $R$-ball $B'$ in the tangent space of $D2^i S^{d-1}$ at $w$ onto $B_w(R) \subset D2^i S^{d-1}$. Call this normal exponential map $e_w\co \R^{d-1} \to 2^i D S^{d-1}$. When $B'$ is equipped with the induced Euclidean metric, this map becomes a bi-Lipschitz map $B' \to B_w(R)$. (See the calculation for lens spaces in section \ref{examples}, which can be generalized to this setting using polar co-ordinates). Taking products with $\R^n$ we get a bi-Lipschitz map $F\co \R^n \times B' \to N_R(\R^n \times \{w\})$, where $N_R(\R^n \times \{w\})$ is the $R$-tubular neighborhood of $\R^n \times \{w\}$ in $\R^n \times 2^i D S^{d-1}$.

The action of $\Gamma_\delta$ induces an action on the normal bundle of the subset $\R^n \times \{w\}$, and thus on $\R^n \times B'$. Further, this action of $\Gamma_\delta$ on $\R^n \times B' \subset \R^n \times \R^{d-1}$ commutes with the action on $\R^n \times 2^i D S^{d-1}$ via the map $F$.  Thus, $F$ induced a bi-Lipschitz map $f \co \R^n \times B' / \Gamma_\delta \to N_{R}(\R^n \times \{w\})$. Use \eqref{invest} to get that  $B_{[p,v]_r}(\delta)$ is contained in $N_R(\R^n \times \{w\}) / \Gamma_\delta$. Thus, the restricted inverse $f^{-1}\co B_{[(q,v)]_\delta}(\delta)  \to \R^n \times B' / \Gamma_\delta$ is a bi-Lipschitz map. The image $\R^n \times B' / \Gamma_\delta$ is an isometric subset of $\R^n \times \R^{d-1} / \Gamma_\delta$ which is a $d-1$-dimensional locally flat vector bundle over $M_\delta$, which is a compact flat manifold of dimension $n$. Thus, it can be embedded using the statement $V-(n, d-1)$, which has smaller total dimension.

The other case is when $V_1$ is a $2k$-dimensional rotational subspace where the action of $L$ is given by the Hopf-action. We can define $\Lambda_L$ and $\Lambda_0=\Lambda_L \cap \Gamma_\delta$ similarly to before. This time the action of $\Lambda_0$ on $\R^n$ will not necessary be trivial, but will be a Hopf-type action. For any point $\xi \in 2^i DS^{2k-1} \cap V_1$ its $\Lambda_0$-orbit is contained in a unique Hopf-fiber $S_\xi^1 \subset 2^i D S^{2k-1}$. The action of $\Gamma_\delta$ conjugates the Hopf-action to itself, and thus any orbit $S_{\xi}^1$ is mapped to another orbit $S_{\psi}^1$. We realize this orbit space as the complex projective space of diameter $2^i D$, i.e. $2^{i}D\mathbb{C P}^{(k-2)/2}$, and the action of $\Gamma_\delta$ descends to $2^{i}D\mathbb{C P}^{(k-2)/2}$. Denote by $S^1_{\xi_1}$ the orbit through $\xi_1$.

Apply the argument in Lemma \ref{loc} with $\Lambda=\Lambda_0$, $\Gamma=\Gamma_\delta$, $l=24\delta\sqrt{d}$, and the seminorm $|g|'=d(gS^1_{\xi_1}, S^1_{\xi_1})$. This gives that $\Lambda_0 \cap \Gamma_\delta$ has finite index in $\Gamma_\delta$ and any coset of $\Lambda_0 $ in $\Gamma_\delta$ can be represented by an element $g \in \Gamma_\delta$ such that

\begin{equation}\label{orbitest}
d(gS^1_{\xi_1}, S^1_{\xi_1}) \leq d(g\xi_1, \xi_1) \leq 24\delta\sqrt{d}c(d)
\end{equation}

Thus, the action of $\Gamma_\delta$ on $2^{i}D\mathbb{C P}^{(k-2)/2}$ is by a finite group $H_0=\Gamma_0 / \Lambda_0$.  By \eqref{orbitest}, we have for any $h \in H_0$ that for the orbit $S^1_{\xi_1}$ corresponding to $\xi_1$:
$$d(hS^1_{\xi_1}, S^1_{\xi_1}) \leq d(g\xi_1, \xi_1) \leq 24c(d)\delta\sqrt{d}< D2^i \pi/4.$$ 

\begin{figure}[h!]
  \centering
    \includegraphics[width=\textwidth]{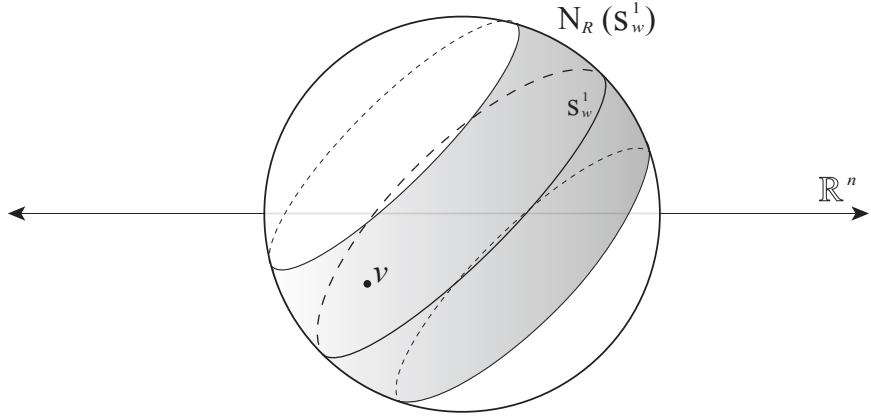}
    \caption{The normal neighborhood of an orbit of the Hopf action is presented in this simplified image, where $S^1_w$ is the invariant fiber, and the shaded region its normal neighborhood. We have simplified the setting by reducing the sphere to a 2-dimensional sphere, and by indicating the $\R^n$ factor abstractly as a line. The point $v$ is depicted as belonging to the normal neighborhood.}
    \label{fig:normals1}
\end{figure}

Lemma \ref{fixedpoint} can be used to give an orbit $S^1_{w}$ for some $w \in 2^i D S^{d-1}$, which is fixed under the action of $H_0$. In particular, the orbit as a subset of $S^{d-1}$ is invariant under $\Gamma_\delta$. Again, by Lemma \ref{fixedpoint} and \eqref{orbitest} we get 
\begin{equation}\label{invest2}d_{S^{d-1}}(S^1_w,v) \leq  D2^i(\pi/2 - \frac{1}{10\sqrt{d}}).
\end{equation} 
This time define $M_\delta = \R^n \times S^1_w /\Gamma_\delta$ and consider it as a subset of $K$. Define $R=D2^i(\pi/2 - \frac{1}{20\sqrt{d}})$. The ball $B_{(q,v)}(\delta) \subset K$ is contained in the $R$-tubular neighborhood $N_R(\R^n \times S^1_w / \Gamma_\delta)$ by \eqref{invest2}. See Figure \ref{fig:normals1} for a simplified image of this normal neighborhood.  The $R$-tubular neighborhood can be identified in a bi-Lipschitz fashion with $\R^n \times S^1_w \times B' / \Gamma_\delta$ for $B'$ a Euclidean ball of radius $R$ and dimension $d-2$. The identification is via the normal exponential map of the submanifold $\R^n \times S^1_w / \Gamma_\delta \subset K$. The quotient $\R^n \times S^1_w \times B' / \Gamma_\delta$ is isometric to a subset of a $d-2$-dimensional locally flat vector bundle $V''$ over $M_\delta$. In particular, we get a map $F\co B_{(q,v)}(\delta) \to V''$. The conclusion thus follows from the statement $V-(n + 1, d - 2)$.

The previous paragraphs give bi-Lipschitz maps for each $T_i$ when $i \geq 1$. For $i = 0$ we need a different argument. First assume $S = 1$. Similar to the case $V-(n,0)$ when $S=1$, we can use Lemma \ref{doublingarg} and embed each ball $B_{[p]}(\delta)$ since for $\delta=l_1/(320c(n)L) \sim D$ we have $\Gamma_\delta = \{e\}$. In other words, the space has a lower bound on its injectivity radius. Also, note $\diam(T_0) \leq 6D$.

 Assume instead that $S > 1$. Consider an arbitrary $B_{[p]}(\delta)$ with $[p] \in T_0$ for $\delta=\max\left(\frac{l_S}{200Lc(n)}, 2 l_{S-1}\right)$. Here $[p]$ is an orbit of the $\Gamma$ action for $p=(p^*, v) \in \R^n \times \R^d$. Since $\diam(T_0) \leq 6D \leq 200Lc(n)C(n) \delta$, we can use the Lemma \ref{doublingarg} to construct embeddings once we can embed all of the balls $B_{[p]}(\delta)$. Again define $\Gamma_\delta = \langle g \in \Gamma | d(gp,p) \leq 8\delta \rangle$. Also, define $\Gamma'_\delta=\langle g \in \Gamma | d(gp^*, p) \leq 8 \delta \rangle$. We have $\Gamma_\delta < \Gamma'_\delta $. Denote their orbits containing $x$ by $[x]_\delta$ and $[x]_\delta'$ respectively. Then, we have $B_{[p]}(\delta) \subset V$ is isometric to $B_{[p]_\delta}(\delta) \subset \R^n \times \R^d / \Gamma_\delta$ and $B_{[p]_\delta'}(\delta) \subset \R^n \times \R^d / \Gamma_\delta'$. Here we used Lemma \ref{small}. 

By the second case in the proof, there is an invariant affine subspace $O \subset \R^n$ of dimension $n-k$ such that $\Gamma'_\delta O \subset O$. In particular $\R^n \times \R^d / \Gamma_\delta'$ is isometric to $O \times \R^{d+k} / \Gamma_\delta'$. The space $O \times \R^{d+k} / \Gamma_\delta'$ is a locally flat vector bundle over $O / \Gamma_\delta$, and since $B_{[p]}(\delta)$ is isometric to a subset of it, we can embed it by the statement $V-(n-k,d+k)$ for some $k>0$.

\QED

\subsection{Flat orbifolds}

We next modify the proof slightly to allow us to embed complete flat orbifolds. 

\begin{theorem} \label{orbifoldflat}Every complete connected flat orbifold $O$ of dimension $n$ admits a bi-Lipschitz embedding into Euclidean space with distortion and dimension depending only on the dimension $n$. Further, every locally flat orbivector bundle over such a base with its natural metric admits such an embedding with distortion and target dimension depending only on the total dimension.
\end{theorem}

\paragraph{Proof:}  We first give a description of the induction argument, and then move onto defining the relevant notation and proving the induction claims.

\paragraph{Summary of induction argument:} The arguments are sufficiently similar to the proof of Theorem \ref{flatmanifemb} that we only indicate the main differences. We will need to modify the induction statements slightly to add the case $V-(0,d)$, which corresponds to an orbivector bundle over a point, i.e. $\R^n / \Gamma$, where $\Gamma$ is a finite group. We assume the same notational conventions as used in the beginning of Theorem \ref{flatmanifemb}. 

\begin{itemize}
\item $V-(n,0):$ There exist constants $L(n,0)$ and $N(n,0)$ such that the following holds. Assume $M$ is an arbitrary connected compact flat orbifold with dimension $n$. Then, there exists a bi-Lipschitz embedding $f\co M \to \R^N$ with distortion $L \leq L(n,0)$ and target dimension $N \leq N(n,0)$ depending only on $n$.
\item $V-(n,d), n\geq 1, d>1:$ There exist constants $L(n,d)$ and $N(n,d)$ such that the following holds. Assume $V$ is an arbitrary $d$-dimensional locally flat Riemannian orbivector bundle over a compact flat orbifold $M$ of dimension $n$. Then, there exists a bi-Lipschitz embedding $f\co V \to \R^N$ with distortion $L \leq L(n,d)$ and target dimension $M \leq N(n,d)$ depending only on $n,d$.
\item $V-(0,d):$ There exist constants $L(0,d)$ and $N(0,d)$ such that the following holds. Assume $V=\R^d / \Gamma$ where $\Gamma < O(d)$ is a finite group. Then, there exists a bi-Lipschitz embedding $f\co V \to \R^N$ with distortion $L \leq L(0,d)$ and target dimension $N \leq N(0,d)$ depending only on $d$.
\end{itemize}

\paragraph{Overview of induction:}
We need to embed all complete connected flat orbifolds. Again, by Biberbach 's theorem every complete connected flat orbifold can be represented as a quotient by a discrete group of isometries \ref{flatmanif}. In particular, any such manifold is either isometric to $\R^n$, a quotient of $\R^n$ by a co-compact quotient or a locally flat vector bundle over such a base. In the first case the embedding is given by the identity mapping. Otherwise, either the space can be represented as $M=\R^n / \Gamma$, or $V= \R^n \times \R^d / \Gamma$, where $\Gamma$ acts co-compactly, properly discontinuously and by isometries on $\R^n$ and by isometries on $\R^d$. The case $V-(n,0)$ covers the problem of embedding $M$, and $V-(n,d)$ covers the problem of embedding $V$. We allow for $n=0$ in these statement. 

The induction works in the same way as for Theorem \ref{flatmanifemb}, except we have new base cases for the induction $V-(0,1)$ and $V-(0,2)$.

\paragraph{Remark on notation:} We will assume the same notation as in Theorem \ref{flatmanifemb}, except for the slight variation that the group actions are not assumed to be free.

\textbf{Special cases $V-(0, 1)$, $V-(0, 2)$ and $V-(1,0)$}: The spaces that result can all be enumerated and checked individually: a ray, two dimensional cone over a circle, a circle, or a cone over an interval. Their bi-Lipschitz embeddings can be explicitly constructed.

The rest of the proof covers various inductive steps similar to Theorem \ref{flatmanifemb}, but including the degenerate one where the base is zero-dimensional. 

\textbf{Reduction $V-(0, n)$ to $V-(n - k, k - 1)$:} Let $\Gamma < O(n)$ be a finite subgroup. The space $M = \R^n / \Gamma$ is isometric to $C(S^{n-1}/\Gamma)$, i.e. the cone over $S^{n-1}/\Gamma$. By Lemma \ref{cone}, it is sufficient to embed $S = S^{n-1}/ \Gamma$. By the Jordan theorem\footnote{In this case $\Lambda$ is the trivial group, and the statement is the classical Jordan theorem \cite{tao2014hilbert}} \ref{jordan} we have a finite index normal subgroup $K \triangleleft \Gamma$ which is abelian with index bounded by $[\Gamma : K] \leq c(n)$.  Decompose $\R^n$ to canonical invariant subspaces $V_i$ for $K$. 

 Let $\delta = \frac{1}{10^5 c(n)\sqrt{n}}$, and fix an arbitrary $[v] \in S^{n-1} / \Gamma$ which is represented by an element $v \in S^{n-1}$. By Lemma \ref{doublingarg} it will be sufficient to embed $B_{[v]}(\delta) \subset S^{n-1} / \Gamma$. Define $\Gamma_\delta = \langle g \in \Gamma | d(gv,v) \leq 8\delta \rangle$, $K_\delta = \langle g \in K | d(gv,v) \leq 8\delta \rangle$ and $K'_\delta = \Gamma_\delta \cap K$. Decompose $v$ into it's components $v_i \in V_i$. There is an index $i$ such that $||v_i|| \geq \frac{1}{3\sqrt{n}}$. To reduce clutter, reindex so that $i = 1$, and define $\xi_1 = \frac{v_1}{||v_1||}$. But then $d(g\xi_1, \xi_1) \leq 24\sqrt{n}\delta$. Repeating the argument in the proof of Theorem \ref{flatmanifemb} we get either an element $w \in S^{n-1}$ fixed by $\Gamma_\delta$, or a Hopf-orbit $S_w^1 \subset V_1 \cap S^{n-1}$ through $w$ which is invariant under $\Gamma_\delta$. By using the normal exponential map, we get that $B_{[v]}(\delta) \subset S^{n-1}$ is bi-Lipschitz to a subset of either a $n-1$ dimensional locally flat orbivector bundle over a point (i.e. $w$), or a $n-2$ dimensional orbivector bundle over a one dimensional compact flat orbifold (i.e. $S_w^1 / \Gamma_\delta$). These can be embedded using the statements $V-(0,n-1)$ or $V-(1,n-2)$.  

\paragraph{Reduction $V-(n,0)$ to $V-(0, n)$ when $S=1$:}  This is very similar to the $V-(n,0)$ case in Theorem \ref{flatmanifemb}. We still have the estimate $\diam(M)<C(n)|l_1|$. Let $p \in \R^n$ be arbitrary. Set $\delta=l_1/(320c(n)L)$. Define the group $\Gamma_\delta=\langle g \in \Gamma | d(gp,p)<8\delta \rangle$. On the other hand, by Lemma \ref{loc} for any $[p] \in M$ we have $B_{[p]}(\delta) \subset M$ is isometric to $B_{[p]_\delta}(\delta)\subset \R^n / \Gamma_\delta(p)$. Thus, by Lemma \ref{doublingarg}  we can reduce the embedding problem to finding embeddings for $B_{[p]_\delta}(\delta)$. Define $\Lambda_0=\langle \lambda \in \Lambda | d(\lambda p,p)<8\delta \rangle = \{e\}$, and similarly to before we get

$$0 \to 0 \to \Gamma_\delta \to H_0 \to 0,$$
where $H_0$ and thus $\Gamma_\delta$ is finite. The group $\Gamma_\delta$ has a fixed point $o$ by Lemma \ref{fixedpoint}, and thus we can identify $\R^n / \Gamma_\delta$ by a quotient of isometries fixing $o$. Thus, $\R^n / \Gamma_{\delta}$ can be embedded by reduction to the $V-(0,n)$-case. 

\paragraph{Reduction $V-(n,0)$ to $V-(n - k, k)$ when $S>1$:}  The proof for manifolds translates almost verbatim.  Consider $\delta =\max\left(\frac{l_S}{200Lc(n)}, 2 l_{S-1}\right)$. Again $\diam(M)<C(n) l_S \leq 200Lc(n) C(n) \delta$ by Lemma \ref{nilpdiam}. Thus, it is sufficient by Lemma \ref{doublingarg} to embed each ball $B_{[p]}(\delta)$. Next fix $[p] \in M$, and its representative $p \in \R^n$. Repeating the arguments from the proof of Theorem \ref{flatmanifemb} we can show that $B_{[p]}(\delta)$ is isometric to an orbivector bundle over a $(n-k)$-dimensional flat orbifold $O$. The only difference is that the action is not free. 

\paragraph{Reduction $V-(n, d)$ to $V-(n - k, k + d), V - (a, b)$ with $k>0$ or $a+b<n+d$.} There is no major change to the proof in Theorem \ref{flatmanifemb} as we didn't use the fact that the action is free in any substantive way. The same decomposition into $T_i$'s is used, followed by a consideration of the action of $\Gamma_\delta$ on $S^{d-1}$, and finding a fixed point $w$ or a fixed Hopf-orbit $S_\xi^1$.

\QED

\textbf{Remark on bounds:} The methods are likely not optimal, so we only give a rough estimate of the size of the dimension and distortion in the previous construction. Similar bounds could be obtained for the $\epsilon(n)$-quasiflat spaces considered below. Let $D(n)$ be the worst case distortion for a flat orbifold of dimension $n$, and let $N(n)$ be the worst case dimension needed. In particular $D(n) \leq \max_{k=0, \dots, n} L(n-k,k)$, and $N(n) \leq \max_{k=0, \dots, n} N(n-k,k)$, where $L,N$ are defined in the statement of $V-(n,k)$ after Theorem \ref{orbifoldflat}. For dimension $n=1$, we can bound the distortion by constants $L(1,0)$, and $L(0,1)$, and they can be embedded in $3$ dimensions. 

Next, $V-(n,d)$ was reduced to $V-(n-k,d+k)$ and $V-(a,b)$ with $a+b<n+d$. Each such reduction is done by using Lemma \ref{doublingarg} at the scale $R/r \sim n!^{n+1} \sim O\left(e^{Cn^2 \ln(n)}\right)$, which increases distortion by a multiplicative factor of $O\left(e^{Cn^2 \ln(n)}\right)$. In order to reduce $V-(n,d)$ to lower dimensional spaces, this reduction might need to be repeated up to $n$ times. Thus, $D(n) = O\left(e^{Cn^3 \ln(n)}D(n-1)\right)$, which gives $D(n) =O\left(e^{Cn^4 \ln(n)}\right)$.  By a similar analysis, we get that the dimension of the embedding can be bounded by $N(n)=O\left(e^{Cn\ln(n)}\right)$. Our arguments are not optimal and probably some exponential factor could be removed, but the methods don't seem to yield better than super exponential distortion.

\subsection{Quasiflat orbifolds}

Below, the proof of Theorem \ref{thm:vector} will proceed by approximating the orbifold by a flat orbifold. This will involve quotienting away the commutator subgroup $[N,N]$ and using Lemma \ref{gromhaus}. Thus, we begin by studying lattices and their commutators.

\begin{lemma} \label{commut} Let $\Lambda<N$ be a lattice. Then $[\Lambda, \Lambda]$ is a lattice in $[N,N]$. 
\end{lemma}

\noindent \textbf{Proof:} $\Lambda'=[\Lambda,\Lambda]$ is a torsion free nilpotent Lie group and by Malcev it corresponds to a connected Lie group $N'$ such that $N'/\Lambda'$ is compact (see \cite{malcevclass} and \cite{chgrfu}). Further, by Malcev there is a unique homomorphism $f\co N' \to N$ such that $f_*(\Lambda')=\Lambda'$. We want to show that $f(N')=[N,N]$. First of all, since $\Lambda'<[N,N]$, we get that $f(N')<[N,N]$. 

The image $f(N')$ is normal in $N$, and $N/f(N')$ is a nilpotent Lie group. The normality follows since $\Lambda'$ is normal in $\Lambda$, and thus $\lambda f(N') \lambda^{-1} < N'$ for all $\lambda \in \Lambda$. From this we can conclude, using the nilpotency of $N'$ and the fact that $\Lambda$ is a lattice, that $N f(N') N^{-1} < N'$. This follows since every element can be represented as $n = \prod_i e^{t_i X_i}$. The derivative of $\text{Ad}_n(x)=nxn^{-1}$ at the identity is $\text{ad}_n(t_i)$, which is a matrix-polynomial in $t_i$, which has the tangent space of $f(N')$ as an invariant subspace whenever $t_i$ is integer. Thus, the same tangent space is left invariant by $\text{ad}_n(t_i)$ for any $t_i \in \R$.

The projection $\pi\co N / \Lambda' \to N / f(N')$ has compact fibers and thus we can show that the action of $\Lambda / [\Lambda, \Lambda]$ on $N / f(N')$ is properly discontinuous and co-compact. However, $\Lambda / [\Lambda, \Lambda]$ is abelian, and thus $N / f(N')$ is also abelian (by another application of Malcev's theorems). Therefore, $[N,N]<f(N')$, which completes the proof.

 \QED

\begin{lemma}\label{diameterest} Let $n$ be an integer, and $\epsilon(n)$ be sufficiently small\footnote{See the remark at the beginning of Section \ref{prelimin}.} and $l \leq \epsilon(n)$ fixed. There exists a dimension dependent constant $C(n)$ such that the following holds. If $N$ is a $n$-dimensional simply connected nilpotent Lie group with a left invariant metric $g$ and sectional curvature $|K| \leq 2$, and if $\Lambda<N$ is a lattice which is generated by elements of length $|l|_e  \leq l \leq \epsilon(n)$, then $\diam(N / \Lambda)\leq C(n) l$.
\end{lemma}

\noindent \textbf{Proof:} As in Lemma \ref{strat}, we choose a short basis $\gamma_i$ for the lattice $\Lambda$. By construction $|\gamma_i|_e \leq l$ for all $i=1, \dots, n$. The result then follows from Lemma \ref{nilpdiam}, with the same choice of $C(n)$. 

\QED

\begin{lemma}\label{approx} Let $S$ be an $\epsilon(n)$-quasiflat $n$-dimensional orbifold. Then there exist constants $0<\delta<1$ depending only on the dimension $n$ and a constant $A(n)$, such that the following two properties hold.

\begin{itemize}
\item For every $[p] \in S$ there exists a flat Riemannian orbivector bundle $V=V(S')$ over a $k$-dimensional $\epsilon(k)$-quasiflat orbifold $S'$ with $k < n$, with dimensions of the fibers $n-k$, and a $A(n)$-bi-Lipschitz map $f_{[p]} \co B_{[p]}(\delta \text{diam}(S)) \to V(S')$.
\item There exists a flat orbifold $S''$ such that $d_{GH}(S,S') < \delta \text{diam}(S)$.
\end{itemize}
\end{lemma}

\noindent \textbf{Proof:} 
Represent $S$ as a quotient $S=N / \Gamma$, where $\Gamma$ is the discrete group guaranteed by Definition \ref{quasiflatdef}. Further, take the finite index subgroup $\Lambda = N \cap \Gamma$ with $[\Gamma : \Lambda ] < c(n)$, which is guaranteed by Lemma \ref{finiteindex}.  Let $l=400c(n), L=2l^n,I_j, J_k, \gamma_j, \Lambda_k, L_k$ and $s$ be as in Lemma \ref{strat} and Definition \ref{shortbasis}. As in Lemma \ref{strat}, we need to assume $\epsilon(n)$ small enough. For this proof, $\epsilon(n) \leq \frac{\delta(n)}{C(n)10^3 c(n)C'(n) L}$ is enough. See Section \ref{prelimin} for the definitions of the various constants. We will first prove the existence of a near-by flat orbifold.

As in Lemma \ref{strat}, we can generate $\Lambda$ by elements $\alpha_i$ such that $|\alpha_i|_e \leq 4c(n)\diam(S) \leq 2c(n) \epsilon(n)$. As such, from nilpotency we see that  $[\Lambda, \Lambda]=\{[\mu,\nu] | \mu, \nu \in \Lambda \} < \Lambda$ is generated by elements $\beta_{i,j}=[\alpha_i, \alpha_j] \in N$. From \eqref{eq:commest} we get $|\beta_{i,j}|_e<\diam(S)/(10C(n) L^2) < \diam(S)/(8LC(n)c(n))$. 

Choose the largest index $a$ such that $|\gamma_a|_e<\diam(S)/(8LC(n)c(n))$, and let $a \in I_b$ for some $b$. By construction $l_b \leq L|\gamma_a|_e \leq \diam(S)/(8C(n)c(n))$. Also, $l_b < l_s$ because $\text{diam}(S)\leq C(n) l_s$ by Lemma \ref{nilpdiam}. Thus, $L_b \neq N$.

We have that $\beta_{i,j} \in \Lambda_b$, because if $\beta_{i,j} \not\in \Lambda_b$, then, by construction, $i_{b}+1 \not\in J_b$ and we have (by construction of the basis in Lemma \ref{strat})

$$|\gamma_{i_{b}+1}|_e \leq |\beta_{i,j}|_e < \text{diam}(S)/(8Lc(n)C(n)),$$
which would contradict the choice of $a$. Thus, we see that $[\Lambda, \Lambda] < L_b$ and that $l_b \leq \text{diam}(S)/(8C(n)c(n))$.  Since $L_b$ is invariant under conjugation by elements in $\Gamma$ and the holonomy of $\Gamma$, we see that the action of $\Gamma$ descends to an action of a discrete group $\Gamma'$ on $N / L_b=N'$, where $\Gamma'=\Gamma / \Gamma_b$ and $\Gamma_b$ is the isotropy group of $L_b$. Define $N' / \Gamma' = S'$. Note that $\Gamma'$ has a translational subgroup $\Lambda' = \Lambda / \Lambda_b$ which is abelian and acts co-compactly on $N'$. Therefore, $N'$ is in fact flat and $S'$ is a flat orbifold. We will next show that $d_{GH}(S',S)<\delta \text{diam}(S)$ for $\delta=1/(2c(n))$.

The orbits of the $L_b$-action on $N $ descend to orbits on $S$, and $S'$ can be identified by the orbit-space of the $L_b$-orbits on $S$. This gives a map between orbifolds $\pi\co S \to S'$, which lifts to a Riemannian submersion on their orbifold universal covers, which are $N$ and $N'$ respectively. We have $\diam(S') \leq \diam(S) \leq \epsilon(n)$. Choose a coset $[p]' \in S'$ represented by an element $p\in N$ with 
\begin{equation} \label{lengt} |p| \leq 4c(n)\epsilon(n).
\end{equation} Denote by $[p]=\pi^{-1}([p]')$ the coset in $S$.

Thus, it will be sufficient to prove that $\diam(\pi^{-1}([p]))<\delta \text{diam}(S)/2$. We see that $\pi^{-1}([p]') = O_{[p]}$, where $O_{[p]} \subset S$ is the orbit of $L_b$ passing through $[p] \in S$. The set $O_{[p]}$ is isometric to $L_b / \Gamma_p$ where $\Gamma_p = \{p^{-1} \gamma p | \gamma \in \Gamma \}$ with a translational subgroup $\Lambda_p = p^{-1}\Lambda_b' p$. This translational subgroup has a basis given by $p^{-1} \gamma_i p$ for $i=1, \dots, i_{b+1}-1$. We have, by \eqref{eq:commest} and \eqref{lengt} that 
$$|p^{-1} \gamma_i p|_e \leq 2|\gamma_i|_e \leq \text{diam}(S)/(4C(n)c(n)).$$ By Lemma \ref{diameterest} we have $\diam(L_b / \Gamma_p) \leq \text{diam}(S)/(4c(n))$. Thus, $d_{GH}(S',S) \leq \delta \text{diam}(S)$ follows by our choice $\delta=1/(2c(n))$ and using Lemma \ref{closeness}. 

Next, we prove the first part, that at the scale of $\delta'=\delta \text{diam}(S)$ the space has a description in terms of a simpler vector bundle. Define $\Gamma_{\delta'}(p)=\langle \gamma \in \Gamma | d(\gamma p, p)\leq 8\delta'\rangle$. By Lemma \ref{small} we know that $B_{[p]}(\delta') \subset S$ is isometric to $B_{[p]_{\delta'}}(\delta') \subset N / \Gamma_{\delta'}(p)$. Further by Lemma \ref{loc}, and observing that $\Lambda_{\delta'} = \langle \gamma \in \Lambda | d(\gamma p, p)\leq 8\delta'\rangle = \Lambda_b$ we have a finite group $H_0$ and a short exact sequence

$$0 \to \Lambda_b \to \Gamma_{\delta'} \to H_0 \to 0.$$

The action of $H_0$ descends to an action on $N / L_b$ by $h[x]_{L_b}=[hx]_{L_b}$. Further, by Lemma \ref{loc}, for every element $h \in H_0$ we can choose a representative $\gamma_h \in \Gamma_{\delta'}$  with $d(\gamma_h p, p) \leq 8(c(n)+1)\delta'$. In particular $d(h[p]_{L_b}, [p]_{L_b}) \leq 8(c(n)+1)\delta'$. Thus, by Lemma \ref{fixedpoint} used on $N/L_b$ (which is flat) there is an orbit $[q]_{L_b}$ which is fixed by the action of $H_0$ and $d([q]_{L_b}, p)\leq 8(c(n)+1)\delta'$. Further, the orbit $[q]_{L_b} \subset N$ is invariant under the action of $\Gamma_{\delta'}$. 

Consider the normal bundle $\nu([q]_{L_b})$ of $[q]_{L_b}$ in $N$. The bundle has a trivialization and can be identified with $\R^t \times L_b$ for some $t$. This trivialization can be obtained similarly to Lemma \ref{metricnilp}. By Lemma \ref{metricnilp} there is a map $f \co N_{2\delta(n)}([q]_{L_b})\to \R^t \times L_b$ from the $2\delta(n)$-tubular neighborhood $N_{2\delta(n)}([q]_{L_b})$ of $[q]_{L_b}$. This map is also $A(n)$-bi-Lipschitz on the ball $B_q(2\delta(n))$. Further, the group $\Gamma_{\delta'}$ acts on $N_{2\delta(n)}([q]_{L_b})$ and $\R^t \times L_b$ by isometries, and the action commutes with $f$. We can thus show that restricted to the ball $B=B_{[q]_{\delta'}}(\delta(n)) \subset N / \Gamma_{\delta'}$ the induced map $\overline{f}\co N_{2\delta(n)}([q]_{L_b}) / \Gamma_{\delta'} \to \R^t \times L_b / \Gamma_{\delta'}$ is also $A(n)$-bi-Lipschitz. 

Denote the orbifold by $L_b / \Gamma_\delta = S'$, and its dimension by $k<n$. The orbifold $S'$ is lower dimensional because $L_b < N$ strictly. The space $\R^t \times L_b / \Gamma_\delta = V(S')$ is a locally flat Riemannian orbivector bundle over the $k$-dimensional $\epsilon(k)-$quasiflat orbifold $S'$\footnote{$S'$ will have the same bound on its curvature and diameter as $S$. We use here the assumption that $\epsilon(k)\geq\epsilon(n)$ for $k \leq n$.}. Also, $B_{[p]_{\delta'}}(\delta') \subset  B_{[q]_{\delta'}}(\delta(n))$, so further restricting we get a $A(n)$-bi-Lipschitz map $\overline{f}\co B_{[p]_{\delta'}}(\delta')  \to V(S')$, which is our desired map.

\QED

A similar statement is needed for orbivector bundles.

\begin{lemma}\label{approxvect} Let $\epsilon(n)$ be small enough. For every $d$ and $n$ there are constants $A(n,d)$ such that the following holds. Let $V=V(S)$ be a locally flat Riemannian $d$-dimensional orbivector bundle over an $\epsilon(n)$-quasiflat $n$-dimensional orbifold $S$. Then there exist constants $0<\delta,\delta'<1$, which depend only on the dimension $n$ and $d$, such that the following two properties hold.

\begin{itemize}
\item For every $p \in V(S)$ there exists a flat Riemannian orbivector $V=V(S')$ bundle over a $k$-dimensional $\epsilon(k)$-quasiflat orbifold $S'$ with $k < n$ and a $A(n,d)$ bi-Lipschitz map $f_p \co B_{p}(\delta\diam(S)) \to V(S')$.
\item There exists a locally flat Riemannian orbivector bundle $V''=V(S'')$ over a flat basis $S''$ such that $d_{GH}(V,V'') < \delta'\text{diam}(S)$.
\end{itemize}
\end{lemma}

\noindent \textbf{Proof:} For this proof $\epsilon(n) \leq \frac{\delta(n)}{C(n)10^3 c(n)C'(n) L}$ is enough (see Definition \ref{shortbasis} for the definition of $L$). We can represent the spaces as $S=N / \Gamma$ and $V= N \times \R^d / \Gamma$, where $\Gamma$ acts by rotations on $\R^d$ and co-compactly on $N$ so that $S$ is $\epsilon(n)$-quasiflat. By Lemma \ref{finiteindex}, there is a lattice $\Lambda = \Gamma \cap N$  of $N$ with $[\Gamma : \Lambda] \leq c(n)$. Further, there is an induced action $h\co \Lambda \to O(d)$. By Lemma \ref{jordan} we have an index $c(d)$ subgroup $\Lambda' < \Lambda$ with $[\Lambda : \Lambda'] \leq c(d)$ and $h(\Lambda')$ abelian. In particular $[\Lambda', \Lambda']=\Lambda'_c < \Lambda'$.

The group $\Lambda'_c$ acts trivially on $\R^d$ and by translations on $N$. Consider the left action of $N_c$ on $N\times \R^d$, which acts trivially on $\R^d$. The group $\Gamma$ maps orbits of $N_c$ to similar orbits, and thus $\Gamma_c=\Gamma / \Lambda'_c$ is a discrete group which acts properly discontinuously on $N  / N_c \times \R^d$. Further, the orbits of $N_c$ descend to $V$. Denote by $V_c=(N / N_c \times \R^d) / \Gamma_c$ the group quotient, and $\pi\co V \to V_c$ the induced projection map. The group action of $\Gamma_c$ is co-compact on $N  / N_c$, and by rotations on $\R^d$. Further, $N  / N_c$ is abelian, and thus isomorphic and isometric to $\R^l$ for some $l$. In other words, up to an isometry $V_c = \R^l \times \R^d / \Gamma_c$, which is a $d$-dimensional locally flat orbivector bundle over the flat orbifold $\R^l  / \Gamma_c$. 

Next, we give a diameter bound needed to apply Lemma \ref{closeness}. Take any $q=[(p,v)] \in V$ and the $N_c$ orbit $O_q$ passing through that point. We can assume, as in Lemma \ref{approx}, that $|p|_e \leq 2c(n)\epsilon(n)$. Since $\Lambda'$ is a lattice of $N$ by Lemma \ref{commut}, also $\Lambda_c'$ is a lattice in $[N,N]=N_c$.  Also $O_q = N_c / p^{-1}\Lambda_cp$. The compact manifold $O_q$ is a subset of $N / p^{-1}\Lambda'p$, and $N_c / p^{-1}[\Lambda, \Lambda]p$ is a compact manifold with $\diam(N_c / [\Lambda, \Lambda]) \leq 20C(n)\diam(S)c(n)$ by computations from Lemma \ref{approx}. Further, $N_c / p^{-1}[\Lambda, \Lambda]p$ covers $O_q$ with index at most $c(d)$. Thus,

\begin{equation}\label{diamestimate}
\diam(O_q) \leq 50C(n)\diam(S)c(n)c(d).
\end{equation}

 Note that $O_q $ can be covered by $N \times \{(p,v )\} / \Lambda_c$, where $(p,v)$ is a representative of $q$ in $N \times \R^d$. Also, the action of $\Lambda_c$ is trivial on $\R^d$. Define $\delta'=100C(n)c(n)c(d)$ and $\Delta = 100C(n)c(n)c(d) \diam(S)$. As before, we can show that $\pi$ is a $\Delta$-Gromov-Hausdorff approximation. In particular $d_{GH}(V_c, V) \leq \Delta$. This shows the second part.

We will next show the first statement by a similar argument as in Lemma \ref{approx}. Fix an arbitrary $p=(q,v) \in N \times \R^d$. Let $l=400c(n), L=2l^n,I_j, J_k, \gamma_j, \Lambda_k, L_k$ and $s$ be as in Lemma \ref{strat} and Definition \ref{shortbasis}. As in Lemma \ref{strat}, we need to assume $\epsilon(n)$ small enough.

As in Lemma \ref{strat}, choose the largest index $a$ such that $|\gamma_a|_e<\diam(S)/(8Lc(n)C(n))$, and let $a \in I_b$ for some $b$. By construction $l_b \leq L|\gamma_a|_e \leq \diam(S)/(8c(n)C(n))$ and $l_b < l_s$.

Let $\delta = 1/(8c(n))$ and redefine $\Delta=\delta \diam(S)$. Define $\Gamma_{\Delta}(q)=\{\gamma \in \Gamma | d(\gamma q, q)\leq 8\Delta\}$ and  $\Gamma_{\Delta}(p)=\{\gamma \in \Gamma | d(\gamma p, p)\leq 8\Delta\}$. Since $\Gamma_{\Delta}(p) < \Gamma_{\Delta}(q)$, by Lemma \ref{small} we know that $B_{[p]}(\Delta) \subset V$ is isometric to $B_{[p]_{\Delta}}(\Delta) \subset N \times \R^d / \Gamma_{\Delta}(q)$. Note, we first show that it is isometric to a ball in $N / \Gamma_{\Delta}(p)$, but by the inclusion $\Gamma_{\Delta}(p)<\Gamma_{\Delta}(q)$ it is easy to see that this induces the desired isometry. Further, by Lemma \ref{loc} we have a finite group $H_0$ and a short exact sequence

$$0 \to \Lambda_b \to \Gamma_{\Delta}(q) \to H_0 \to 0.$$

Here, $\Lambda_b<N$ is defined similarly to Lemma \ref{approx}. The group $\Lambda_b$ is a lattice in the connected Lie subgroup $L_b < N$ associated to it (see Lemma \ref{lem:basis}). By the same argument as in Lemma \ref{approx}, we can find a $L_b$-orbit $[x]_{L_b}$ passing through $x\in N$ with $d(x,q)<8(c(n)+1)\delta'$ in $N$ which is invariant under the action of $\Gamma$. By Lemma \ref{metricnilp}, we find for the $\delta(n)$-tubular neighborhood $N_{\delta(n)}([x]_{L_b})$ of $[x]_{L_b}$ in $N$ a map $f\co N_{\delta(n)}([x]_{L_b}) \times \R^d \to L_b \times \R^t \times \R^d$, which is identity on the $\R^d$ factor. Further, $f$ is $A(n)$-bi-Lipschitz when restricted to $B_x(2\delta(n)) \times \R^d$. 

The action of the group $\Gamma_{\Delta}(q)$ induces canonically an action on $L_b \times \R^t \times \R^d$. Quotienting by this action induces a map $\overline{f}\co N_{\delta(n)}(([x]_{L_b})) \times \R^d / \Gamma_{\Delta}(p) \to L_b \times \R^t \times \R^d / \Gamma_{\Delta}(p)$. Let $B(\delta(n))=N_{\delta(n)}(([x]_{L_b})) \times \R^d / \Gamma_{\Delta}(p)$. The space $V''= L_b \times \R^t \times \R^d / \Gamma_{\Delta}(p)$ is a locally flat Riemannian orbivector bundle over $L_b / \Gamma_{\Delta}(p)$, which is an $\epsilon(n-l)$-quasiflat $n-l$-dimensional orbifold for some $n\geq l>0$. Note that $\epsilon(n-l) \geq \epsilon(n)$ by assumption.

Finally, consider the ball $B_{[p]_{\Delta}}(3\Delta) \subset B(\delta(n))$. This inclusion can be seen by considering lifted balls and noting that $B_p(2\Delta) \subset A(\delta(n))$. Note that we have chosen $\epsilon(n)$ so small that $8(c(n)+2)\Delta<\frac{\delta(n)}{10}$. The restricted map $\overline{f}|_{B_{[p]_{\Delta}}(\Delta)}$ is also $A(n)$-bi-Lipschitz.

\QED

Finally, we can prove the embedding theorem for quasiflat manifolds and vector bundles.

\noindent \textbf{Proof of Theorem \ref{thm:vector}:} The proof will proceed by a reverse induction argument initiated by the flat case of the theorem. We will define the induction statements $E(d,n)$ as follows. 

\begin{itemize}
\item $E(0,n)$: There exist $N(0,n)$ and $L(0,n)$ such that every $\epsilon(n)$-quasiflat $n$-dimensional orbifold $S$ admits a bi-Lipschitz map $f\co S \to \R^{N(0,n)}$ with distortion $L(0,n)$.
\item $E(d,n)$ for $n \geq d\geq 1$: There exist $N(d,n)$ and $L(d,n)$ such that every $d$-dimensional locally flat orbivector bundle $V=V(S)$ over an $\epsilon(n-d)$-quasiflat $n-d$-dimensional orbifold $S$ admits a bi-Lipschitz map $f \co V \to \R^{N(d,n)}$ with distortion $L(d,n)$.
\end{itemize}

The statement will be proved by an induction on $n$. The case $n=1$ is completed since the cases included in $E(0,1)$ and $E(1,1)$ are flat orbifolds and can be embedded using Theorem \ref{orbifoldflat}. Next, fix $n>1$ and assume we have proved $E(t,s)$ for all $1 \leq s \leq n$ and $0 \leq t \leq s$. 

We proceed to show $E(d,n)$ by reverse induction on $d$. The  base cases are $d=n,n-1$. The statement $E(n,n)$ corresponds to a flat orbifold and is covered in Theorem \ref{orbifoldflat}. Further, the case $E(n-1,n)$ corresponds to a locally flat orbivector bundle over a flat base, so is also covered in Theorem \ref{orbifoldflat}. 

Next, fix $0 \leq d \leq n-2$. Assume the statement $E(k,n)$ has been proved for $k>d$. There are two slightly different cases $d=0$ and $d\geq 1$. Assume the first. In this case by Lemma \ref{approx} there is a $\delta$ such that $B_p(\delta \diam(S))$ is bi-Lipschitz to a case covered by $E(s,n)$ for some $s>0$. Thus, $B_p(\delta \diam(S))$ admits a bi-Lipschitz embedding by induction. Also, $d_{GH}(S,S')\leq \delta \diam(S)$ where $S'$ is a flat orbifold of lower dimension. Also, $S'$ admits a bi-Lischitz embedding by Theorems \ref{orbifoldflat}. The statement then follows from Lemma \ref{gromhaus}. The case $d \geq 1$ is similar, but one uses Lemma \ref{approxvect} instead and concludes that $B_p(\delta \diam(S))$ is bi-Lipschitz to a case covered by $E(d+s,n)$ for some $s>0$, and $d_{GH}(V,V')\leq \delta' \diam(S)$ for a locally flat orbivector bundle $V'$ over a flat base. The constants $\delta, \delta'$ here depend only on the dimension.

We remark, that to apply Lemma \ref{gromhaus} we technically need to show that every $\epsilon(n)$-quasiflat space is doubling. This could be shown by combining techniques from \cite[IV.5.8]{varopoulos2008analysis} (and \cite{nagelstein}). Also, Lemmas \ref{approxvect} and \ref{approx} combined with induction actually can be used to prove doubling directly. However, in order to apply the proof of Lemma \ref{gromhaus} we only need the doubling estimate on all balls $B_p(R) \subset V$ and $B_p(R) \subset R$, where $R>0$ has a fixed upper bound depending on $\epsilon(n), \delta, \delta'$ and the constants $N(s,t), L(s,t)$ for cases $E(s,t)$ assumed true by induction. This doubling estimate follows by the sectional curvature estimate $|K| \leq 2.$

\QED

\subsection{Proofs of main theorems}

With the previous special cases at hand we can prove the main theorems.

\textbf{Proof of Theorem \ref{thm:bilip}:} Fix $\epsilon=\epsilon(n)$ and $\eta=1/4$. By Fukaya's fibration Theorem  \ref{thm:fukayamanifold}, every point $p \in M^n$ has a definite sized ball
$B_p(\delta)$, where $\delta=5\rho(n,\eta,\epsilon)$, equipped with a metric $g'$ such that $(B_p(\delta), g')$ is isometric to a subset of a compact $\epsilon(n)$-quasiflat manifold $M'$ or a locally flat vector bundle $V$ over such a base. Since $||g - g'
|| < 1/4$, then the ball $B_p(\delta/2)$ is bi-Lipschitz to a subset of $V$, or $M'$, with distortion at most $\frac{1+1/4}{1-1/4}$.
Since the spaces $M'$ and $V$ can be embedded by Theorem \ref{thm:vector}, we can embed $M^n$ using a doubling Lemma \ref{doublingarg} with $R=\diam(A)$ and $r=\delta$. Note, by our choice of $\epsilon$ and $\eta$, the constant $\delta$ depends only on the dimension $n$.

\QED

\textbf{Proof of Theorems \ref{thm:orbifoldmain} and \ref{thm:flatellipt}:}
The only difference to the proof of Theorem \ref{thm:bilip} is using Fukaya's fibration theorem for orbifolds \ref{thm:fukayaorbifold} with $\epsilon=\epsilon(n)$ and $\eta=1/4$.

\QED

\section{Appendix} \label{appendix}

\subsection{Collapsing theory of orbifolds}

Much of the terminology here is presented in \cite{Kleiner2014,thurston2002geometry}. Recall, a \textit{topological orbifold} is a second countable Hausdorff space with an atlas of charts $U_i$, which are open sets satisfying the following.

\begin{itemize}
\item Each $U_i$ has a covering by $\hat{U}_i$ with continuous open embeddings $\phi_i\co \hat{U}_i \to \R^n$ and a finite group $G_i$ acting continuously on $\hat{U}_i$ such that $U_i$ is homeomorphic to $\hat{U}_i / G_i$. Denote the projection by $\pi_i\co \hat{U}_i \to U_i$.
\item There are transition homeomorphisms $\phi_{ij}\co \pi_i^{-1}(U_i \cap U_j) \to \pi_j^{-1}(U_i \cap U_j)$, such that $\pi_j \circ \phi_{ij} = \pi_i$.
\end{itemize}

A chart is denoted simply by $U_i$, without explicating the covering $\hat{U}_i$ or other aspects. An \textit{orbifold} is defined as one for which $\phi_j \circ \phi_{ij}\circ \phi_{i}^{-1}$ are smooth maps. Further, a \textit{Riemannian orbifold} is an orbifold with a metric tensor $g_i$ on each $U_i$ for which $d\phi_{ij}^* (g_j) =g_i$. Any other notion like geodesics, sectional curvature or Ricci-curvature, is defined by considering the lifted metric on $\hat{U}_i$. One may also define the notion of orbifold maps as those that lift to the covering charts $\hat{U}_i$ and are compatible with the projections. A smooth orbifold map is also naturally defined by considering charts. Such a map is non-singular if in charts it can be expressed as a local diffeomorphism. A non-singular orbifold map which is surjective is called an orbifold covering map.

We also need a notion of an orbivector bundle and an orbiframe bundle. An orbivector bundle is a triple $(V,O,\pi)$, where $V,O$ are smooth orbifolds and $\pi\co V \to O$ is a smooth map with the following properties.

\begin{itemize}
\item There is an atlas of $O$ consisting of charts $U_i$ such that for all of the charts $(U_i,\hat{U}_i,G_i)$ of $O$ there is a chart $V_i = \pi_i^{-1}(U_i)$
\item There are covering charts $\hat{V}_i = \hat{U}_i \times \R^k$ and a $G_i$-action, such that $V_i = \hat{V}_i / G_i$, and maps $\pi^V_i\co \hat{V}_i \to V_i$. 
\item The projections $\hat{\pi}_i\co \hat{V}_i \to \hat{U}_i$ satisfy $\hat{\pi}_i = \pi \circ \pi^V_i$. 
\item The action of $G_i$ on $\hat{U}_i \times \R^k$ splits as $g(u,x) = (gu,g_u x)$, where $G_i$ acts linearly on $\R^k$ ( the action can depend on the base point $u$), and the action on the first component is the original action.
\item The actions of $G_i$ and $G_j$ on overlapping charts commute.
\end{itemize}

If there is a bundle metric on each $\hat{V}_i$, such that $G_i$ acts by isometries on the fibers $\R^k$ and such that it is invariant under transition maps, then we call the orbivector bundle a Riemannian orbivector bundle. A connection on such bundle is defined as an affine connection on each $\hat{V}_i$ which is invariant under $G_i$ and the compatible via the transition maps. Such a connection is flat if it is flat on each chart. Such a connection together with a Riemannian structure on the base induces a metric on $V$ making it into a Riemannian orbifold. This allows us to define a central notion in this chapter.

\begin{definition}
$(V,O,\pi)$ is called a locally flat orbivector bundle, if it is a Riemannian orbivector bundle over a Riemannian orbifold with a flat connection.
\end{definition}

To every smooth orbifold we can associate a tangent orbivector bundle $TO$, where the group actions on the fibers are given by differentials. The fiber of a point $\pi^{-1}(p) \subset TO$ is denoted $C_p(O)$. For Riemannian orbifolds, this is easily seen to correspond with the metric tangent cone at $p$. In that case $C_p(O) = \R^n / G$ for some group $G$ acting by isometries on $\R^n.$ One can also define orbifold principal bundles by considering a fiber modeled on a Lie group $G$ and assuming transition maps are given by left-multiplication. This allows us to define for Riemannian orbifolds an orbiframe bundle $FO$. The orbiframe bundle $FO$ is always a manifold.

In order to do collapsing theory, we need to define a local  orbifold cover. A connected Riemannian orbifold has a distance function and we can define a complete Riemannian orbifold as one which is complete with respect to its distance function. On such orbifolds we may define an exponential map $\text{exp}\co TO \to O$ by considering geodesics and their lifts. Fix now a point $p \in O$ and take the exponential map at that point $\text{exp}_p \co C_p(O) \to O$, which is surjective. As $C_p(O) = \R^n / G$, we can lift the exponential map to a map $\widetilde{\text{exp}}_p\co \R^n \to O$. Next we state a crucial property of such an exponential map.

\begin{lemma}\label{nonsing} If a complete Riemannian orbifold has sectional curvature bounded by $|K|\leq 1$, then $\widetilde{\text{exp}}_p$ is a non-singular orbifold map on the ball $B_0(\pi) \subset \R^n$. 
\end{lemma}

The map $\text{exp}_p \co C_p(O) \to O$ will define locally a homeomorphism in the neighborhood of $p \in C_p(O)$. We call the largest $\delta$ such that $\text{exp}_p|_{B(p,\delta)}$ is injective on $B(p,\delta )\subset C_p(O)$ the \emph{injectivity radius} at $p$, and denote it $\text{inj}_p(O)$. If $O' \subset O$ is a sub-orbifold, then we can also define a normal bundle of it and a normal injectivity radius by the largest radius such that the normal exponential map is injective. These terms are used in the proof of Theorem \ref{thm:fukayaorbifold} below.

Lemma \ref{nonsing} is alluded to in the appendix of \cite{ding2011restriction} and can be proved by using local co-ordinates and standard Jacobi-field estimates. Using this terminology, we can now state and outline the proof of Fukaya's fibration theorem for collapsed orbifolds

\begin{theorem} Let $(O,g)$ be a complete Riemannian orbifold of dimension $n$ and sectional curvature $|K|\leq 1$, and $\epsilon>0$ be an arbitrary constant. For every $\eta > 0$ there exists a universal $\rho(n,\eta,\epsilon)>0$ such that for  any point $p \in O$ there exists a metric $g'$ on the ball $B_p(10\rho(n,\eta,\epsilon))$ and a complete Riemannian orbifold $O'$ with the following properties.

\begin{itemize} 
\item $$||g-g'||_{g}<\eta$$ 
\item  $(B_p(5\rho(n,\eta,\epsilon)),g')$ is isometric to a subset of $O'$.
\item $O'$ is either a $\epsilon$-quasiflat orbifold or a locally flat Riemannian orbivector bundle over a $\epsilon$-quasiflat orbifold $S$.
\end{itemize}
\end{theorem}

\textbf{Proof:} Fix $\epsilon,\eta>0$. Without loss of generality, we can assume $\epsilon$ is very small. If the result were not to hold, then it would have to fail for a sequence of pointed orbifolds $(O_i,p_i)$. Thus, let $O_i$ be a sequence of examples such that there doesn't exist any $\delta>0$ such that for any $i$ the ball $B_{p_i}(\delta)$ has a metric $g''_\eta$ with the following two properties. On the one hand $||g''_\eta-g||_g < \eta$, and $(B_{p_i}(\delta),g''_{\eta})$ is isometric to a subset of an $\epsilon$-quasiflat orbifold, or an orbivector bundle over such a base. We will see a contradiction by establishing the existence of such a $\delta$.  

One should first regularize the orbifold slightly as to attain uniform bounds on all co-variant derivatives. This is done in \cite[Second part, Section 5]{Kleiner2014}. The action of the local group on the frame bundle is free and thus Fukaya's proofs in \cite{fukaya1988,fukaya1989} (see also \cite{chgrfu}) show that a subsequence of the orbifold frame bundles $FO_i$ converges to a smooth manifold $Y$ with bounds on all covariant derivatives. Assume without loss of generality that $FO_i$ is already this sequence.

Using standard collapsing theory from Cheeger, Gromov and Fukaya \cite{chgrfu}, we obtain an equivariant fibration of $FO_i$ by nilpotent fibers, and a metric $g_\eta$ which is invariant under the local nilpotent action and the local $O(n)$ action. It is possible for those orbits to be  finite. Further, we have bounds on their second fundamental forms and normal injectivity radii of the nilpotent fibers. (In the case of a point, the normal injectivity radius is just the injectivity radius.) By a remark in \cite{rong1996fundamental}, we can also control the sectional curvatures of $g_\eta$. By equivariance, this structure descends onto a fibration of $O_i$ by orbits of nilpotent Lie groups. (Again, the fibers may be points, in which case the normal injectivity radius bound means the injectivity radius bound). Corresponding to this fibration we have, after possibly passing to another subsequence and re-indexing, an invariant metric $g'_\eta$, which satisfies $||g-g'_{\eta}||_g<\eta'(i)$ for some $\eta'(i)$ going to zero with $i$. Call the fiber passing through $p$ by $\mathcal{O}_p$. 

By the argument in \cite[Appendix 1]{chgrfu} we can find (for large enough $i$) a constant $\delta$ and a fiber $\mathcal{O}_{q_i} \subset O_i$ with $d(q_i,p_i) < \delta$, such that the orbit $\mathcal{O}_{q_i}$ has normal injectivity radius at least $20\delta$. One can also bound the second fundamental form of $\mathcal{O}_{q_i}$. The normal bundle has a natural metric inherited from $g'_\eta$, and a locally flat connection inherited from the nilpotent action on the fiber $\mathcal{O}_{q_i}$. By the second fundamental form bounds, for $i$ sufficiently large, the restricted metric $g'_\eta$ makes $\mathcal{O}_{q_i}$ a $\epsilon$-quasiflat orbifold. Further, the metric induced by the locally flat connection makes the normal bundle to $\mathcal{O}_{q_i}$ a locally flat Riemannian orbivector bundle over an $\epsilon$-quasiflat base. Denote by $g''_\eta$ the metric tensor which is induced on this normal bundle by the locally flat connection.

Denote by $N_{10\delta}(\mathcal{O}_{q_i})$ the vectors in the normal bundle of the fiber $\mathcal{O}_{q_i}$ with length at most $10\delta$. One can show, for $i$ sufficiently large, that the normal-exponential map $f\co N_{10\delta}(\mathcal{O}_{q_i}) \to O_i$ induces a map onto the $10\delta$-tubular neighborhood (with respect to $g'_\eta$ of $O_i$) with $||f^*g - g''_\eta||_g < \eta''(i)$ for some $\eta''(i)$ tending to zero with $i$. Further, the image of $f$ contains $B_p(\delta)$ and we can consider the restricted map $f^{-1} \co B_p(\delta) \to N_{10\delta}(\mathcal{O}_{q_i})$. This gives a metric $f^{-1*}(g''_\eta)=g'$ on $B_p(\delta)$ with the desired properties. The model space $O'$ is the normal bundle of $\mathcal{O}_{q_i}$ with the metric $g''_\eta$ and the isometry of $B_{p}(\delta)$ is given by $f$. However, we assumed that no such structure could exist, and thus the statement follows by contradiction.

\QED

\subsection{Diameter estimate}

\begin{lemma}\label{nilpdiam} Let $M = N / \Gamma$ be a $n$-dimensional $\epsilon(n)$-quasiflat orbifold and $\Lambda=N \cap \Gamma < \Gamma$ the translational normal subgroup of $\Gamma$. Further, let $\gamma_i$ be a short basis for $M$ as in Lemma \ref{strat}. Then we have $ \diam(M) \leq C(n) \max_{i=1, \dots, n}|\gamma_i|$.
\end{lemma}

Here, we use the notation $|\gamma|=|\gamma|_e=d(\gamma, e)$, where $e$ is the identity element of $N$, and $d$ is the Riemannian distance function on $N$.

\textbf{Proof:} Denote $D=\max_{i=1, \dots, n}|\gamma_i|$. To prove the lemma we can assume $\Gamma=\Lambda$, because taking a finite cover will only increase the diameter. We will show a more general statement. The statement we prove is that, if a cocompact lattice $\Lambda < N$ admits a triangular basis $\gamma_1, \dots, \gamma_n$  in the sense of definition \ref{def:basis}, then $ \diam(M) \leq C(n) \max_{i=1, \dots, n}|\gamma_i|$. The statement is obvious for $n=1$, and we will proceed to assume that $n>1$ and that the statement has been shown for $n-1$-dimensional simply connected nilpotent Lie groups. 

Assume that $\Lambda \subset N$ is a discrete co-compact lattice of a $n$-dimensional nilpotent Lie group with a triangular basis $\gamma_i$ and set $\gamma_i = e^{X_i}$ with $X_i \in T_{e}N$ a basis and $[X_i, X_j] \in \mathrm{span}(X_{1}, \dots, X_{i-1})$ for $i<j$. Since $\gamma_i$ is assumed to be a basis, the vectors $X_i$ form a basis for $T_{e}N$. (Recall Definition \ref{def:basis}). In particular, we can introduce coordinates $(t_1, \dots, t_n) \to \prod_{i=1}^n e^{t_{n-i}X_{n-i}}$, such that the product can be expressed as 

 $$\prod_{i=1}^n e^{t_{n-i}X_{n-i}} \times \prod_{i=1}^n e^{s_{n-i}X_{n-i}} = \prod_{i=1}^n e^{(s_{n-i}+t_{n-i}+p_{n-i}(s,t)) X_i},$$ 
where $p_i(s,t)$ depends only on $s_j,t_j$ for $j < i$. We will show that $\diam(N) \leq C(n) D$. 

Since $\gamma_1$ is in the center of $N$, we have that $X_1$ commutes with all left-invariant fields. Further, by the Koszul formula $t\to e^{tX_1}$ is seen to be a geodesic in $N$. Thus, $d(e^{s X_1},e) = s d(e^{X_1},e)=s|\gamma_1|$. Consider the Lie group $N'=N / L$, where $L = \{e^{tX_1}\}$, and a discrete subgroup $\Gamma' = \langle \gamma_i', i=1,\dots, {n-1}\rangle = \Gamma / \langle \gamma_1^n, n \in \Z \rangle $. Here, the generators correspond to equivalence classes of generators of $\Gamma$ in $\Gamma / \langle \gamma_1^n, n \in \Z \rangle $: $\gamma_i' = [\gamma_{i+1}]$, and if we denote $e^{X_i'}=\gamma'_i$, we have co-ordinates $(t_1, \dots, t_{n-1}) \to \prod_{i=1}^{n-1}e^{t_{n-i}X'_{n-i}}$. Here the exponential map is in $N'$. Take as the metric on $N'$ the quotient metric, which is induced by the restricted Riemannian metric on a vector space perpendicular to $\mathrm{span}(X_{1})$. With respect to this metric $|\gamma_i'| \leq |\gamma_{i+1}| \leq |\gamma_n| \leq D$ for all $1 \leq i \leq n-1$. The norm $|\cdot|$ on $N'$ is defined as $|n'|=d(n'e',e')$, where $n' \in N'$ and $e'$ is the identity element of $N'$. 

Define $M' = N' / \Gamma'$. The group $L$ acts naturally on M by isometries and $M/L$ is isometric to $M'$. Further, the $L$-orbits in $M$ have length $|\gamma_1|$, so by Lemma \ref{closeness} $d_{GH}(M,M')\leq 3|\gamma_1|$. By the inductive hypothesis, $\diam(M') \leq C(n-1)D$. Thus,

\begin{eqnarray*}
\diam(M) &\leq& \diam(M') + 2d_{GH}(M,M') \leq C(n-1)D+ 6D \\
                           &\leq& \left(6+C(n-1)\right)D.
\end{eqnarray*}

In particular, we can choose $C(n)=6n$.

\QED

\begin{lemma}\label{metricnilp} Let $N^n$ be a simply connected nilpotent Lie group and $M^k \triangleleft N^n$ a connected normal subgroup thereof. Assume that $N$ is equipped with a left-invariant metric which satisfies the sectional curvature bound $|K| \leq 2$. Then for some universal $\delta(n)> 0$ and the $4\delta(n)$-tubular neighborhood $N_{4\delta(n)}(M^k)$ of $M^k$ in $N^n$ there is a map $F\co N_{\delta(n)}(M^k) \to M^k \times \R^{n-k}$ whose restriction onto a ball $F\co B_e(2\delta(n)) \to M^k \times \R^{n-k}$ is $A(n)$-bi-Lipschitz. Also, $F$ is an isometry when restricted to $F^{-1}(M^k \times \{0\})$. Furthermore, for any affine isometry $I$ of $N^n$ which preserves $M^k$, also $F \circ I \circ F^{-1}$ is a restriced affine isometry of $M^k \times \R^{n-k}$.
\end{lemma}

\textbf{Proof:} Introduce an orthonormal triangular basis $X_1, \dots, X_k$ for $M^k$, and extend it to an orthonormal basis for $N^n$ with vectors $X_{k+1}, \dots, X_n$. Consider the map $G\co  M^k \times \R^{n-k} \to N $, which is defined by 

$$G(m,x)= me^{\sum_{i=k+1}^n x_i X_i}.$$

This is smooth. The exponential map used here is the Lie group exponential. Clearly, any affine isometry $I$ preserving $M^k$ will conjugate to be an isometry of $M^k \times \R^{n-k}$ by considering its action on $M^k$ and its normal bundle. Apply the Campbell-Hausdorff  formulas to compute the quantities $\partial_{x_i} G$.

$$\partial_{x_i} G(m,x) = X_i + \sum_{k=1}^n \sum_{i_1, \dots, i_k} M^k_{i_1,\dots, i_k}x_{i_1} \cdot x_{i_2} \cdots x_{i_k}[X_{i_1},[X_{i_2},[ \dots, [X_{i_k}, X_i]]]].$$

Denote the error terms in the sum by $Y_i=\partial_{x_i} G(m,x)-X_i$. Now assume $x_i \leq \delta$ for any $i$. The metric tensor can now be expressed in the given co-ordinates as follows.

\begin{eqnarray*}
\left\langle \partial_{x_i} G(m,x), \partial_{x_j} G(m,x) \right\rangle &=& \left\langle X_i + Y_i, X_j + Y_j \right\rangle \\
                                          &=& \delta_{ij} + \left\langle Y_i, Y_j \right\rangle + \left\langle X_i, Y_j \right\rangle + \left\langle Y_i, X_j \right\rangle
\end{eqnarray*}

Since $||Ad|| \leq C'(n)$, we get $||Y_i|| \leq M(n)\delta^k C'(n)^k$. Here, $\delta_{ij}$ is the Dirac delta. Choose $\delta(n)<\frac{1}{10^3 M(n)C'(n)}$ and apply Cauchy's inequality to give

\begin{equation}\label{estimatemetr}
1 - \frac{3}{8} \leq \left\langle \partial_{x_i} G(m,x), \partial_{x_i} G(m,x) \right\rangle \leq 1 + \frac{3}{8},
\end{equation}
whenever $x_i \leq 10\delta(n)$. Let $C(\delta(n))=\{(m,x) | x_i \leq 10\delta(n)\} \subset M^k \times \R^{n-k}$. From its definition and \eqref{estimatemetr},  we can see that $G(C(\delta(n))$ contains the $4\delta(n)$-tubular neighborhood $N_{4\delta(n)}(M^k)$ of $M^k$. Also, when restricting $G^{-1}$ onto $B_e(2\delta(n)) \subset N$ we get a  $\frac{1+\frac{3}{8}}{1-\frac{3}{8}}$-bi-Lipschitz into $M^k \times \R^{n-k}$. Thus, the inverse $G^{-1}$ is our desired map. We remark, that the distortion can be made arbitrarily small by choosing a smaller $\delta(n)$.

\QED

\bibliographystyle{amsplain}\bibliography{geometric}

\providecommand{\bysame}{\leavevmode\hbox to3em{\hrulefill}\thinspace}
\providecommand{\MR}{\relax\ifhmode\unskip\space\fi MR }
\providecommand{\MRhref}[2]{%
  \href{http://www.ams.org/mathscinet-getitem?mr=#1}{#2}
}
\providecommand{\href}[2]{#2}
\begin{thebibliography}{10}

\bibitem{alexbook}
S.~Alexander, V.~Kapovitch, and A.~Petrunin, \emph{Alexandrov geometry}, (In
  preparation).

\bibitem{andoni2015snowflake}
A.~Andoni, A.~Naor, and O.~Neiman, \emph{Snowflake {Universality} of
  {Wasserstein} spaces}, arXiv preprint arXiv:1509.08677 (2015).

\bibitem{assouad}
P.~Assouad, \emph{Plongements {lipschitziens} dans $\mathbb{R}^n$}, Bulletin de
  la Soci{\'e}t{\'e} Math{\'e}matique de France \textbf{111} (1983), 429--448.

\bibitem{auslander1960bieberbach}
L.~Auslander, \emph{Bieberbach's {Theorems} on {Space} {Groups} and {Discrete
  Uniform Subgroups of {Lie} Groups}}, Annals of Mathematics (1960), 579--590.

\bibitem{bonk2003bi}
M.~Bonk and U.~Lang, \emph{Bi-{Lipschitz} parameterization of surfaces},
  Mathematische Annalen \textbf{327} (2003), no.~1, 135--169.

\bibitem{bourgain85}
J.~Bourgain, \emph{On {lipschitz} embedding of finite metric spaces in
  {Hilbert} space}, Israel Journal of Mathematics \textbf{52} (1985), no.~1-2,
  46--52.

\bibitem{BuragoGeom}
D.~Burago, Y.~Burago, and S.~Ivanov, \emph{A course in metric geometry},
  Graduate Studies in Mathematics. AMS, vol.~33, 2001.

\bibitem{BGP}
Y.~D. Burago, M.~L. Gromov, and G.~Ya. Perel'man, \emph{Alexandrov spaces with
  curvature bounded below}, Uspekhi Mat. Nauk. (1992), no.~2(284), 3--51,
  Transl. Russian Math. Surveys, 1992, vol 47(2), 1--58.

\bibitem{bibergromov}
P.~Buser and H.~Karcher, \emph{The {Bieberbach} case in {Gromov}'s almost flat
  manifold theorem}, Global Differential Geometry and Global Analysis
  (Dirk~Ferus et~al., ed.), vol. 838, Springer, 1981, pp.~82--93.

\bibitem{buser1981gromov}
\bysame, \emph{Gromov's almost flat manifolds}, Soci{\'e}t{\'e}
  math{\'e}matique de France, 1981.

\bibitem{ChDiff99}
J.~Cheeger, \emph{Differentiability of {Lipschitz} {Functions} on {Metric
  Measure Spaces}}, Geometric \& Functional Analysis (GAFA) \textbf{9} (1999),
  no.~3, 428--517.

\bibitem{chgrfu}
J.~Cheeger, K.~Fukaya, and M.~L. Gromov, \emph{Nilpotent structures and
  invariant metrics on collapsed manifolds}, Journal of the American
  Mathematical Society \textbf{5} (1992), no.~2, 327--372.

\bibitem{collins2007jordan}
M.~J. Collins, \emph{On {Jordan's} theorem for complex linear groups}, Journal
  of Group Theory \textbf{10} (2007), no.~4, 411--423.

\bibitem{corwin2004representations}
L.~Corwin and F.~P. Greenleaf, \emph{{Representations of Nilpotent {Lie} Groups
  and Their Applications: {Basic} Theory and Examples}}, Cambridge university
  press, 2004.

\bibitem{ding2011restriction}
Y.~Ding, \emph{A restriction for singularities on collapsing orbifolds}, ISRN
  Geometry (2011).

\bibitem{ding2011restrictionarxiv}
\bysame, \emph{A restriction for singularities on collapsing orbifolds}, Arxiv
  preprint arXiv:1101.4444 (2011).

\bibitem{fukaya1988}
K.~Fukaya, \emph{A boundary of the set of the {Riemannian} manifolds with
  bounded curvatures and diameters}, Journal of Differential Geometry
  \textbf{28} (1988), no.~1, 1--21.

\bibitem{fukaya1989}
\bysame, \emph{Collapsing {Riemannian} manifolds to ones with lower dimension.
  {II}}, J. Math. Soc. Japan \textbf{41} (1989), no.~2, 333--356.

\bibitem{fukayaICM}
\bysame, \emph{Collapsing {Riemannian} manifolds and its applications},
  Proceedings of the International Congress of Mathematicians, vol.~1, 1990,
  pp.~491--500.

\bibitem{fukaya90vect}
\bysame, \emph{Hausdorff convergence of {Riemannian} manifolds and its
  applications.}, Recent topics in differential and analytic geometry
  (T.~Ochiai, ed.), Adv. Stud. Pure Math., vol.~18, Academic Press, 1990,
  pp.~143--238.

\bibitem{ghanaat}
P.~Ghanaat, M.~Min-Oo, and E.~A. Ruh, \emph{Local structure of {Riemannian}
  manifolds}, Indiana University Mathematics Journal \textbf{39} (1990), no.~4,
  1305--1312.

\bibitem{greene1988lipschitz}
R.~Greene and H.-H. Wu, \emph{Lipschitz convergence of {Riemannian} manifolds},
  Pacific journal of mathematics \textbf{131} (1988), no.~1, 119--141.

\bibitem{gromov78almost}
M.~L. Gromov, \emph{Almost flat manifolds}, Journal of Differential Geometry
  \textbf{13} (1978), no.~2, 231--241.

\bibitem{Gbook}
\bysame, \emph{Metric structures for {Riemannian} and non-{Riemannian} spaces},
  Modern Birkh\"auser Classics, 1999.

\bibitem{grove1973conjugatec}
K.~Grove and H.~Karcher, \emph{How to conjugate {$C^1$}-close group actions},
  Mathematische Zeitschrift \textbf{132} (1973), no.~1, 11--20.

\bibitem{haviv2013euclidean}
I.~Haviv and O.~Regev, \emph{The {Euclidean} distortion of flat tori}, Journal
  of Topology and Analysis \textbf{5} (2013), no.~02, 205--223.

\bibitem{Heinonen2000}
J.~Heinonen, \emph{{Lectures on Analysis on Metric Spaces}}, Springer, 2000.

\bibitem{heinonen2005lectures}
\bysame, \emph{Lectures on {Lipschitz} analysis}, Rep. Univ. Jyv\"askyl\"a
  Dept. Math. Stat. 100, 2005.

\bibitem{heinonen2011flat}
J.~Heinonen and S.~Keith, \emph{Flat forms, bi-{Lipschitz} parametrizations,
  and smoothability of manifolds}, Publications math{\'e}matiques de l'IH{\'E}S
  \textbf{113} (2011), no.~1, 1--37.

\bibitem{indyk2004}
P.~Indyk and J.~Matou\v{s}ek, \emph{{Chapter 8, Low-Distortion Embeddings of
  Finite Metric Spaces}}, Handbook of Discrete and Computational Geometry (ed.
  C. D. Toth, J. O'Rourke, J. E. Goodman), vol.~37, 2004, pp.~177--197.

\bibitem{kapovitch2007perelman}
V.~Kapovitch, \emph{Perelman's {Stability} theorem}, Surveys in differential
  geometry: {Metric} and {Comparison} geometry, vol.~XI, Intl Press, 2007.

\bibitem{center}
H.~Karcher, \emph{Riemannian center of mass and mollifier smoothing},
  Communications on pure and applied mathematics \textbf{30} (1977), no.~5,
  509--541.

\bibitem{naorfourier}
S.~Khot and A.~Naor, \emph{Nonembeddability theorems via {Fourier} analysis},
  Mathematische Annalen \textbf{334} (2006), no.~4, 821--852.

\bibitem{Kleiner2014}
B.~Kleiner and J.~Lott, \emph{{Local Collapsing, Orbifolds, and
  Geometrization}}, Ast{\'e}risque, vol. 365, Soci{\'e}t{\'e} math{\'e}matique
  de France, 2014.

\bibitem{Laakso2000}
T.~J. Laakso, \emph{Ahlfors $q$-regular spaces with arbitrary $q>1$ admitting
  weak {Poincar\`{e}} inequality}, Geomet. func. anal \textbf{10} (2000),
  111--123.

\bibitem{langplaut}
U.~Lang and C.~Plaut, \emph{{Bilipschitz Embeddings of Metric Spaces into Space
  Forms}}, Geometriae Dedicata \textbf{87} (2001), no.~1-3, 285--307.

\bibitem{kirsz}
U.~Lang and V.~Schroeder, \emph{Kirszbraun's theorem and metric spaces of
  bounded curvature}, Geometric \& Functional Analysis (GAFA) \textbf{7}
  (1997), no.~3, 535--560.

\bibitem{Luosto1996}
K.~Luosto, \emph{Ultrametric spaces bi-{Lipschitz} embeddable in
  $\mathbb{R}^n$}, Fundamenta Mathematicae \textbf{150} (1996), no.~1, 25--42.

\bibitem{Luukkainen1994}
J.~Luukkainen and H.~Movahedi-Lankarani, \emph{Minimal bi-{Lipschitz} embedding
  dimension of ultrametric spaces}, Fundamenta Mathematicae \textbf{144}
  (1994), no.~2, 181--193.

\bibitem{malcevclass}
A.~Malcev, \emph{On a class of homogeneous spaces}, Amer. Math. Soc.
  Translations \textbf{9} (1962), no.~1, 276--304 (Originally in Russian,
  {Izvestia} {Ak.} {Nauk} {SSSR}, {Ser.} {Math.}, 13 (39)).

\bibitem{matsek90}
J.~Matou{\v{s}}ek, \emph{Bi-{Lipschitz} embeddings into low-dimensional
  {Euclidean} spaces}, Commentationes Mathematicae Universitatis Carolinae
  \textbf{31} (1990), no.~3, 589--600.

\bibitem{leftinv}
J.~Milnor, \emph{Curvatures of left invariant metrics on {Lie} groups},
  Advances in Mathematics \textbf{21} (1976), no.~3, 293 -- 329.

\bibitem{nagelstein}
A.~Nagel, E.~M. Stein, and S.~Wainger, \emph{{Balls and metrics defined by
  vector fields I}: {Basic} properties}, Acta Mathematica \textbf{155} (1985),
  no.~1, 103--147.

\bibitem{naorintro}
A.~Naor, \emph{An introduction to the {Ribe} program}, Japanese Journal of
  Mathematics \textbf{7} (2012), no.~2, 167--233.

\bibitem{naorassouad}
A.~Naor and O.~Neiman, \emph{Assouad's theorem with dimension independent of
  the snowflaking}, Revista Matematica Iberoamericana \textbf{28} (2012),
  no.~4, 1123--1142.

\bibitem{naormarkov}
A.~Naor, Y.~Peres, O.~Schramm, and S.~Sheffield, \emph{Markov chains in smooth
  {Banach} spaces and {Gromov}-hyperbolic metric spaces}, Duke Mathematical
  Journal \textbf{134} (2006), no.~1, 165--197.

\bibitem{nikolaev1991}
I.~G. Nikolaev, \emph{Bounded curvature closure of the set of compact
  {Riemannian} manifolds}, Bulletin of the American Mathematical Society
  \textbf{24} (1991), no.~1, 171--177.

\bibitem{pansu89}
P.~Pansu, \emph{M{\'e}triques de {Carnot}-{Carath{\'e}odory} et
  quasiisom{\'e}tries des espaces sym{\'e}triques de rang un}, Annals of
  Mathematics (1989), 1--60.

\bibitem{ratcliffehyper}
J.~Ratcliffe, \emph{Foundations of hyperbolic manifolds}, vol. 149, Springer,
  2006.

\bibitem{romneymatt}
M.~Romney, \emph{Conformal {Grushin} spaces}, Conformal Geometry and Dynamics
  of the American Mathematical Society \textbf{20} (2016), no.~6, 97--115.

\bibitem{rong1996fundamental}
X.~Rong, \emph{On the fundamental groups of manifolds of positive sectional
  curvature}, Annals of mathematics (1996), 397--411.

\bibitem{ruh1982}
E.~A. Ruh, \emph{Almost flat manifolds}, Journal of Differential Geometry
  \textbf{17} (1982), no.~1, 1--14.

\bibitem{semmes1993bilipschitz}
S.~Semmes, \emph{Bilipschitz mappings and strong {$A_\infty$} weights}, Ann.
  Acad. Sci. Fenn. Ser. AI Math. \textbf{18} (1993), no.~2, 211--248.

\bibitem{semmes1996nonexistence}
\bysame, \emph{On the nonexistence of bilipschitz parametrizations and
  geometric problems about {$A_\infty$}-weights.}, Revista Matem{\'a}tica
  Iberoamericana \textbf{12} (1996), no.~2, 337--410.

\bibitem{seo2011char}
J.~Seo, \emph{A characterization of bi-{Lipschitz} embeddable metric spaces in
  terms of local bi-{Lipschitz} embeddability}, Math. Res. Lett. \textbf{18}
  (2011), no.~6, 1179--1202.

\bibitem{Sormani2012}
C.~Sormani, \emph{How {Riemannian} manifolds converge}, Metric and Differential
  Geometry, Springer, 2012, pp.~91--117.

\bibitem{tao2014hilbert}
T.~Tao, \emph{{Hilbert's} fifth problem and related topics}, vol. 153, American
  Mathematical Soc., 2014.

\bibitem{thurston2002geometry}
W.~P. Thurston, \emph{Geometry and topology of 3-manifolds}, Notes, Princeton
  University (1977--1978).

\bibitem{thurston1997three}
W.~P. Thurston and S.~Levy, \emph{Three-dimensional geometry and topology},
  vol.~1, Princeton University Press, 1997.

\bibitem{toro94}
T.~Toro, \emph{Surfaces with generalized second fundamental form in {$L^{2}$}
  are {Lipschitz} manifolds}, Journal of Differential Geometry \textbf{39}
  (1994), no.~1, 65--101.

\bibitem{toro95}
\bysame, \emph{Geometric conditions and existence of bi-{Lipschitz}
  parameterizations}, Duke Mathematical Journal \textbf{77} (1995), no.~1,
  193--227.

\bibitem{varopoulos2008analysis}
N.~T. Varopoulos, L.~Saloff-Coste, and T.~Coulhon, \emph{Analysis and geometry
  on groups}, vol. 100, Cambridge university press, 2008.

\end{thebibliography}

\end{document}